\documentclass[reqno,11pt]{amsart}
\usepackage[T1]{fontenc}
\usepackage[utf8]{inputenc}
\usepackage[english]{babel}
\usepackage{lmodern,microtype,csquotes}
\usepackage{amsmath,amssymb,amsfonts,amsthm,mathtools,mathrsfs}
\usepackage{enumitem,graphicx,dsfont}
\usepackage{cases}
\usepackage[a4paper,margin=2.4cm]{geometry}
\usepackage[colorlinks=true,linkcolor=black,citecolor=black,urlcolor=blue]{hyperref}
\usepackage{aliascnt}
\usepackage[capitalise,nameinlink]{cleveref}
\usepackage[backend=biber,sorting=nyt,maxbibnames=12,doi=true,url=false,isbn=false]{biblatex}
\numberwithin{equation}{section}
\allowdisplaybreaks[2]
\newtheorem{theorem}{Theorem}[section]
\newaliascnt{lemma}{theorem}
\newtheorem{lemma}[lemma]{Lemma}
\aliascntresetthe{lemma}
\newaliascnt{proposition}{theorem}
\newtheorem{proposition}[proposition]{Proposition}
\aliascntresetthe{proposition}
\newaliascnt{corollary}{theorem}
\newtheorem{corollary}[corollary]{Corollary}
\aliascntresetthe{corollary}
\theoremstyle{definition}
\newaliascnt{definition}{theorem}
\newtheorem{definition}[definition]{Definition}
\aliascntresetthe{definition}
\newaliascnt{example}{theorem}

\theoremstyle{remark}
\newaliascnt{remark}{theorem}
\newtheorem{remark}[remark]{Remark}
\aliascntresetthe{remark}
\newcommand{\numberthis}{\addtocounter{equation}{1}\tag{\theequation}}
\DeclareMathOperator{\Ric}{Ric}
\DeclareMathOperator{\Hess}{Hess}
\DeclareMathOperator{\Sec}{Sec}
\DeclareMathOperator{\Supp}{Supp}
\DeclareMathOperator{\Var}{Var}
\DeclareMathOperator{\CD}{CD}
\DeclareMathOperator{\MCP}{MCP}
\newcommand{\Vol}{\mathrm{Vol}}
\newcommand{\m}{\mathfrak m}
\newcommand{\N}{\mathbb N}
\newcommand{\R}{\mathbb R}
\newcommand{\E}{\mathbb E}
\renewcommand{\P}{\mathbb P}
\newcommand{\abs}[1]{\left\lvert#1\right\rvert}
\newcommand{\dd}{\,\mathrm d}

\begin{document}
\raggedbottom
\title[Quantitative Wasserstein bounds for empirical measures on manifolds]{Quantitative $p$-Wasserstein bounds for empirical measures on manifolds: Expectation and concentration}
\author{Jonathan Junn\'e}
\address{Univ Rennes, INRIA (Mingus), Institut de Recherche Math\'ematiques de Rennes,
CNRS UMR 6625, Campus Beaulieu, F-35042 Rennes Cedex, France}
\email{jonathan.junne@inria.fr}
\date{\today}

\begin{abstract}
    We derive quantitative bounds in expectation and
    concentration for the $p$-Wasserstein distance between
    the empirical measure of independent samples and their
    common law on metric measure spaces satisfying the
    measure contraction property $\MCP(K,n)$, with $K\le0$.
    Our results extend the expectation estimates of
    Boissard--Le Gouic (Theorem 1.1 in Ann. Inst. Henri
    Poincar\'e Probab. Stat. 50 (2014), no. 2, 539--563)
    and the concentration estimates of
    Bolley--Guillin--Villani (Theorem 2.1 in Probab.
    Theory Relat. Fields 137 (2007), 541--593).
    In spaces with exponential volume growth, suitable
    exponential moment assumptions allow us to recover
    the Euclidean rate $N^{-1/n}$ in expectation when
    $n>2p$. As an application, we obtain uniform-in-time concentration estimates for the 
    empirical measure of McKean--Vlasov interacting
    particle systems on complete Riemannian manifolds
    with Ricci curvature bounded below by combining comparison theory and parallel coupling to a noninteracting proxy system.
\end{abstract}


\maketitle
\section{Introduction}

The convergence of empirical measures built from independent and
identically distributed samples to their common law in
$p$-Wasserstein distance is well understood on compact spaces, Euclidean spaces,
and more generally spaces with polynomial volume growth, under suitable moment
assumptions; see
\cite{boissard_mean_2014,bolley_quantitative_2007,dedecker_concentration_2026,fournier_rate_2015,weed_sharp_2019}
and the references therein.
On noncompact manifolds, however, negative curvature may lead
to exponential volume growth. This distinction already appears in the quantization of measures
on the hyperbolic plane \cite{iacobelli_asymptotic_2016}, where the Euclidean rate fails in general
without additional exponential moment assumptions compensating
for the growth of volume.

The purpose of this paper is to quantify how lower curvature
bounds affect these convergence rates, both in expectation
and in concentration, following the two complementary approaches
of Boissard--Le Gouic \cite{boissard_mean_2014} and
Bolley--Guillin--Villani \cite{bolley_quantitative_2007}.
Finally, combining these estimates for independent samples
with comparison theory and parallel coupling on manifolds,
we obtain concentration estimates for McKean--Vlasov mean-field interacting particle systems approximating aggregation--diffusion on manifolds with Ricci curvature bounded below,
uniformly over finite time intervals.

\subsection{Related literature}

As explained above, quantitative convergence of empirical measures has
been extensively studied in Euclidean spaces. We refer in particular
to the seminal works of Fournier--Guillin
\cite{fournier_convergence_2023,fournier_rate_2015}
for estimates under moment assumptions, and to
Bolley--Guillin--Villani \cite{bolley_quantitative_2007}
for concentration inequalities.
On general metric spaces, Boissard--Le Gouic
\cite{boissard_simple_2011,boissard_mean_2014}
obtained nonasymptotic expectation bounds in terms of covering
numbers, while Weed--Bach \cite{weed_sharp_2019}
studied sharp convergence rates on compact metric spaces.
More recently, Dedecker--Fischer--Michel
\cite{dedecker_concentration_2026} extended the concentration
estimates of Fournier--Guillin to general Polish spaces,
with the covering dimension of the support replacing the
dimension of the ambient Euclidean space. These results do not cover spaces with exponential volume growth.

On Riemannian manifolds, Wang and collaborators have obtained
Wasserstein convergence rates and limit theorems for empirical
and occupation measures of diffusion processes, using semigroup
and spectral methods
\cite{wang_precise_2021,wang_wasserstein_2022,
wang_wasserstein_2021}.
For McKean--Vlasov interacting particle systems,
Arnaudon--Del Moral \cite{arnaudon_delmoral_2019}
studied nonlinear diffusions, their manifold extensions,
and quantitative propagation of chaos.
Recently, Chae--Cho--Jang--Yang
\cite{chae_cho_jang_yang_2026}
established global well-posedness, uniform moment estimates,
and uniform-in-time propagation of chaos in a modified
Wasserstein distance, using reflection coupling.
These results complement ours, as we focus on concentration estimates for the full empirical measure in $p$-Wasserstein distance, using a different approach based on independent-sample concentration and parallel coupling.

Related results on aggregation and diffusion on the sphere include
the analysis of stationary McKean--Vlasov states and bifurcations
by Shalova--Schlichting \cite{shalova_schlichting_2026},
the study of aggregation--diffusion energy minimizers
by Fetecau--Park--Vaidya \cite{fetecau_park_vaidya_2026},
and the nonlinear diffusion limit of nonlocal aggregation
interactions by Peletier--Shalova
\cite{peletier_shalova_2025}.
These works concern related models, although the global
convexity assumptions used in our particle-system application
cannot hold on a compact manifold.

\subsection{Organization}

We introduce the metric and geometric setting and state the main results in \cref{sec:def and results}. The expectation theorem is proved by applying the covering estimate on annuli, followed by concentration around the expectation and the quantitative Sanov argument in \cref{sec:Proofs for independent random variables}. In \cref{sec:An interacting particle system}, we follow the McKean--Vlasov construction through a noninteracting proxy, prove the moment and increment bounds with the necessary localization, and then use parallel coupling to obtain concentration and propagation of chaos. The well-posedness, radial cutoff, and covering arguments are completed in the appendices.

\section{Definitions and main results}\label{sec:def and results}
\subsection{Wasserstein distance and relative entropy}

We first recall the notions of Wasserstein distance, relative entropy,
and transport-entropy inequalities; we refer to Villani
\cite{villani_optimal_2009} for a detailed account. Let \((E,d)\) be a
Polish metric space. For \(p\ge 1\), we denote by \(\mathcal P_p(E)\) the
set of probability measures on \(E\) with finite \(p\)-th moment.

\begin{definition}[Wasserstein distance]
    Let \(\mu,\nu\in\mathcal P_p(E)\). The \emph{\(p\)-Wasserstein distance}
    between \(\mu\) and \(\nu\) is defined by
    \[
        W_p(\mu,\nu)
        \coloneqq
        \inf_{\pi\in\Pi(\mu,\nu)}
        \left(
            \int_{E\times E} d^p(x,y)\dd\pi(x,y)
        \right)^{1/p},
    \]
    where \(\Pi(\mu,\nu)\) is the set of probability measures on
    \(E\times E\) with first marginal \(\mu\) and second marginal
    \(\nu\).
\end{definition}

In this paper, \(W_p(\widehat\nu^N,\nu)\) measures how far the
empirical measure
\[
    \widehat\nu^N \coloneqq \frac1N\sum_{i=1}^N\delta_{X^i},
    \qquad X^1,\dots,X^N\sim \nu \text{ i.i.d.},
\]
is from its common law \(\nu\).  Large deviations for empirical measures
are governed by relative entropy. 

\begin{definition}[Relative entropy]
    If \(\mu,\nu\in\mathcal P(E)\), the
\emph{relative entropy} of \(\mu\) with respect to \(\nu\) is
\[
    H(\mu\mid\nu)
    \coloneqq
    \begin{cases}
        \displaystyle
        \int_E \log\left(\frac{d\mu}{d\nu}\right)\dd\mu,
        & \text{if } \mu\ll\nu;\\[1.2em]
        +\infty,
        & \text{otherwise.}
    \end{cases}
\]
\end{definition}
Sanov's theorem says, roughly speaking, that deviations of
\(\widehat\nu^N\) from \(\nu\) have probability of order
\[
    \P\left(\operatorname{dist}(\widehat\nu^N,\nu)>\varepsilon\right)
    \approx
    \exp\!\left(
        -N
        \inf_{\operatorname{dist}(\mu,\nu)>\varepsilon}
        H(\mu\mid\nu)
    \right),
\]
for distances compatible with weak convergence. In particular, to obtain
concentration estimates in Wasserstein distance, one needs lower bounds
on \(H(\mu\mid\nu)\) in terms of \(W_p(\mu,\nu)\). This is precisely the
role of transportation-cost inequalities.

\begin{definition}[Talagrand inequality]
    Let \(p\ge 1\) and \(\lambda>0\). We say that
    \(\nu\in\mathcal P_p(E)\) satisfies the \emph{transport-cost inequality},
    or \emph{Talagrand inequality}, \(T_p(\lambda)\) if
    \begin{equation}\label{eq:Tp lambda}
        W_p(\mu,\nu)
        \le
        \sqrt{\frac{2}{\lambda}H(\mu\mid\nu)}
    \end{equation}
    for every probability measure \(\mu\) on \(E\).
\end{definition}

Equivalently, \(T_p(\lambda)\) gives the entropy lower bound
\[
    H(\mu\mid\nu)
    \ge
    \frac{\lambda}{2}W_p^2(\mu,\nu).
\]
Inserted into the Sanov rate function, this suggests the Gaussian
concentration scale
\begin{equation}\label{eq:desired concentration estimate}
    \P\left(W_p(\widehat\nu^N,\nu)>\varepsilon\right)
    \lesssim
    \exp\!\left(-\frac{\lambda}{2}N\varepsilon^2\right),
\end{equation}
up to the metric-entropy corrections needed to make this estimate
nonasymptotic.

The condition \(T_p(\lambda)\) becomes stronger as \(p\) increases, since
\(W_p\ge W_q\) whenever \(p\ge q\) on probability measures with finite
moments. In particular, \(T_2\) implies \(T_1\); see \cite{villani_optimal_2009} for the following properties. On a Riemannian manifold, if the probability measure $\nu = e^{-U}d\Vol$ with $C^2$ potential is such that $\Ric + \Hess U \ge K$ for some $K > 0$, then $\nu$ satisfies $T_2(K)$. \(T_1\) is equivalent to the existence of
a square-exponential moment; there exist \(o\in E\) and \(\alpha>0\) such
that
\[
    \int_E e^{\alpha d^2(o,x)}\dd\nu(x)<+\infty .
\]

\subsection{Geometric setting}
On a metric space with finite diameter, Boissard--Le Gouic~\cite{boissard_mean_2014} obtained a quantitative bound on the expected value of the Wasserstein distance between the empirical and the sampling measure on bounded metric spaces in terms of covering numbers.
\begin{definition}[Covering and packing numbers]
    Let $(E, d)$ be a metric space. 
    The \emph{covering number} $N(S, \delta)$ of a set $S \subset E$ denotes the minimal number of balls needed to cover $S$ by balls of radius $\delta > 0$ with respect to the metric $d$. The \emph{packing number} $P(S,\delta)$ is the maximum cardinality of any $\delta$-packing, namely, any collection of points in $S$ of distance at least $\delta$ from each other.
\end{definition}
\begin{theorem}[{\cite[Theorem 1.1]{boissard_mean_2014}}]\label{thm:BLG}
    Let $(E, d)$ be a metric space with finite diameter $D$ and $\nu$ be a probability measure on $E$. Choose $0<\varepsilon <D/4$ such that $N(E, \varepsilon) < +\infty$. 
    Then
    \begin{equation}\label{eq:estimate BG}
        \begin{aligned}
            \E W_p(\widehat\nu^N,\nu)
            \le C\biggl(\varepsilon
                +N^{-1/(2p)}\int_\varepsilon^{D/4}
                    N(E,\delta)^{1/(2p)}\dd\delta
            \biggr),
        \end{aligned}
    \end{equation}
    with $C \le 64/3$.
\end{theorem}
For an unbounded metric space $(E,d)$, the covering number $N(E, \delta)$ is infinite for every $\delta > 0$.
In that setting, the estimate \eqref{eq:estimate BG} can only be applied to the probability measure truncated on some possibly large ball. 
Most of the samples lie in a large compact under some moment assumptions. 
A typical way to estimate the covering number of a large ball is to control its volume growth.
Such a control holds on metric measure spaces that satisfy the measure contraction property due to Ohta \cite{ohta_measure_2007}.
\begin{definition}[Measure contraction property]
    Let $(E, d, \m)$ be a metric measure space. Let $e_t(\gamma) \coloneqq \gamma(t)$ be the \emph{evaluation map} on the set of minimal geodesics $\Gamma \coloneqq \{ \gamma : [0,1] \to E \: : \: \gamma \text{ is minimal}\}$, and denote the length by $\ell(\gamma) \coloneqq d(\gamma(0), \gamma(1))$. 
    A \emph{dynamical transference plan} $\Pi$ is a Borel probability measure on $\Gamma$ and $\m_t \coloneqq (e_t)_\#\Pi$ is the associated \emph{displacement interpolation}. 
    For $K \in \R$ and $n > 1$, the space $(E, d, \m)$ is said to satisfy a \emph{measure contraction property} $\MCP(K, n)$ if for every $x \in E$ and measurable set $A \subseteq E$, there exists a displacement interpolation $\m_t$ associated with a transference plan $\Pi$ satisfying the following:
    \begin{enumerate}[label=(\roman*)]
        \item We have $\m_0 = \delta_x$ and $\m_1 = \m[A]^{-1} \m|_A$ as measures;
        \item For every $t \in [0, 1]$,
        \begin{equation*}
            d\m \ge (e_t)_\# \left(t \frac{s_{K,n}\left(t \ell(\gamma)\right)}{s_{K,n}\left(\ell(\gamma)\right)}\, \m[A] \: d\Pi(\gamma)\right)
        \end{equation*}
        holds as measures with the convention $0/0 \coloneqq 1$, where
        \begin{equation*}
            s_{K,n}(t) \coloneqq \begin{cases}
                \sin^{n-1}\left(\sqrt{\tfrac{K}{n-1}}t\right) & \text{if } K > 0; \\
                t^{n-1} & \text{if } K = 0; \\
                \sinh^{n-1}\left(\sqrt{\tfrac{|K|}{n-1}}t\right) & \text{if } K < 0.
            \end{cases}
        \end{equation*}
    \end{enumerate}
\end{definition}

The Bonnet--Myers theorem holds for metric measure spaces satisfying the measure contraction property \cite[Theorem 4.3]{ohta_measure_2007}; if $(E, d, \m)$ satisfies $\MCP(K,n)$ for some $K > 0$ and $n > 1$, then
\begin{equation*}
    \operatorname{diam}(E) \le \pi\sqrt{\frac{n-1}{K}}.
\end{equation*}
For the noncompact covering estimates, we will restrict ourselves to $K\le0$.
\begin{remark}[Curvature-dimension and measure contraction]
    The measure contraction property is a weaker notion than the synthetic notion of curvature-dimension condition $\CD(K,n)$; we refer again to the monograph of Villani \cite{villani_optimal_2009} for an account on the subject. More precisely, if the weighted $m$-dimensional Riemannian manifold $(M,g,\m)$ with reference measure $\m = e^{-U}d\Vol$ satisfies $\CD(K,n)$ with $n \ge m \ge 2$,
    \[
    \Ric_{U,n} \coloneqq \Ric + \Hess U - \frac{\nabla U \otimes \nabla U}{n-m} \ge K
    \]
    with the convention that $U$ is constant when $n=m$, so that $\Ric_{U,m}=\Ric$, then it satisfies $\MCP(K,n)$ \cite[Theorem 3.2]{ohta_measure_2007}
\end{remark}
Metric measure spaces satisfying the measure contraction property also satisfy the Bishop--Gromov volume comparison theorem.
\begin{theorem}[{Bishop--Gromov \cite[Theorem 5.1]{ohta_measure_2007}}]\label{thm:Bishop Gromov}
    Let $(E, d, \m)$ be a metric measure space satisfying $\MCP(K,n)$ for some $K \in \R$ and $n > 1$. 
    Then for any $x \in E$, the function
    \begin{equation*}
        r \mapsto \frac{\m\left[B(x, r)\right]}{V_{K,n}(r)}
    \end{equation*}
    is monotone non-increasing; $V_{K,n}(r) \coloneqq \int_0^r s_{K,n} (l) \: dl$ denotes the volume of the ball in the reference space. 
\end{theorem}
As a consequence, the following version of Gromov's packing lemma, which estimates the growth of the covering number of balls under the measure contraction property, holds. The proof is deferred to \cref{sec:appendix covering}.
\begin{lemma}[Covering number estimate]\label{lem:packing}
    Let $(E, d, \m)$ be a metric measure space satisfying $\MCP(K,n)$ for some $K \le 0$ and $n > 1$.
    Then for any $x \in E$, $0 < \delta \le R/2$ with $ R \le \operatorname{diam}(E)$,
    \begin{equation}\label{eq:covering estimate uniform K}
        N(B(x, R), \delta) \le \left(\frac{5R}{\delta}\right)^n \exp\left(2\sqrt{|K|(n-1)}R\right)
    \end{equation}
\end{lemma}
\subsection{Upper bounds for the expected Wasserstein distance} 
A reason to pay attention to the exponential integrability condition is the relation with \emph{quantization of measures} on noncompact Riemannian manifolds investigated by Iacobelli \cite{iacobelli_asymptotic_2016}.
Quantization of measures seeks the best possible approximation of the probability measure $\nu$ among all probability measures $\nu^N$ supported in at most $N$ points; the quantization error is given in $p$-Wasserstein distance as
\begin{equation*}
    e_{N,p}(\nu) := \inf \left\{ W_p(\nu^N, \nu); \: \#\left(\Supp\, \nu^N\right) \le N\right\}.
\end{equation*}
\begin{theorem}[{\cite[Corollary 1.6] {iacobelli_asymptotic_2016}}]\label{thm:iacobelli}
    Let $M$ be a complete connected Riemannian manifold of dimension $m \ge 2$ with sectional curvature bounded below; i.e., $\Sec \ge K/(m-1)$ from some $K \le 0$. Let $p \ge 1$ and let $\nu$ be a probability measure on $M$ with absolutely continuous component $\rho$ with respect to the Riemannian volume such that
    \begin{equation*}
        \int_M d^{p+\delta}(o,x) + e^{p\sqrt{\frac{|K|}{m-1}}\, d(o, x)} \dd\nu(x) < +\infty
    \end{equation*}
    for some $o \in M$ and $\delta > 0$. Then
    \begin{equation}\label{eq:Iacobelli quantisation}
        \lim_{N \to \infty} N^{1/m} e_{N,p}(\nu) = Q_p^{1/p}([0,1]^m)\left(\int_M \rho(x)^{m/(m+p)} \dd\Vol(x)\right)^{(m+p)/(mp)},
    \end{equation}
    where $Q_p([0,1]^m)$ denotes the quantization coefficient.
\end{theorem}
\cite{iacobelli_asymptotic_2016} provides an example of probability measure on the hyperbolic plane where all the polynomial moments are finite but no exponential moment is, and in particular the condition $\int_M e^{p\sqrt{|K|/(m-1)}\, d(o, x)} d\nu(x)$ coming from the sectional curvature is not satisfied. $N^{1/m}e_{N,p}(\nu)$ diverges as $N \to +\infty$ for this measure. The lower bound on the sectional curvature was then replaced by general covering growth estimates in \cite{aydin_asymptotic_2025}. For a sub-Riemannian counterpart, see the recent work \cite{dimarco_quantization_2026} on Carnot groups.

This suggests that we can obtain a quantitative estimate \[\E\left[W_p(\widehat\nu^N, \nu)\right] \lesssim N^{-1/n}\]
under $\MCP(K,n)$ whenever the measure satisfies some sufficiently high exponential moment. This is exactly the content of the following theorem.
\begin{theorem}[Expectation bound under exponential tails]\label{thm:sharp E nu n n K neg}
    Let $(E,d,\m)$ be an unbounded metric measure space satisfying $\MCP(K,n)$ for some $K \le 0$ and $n>2p \ge 2$, and let $\nu$ be a probability measure on $E$ that satisfies the exponential moment condition
    \begin{equation*}
        \mathcal E_\alpha \coloneqq \int_E e^{\alpha d(o,x)}\dd\nu(x)<+\infty
    \end{equation*}
    for some $o\in E$ and
    \begin{equation}\label{eq:assumption on alpha}
        \alpha > \frac{2p}{n-p}\sqrt{|K|(n-1)}.
    \end{equation}
    Then there exists a constant $C$ depending only on $K,n,p,\alpha,$ and on $\mathcal E_\alpha$ such that, for every $N\ge 1$,
    \begin{equation}\label{eq:optimal estimate MCP E alpha}
        \E\big[W_p(\widehat\nu^N,\nu)\big]
        \le C N^{-1/n}.
    \end{equation}
\end{theorem}
\begin{remark}[Quantization and curvature dependence]
    Even though the condition \eqref{eq:assumption on alpha} on the exponential moment of \cref{thm:sharp E nu n n K neg} is stronger than the one of \cref{thm:iacobelli}, \eqref{eq:optimal estimate MCP E alpha} is nonasymptotic while \eqref{eq:Iacobelli quantisation} is asymptotic. Note that a lower bound $\Sec \ge K/(m-1)$ on the sectional curvature implies $\Ric\ge Kg$ for Riemannian manifolds of dimension $m \ge 2$; then for $m > 2p$ and $K < 0$,
    \begin{equation}\label{eq:ratio}
        \left(\frac{2p}{m-p}\sqrt{|K|(m-1)}\right) \div \left(p\sqrt{|K|/(m-1)}\right) = 2\frac{m-1}{m-p} \ge 2.
    \end{equation}

    The factor $2$ in the ratio \eqref{eq:ratio} comes from the exponential $e^{2\sqrt{|K|(n-1)}R}$ in the covering estimate \eqref{eq:covering estimate uniform K}, as it will become apparent in the proof. On homogeneous spaces, one can improve the covering estimate, replacing the $2$ in the exponential factor by $1$. Then in that case $\frac{m-1}{m-p} \ge 1 $, and for $p = 1$ and dimension $m \ge 3$, \eqref{eq:assumption on alpha} becomes $\alpha > \sqrt{|K|/(m-1)}$, which is only stronger than the condition of \cref{thm:iacobelli} for the inequality being strict.

    Finally, the curvature-dependent shell contribution to the
    constant $C$ in \eqref{eq:optimal estimate MCP E alpha}
    contains the factor
    \begin{equation*}
        \frac{e^{2\sqrt{|K|(n-1)}/n}}{(1-e^{-\beta_K})^2},
        \qquad
        \beta_K\coloneqq
        \alpha\left(\frac1p-\frac1n\right)
        -\frac{2\sqrt{|K|(n-1)}}{n}.
    \end{equation*}
    For fixed $\alpha$ and its moment bound, this choice of $C$
    depends continuously on $K$ as $K\uparrow0$.
    At $K=0$, condition \eqref{eq:assumption on alpha}
    becomes $\alpha>0$; weaker moment assumptions may also
    suffice, but a finite $p$-th moment alone does not guarantee
    an algebraic rate.
\end{remark}

\begin{corollary}[Concentration around the expectation under a transport-cost inequality]
\label{cor:squareexp-expectation}
\label{thm:centered-concentration}
    Under the assumptions of
    \cref{thm:sharp E nu n n K neg}, suppose, in addition, that $\nu$ satisfies
    $T_p(\lambda)$ for some $\lambda>0$.
    Then, for every $N\ge1$ and $\varepsilon>0$,
    \begin{equation}\label{eq:concentration-shifted}
        \P\left(W_p(\widehat\nu^N,\nu)>\varepsilon\right)
        \le
        \exp\left[
            -\frac{\lambda N}{2}
            \big(\varepsilon-CN^{-1/n}\big)_+^2
        \right],
    \end{equation}
    where $C$ is the constant in
    \eqref{eq:optimal estimate MCP E alpha}.
\end{corollary}

\begin{proof}
    Since $T_p(\lambda)$ implies $T_1(\lambda)$,
    the measure $\nu$ admits a finite square-exponential
    moment, and hence exponential moments of every order.
    In particular, the moment condition
    \eqref{eq:assumption on alpha} is satisfied for
    a suitable $\alpha>0$.
    Thus, by \cref{thm:sharp E nu n n K neg},
    \[
        \E W_p(\widehat\nu^N,\nu)\le CN^{-1/n}.
    \]

    On the other hand, by tensorization of the transport-cost
    inequality $T_p(\lambda)$, the concentration result of
    Gozlan--L\'eonard \cite[Theorem 12]{gozlan_leonard_2007}
    gives, for every $t\ge0$,
    \[
        \P\left(W_p(\widehat\nu^N,\nu)
            \ge \E W_p(\widehat\nu^N,\nu)+t\right)
        \le
        \exp\left(-\frac{\lambda N}{2}t^2\right).
    \]
    Combining these two estimates yields
    \eqref{eq:concentration-shifted}.
\end{proof}

\subsection{Quantitative Sanov's theorem}
The covering number of a metric space is also related to the covering number of probability measures supported on that space with respect to the Wasserstein distance, as the following proposition shows.
\begin{proposition}[Metric entropy of probability measures]\label{thm:Metric entropy}
    Let $(E, d)$ be a Polish metric space with finite diameter $D$. By \cite[Theorem A.1]{bolley_quantitative_2007}, for all $\delta \in (0, D)$, the space $\mathcal{P}(E)$ of probability measures on $E$ can be covered by $\mathcal{N}_p(E, \delta) \coloneqq N(\mathcal{P}(E), \delta)$ balls of radius $\delta$ in $W_p$ distance, with
    \begin{equation}\label{eq:estimate Wasserstein covering vs metric}
        \mathcal{N}_p(E, \delta) \le \left(\frac{8eD}{\delta}\right)^{pN\left(E, \frac{\delta}{2}\right)}.
    \end{equation}
\end{proposition}
Thanks to this inequality, Bolley--Guillin--Villani derived the following quantitative version of Sanov's theorem on bounded spaces.
\begin{proposition}[Quantitative Sanov bound on bounded spaces]
    As in \cite[Theorem 2.1, Step 2]{bolley_quantitative_2007}, let $(E, d)$ be a metric space with finite diameter $D$, $\nu \in \mathcal{P}(E)$ and $\mathcal{A}$ a measurable subset of $\mathcal{P}(E)$. Then for all $\delta \in (0, D)$,
    \begin{equation*}
        \begin{aligned}
            \P\left(\widehat\nu^N\in\mathcal A\right)
            \le{}&\left(\frac{8eD}{\delta}\right)^{pN(E,\delta/2)}\\
            &\times\exp\left(-N\inf_{\nu'\in\mathcal A_\delta}
                     H(\nu'\mid\nu)\right),
        \end{aligned}
    \end{equation*}
    where $\mathcal{A}_\delta$ denotes the $\delta$-thickening in $W_p$ distance of $\mathcal{A}$ given by
    \begin{equation*}
        \mathcal{A}_\delta \coloneqq \Big\{\nu' \in \mathcal{P}(E) \: : \: \exists \nu_\mathcal{A} \in \mathcal{A},  W_p(\nu', \nu_\mathcal{A}) \le \delta \Big\}.
    \end{equation*}
\end{proposition}
For unbounded spaces, the diameter is infinite. But if the probability measure decays sufficiently fast, then most of the samples lie in some large ball for which this quantitative version of Sanov's theorem holds. The idea is to truncate the space into a large ball for some optimal radius balancing the metric entropy contribution with the tail error.

We generalize the quantitative concentration inequality \cite[Theorem 2.1]{bolley_quantitative_2007} of Bolley--Guillin--Villani for independent random variables satisfying transport-cost inequalities from Euclidean spaces $\R^n$ to the general $\MCP(K,n)$ setting. More precisely, we have the following result.
\begin{theorem}[Quantitative Wasserstein concentration]\label{thm:W_p quantitative concentration}
    Let $(E, d, \m)$ be an unbounded metric measure space satisfying $\MCP(K,n)$ for some $K \le 0$ and $n > 1$. Let $1 \le p \le 2$ and $\nu$ be a probability measure on $E$ that satisfies a $T_p(\lambda)$ inequality \eqref{eq:Tp lambda}. Then, for every $n' > n$ and $0 < \lambda' < \lambda$, there exists a constant $N_0$ depending only on $K,\lambda,\lambda',n,n',p$ and on the square-exponential moment of $\nu$, such that for all $0 < \varepsilon < 1$ and for all $N \ge N_0 \varepsilon^{-(n'+2)}$,
    \begin{equation}\label{eq:W_p(nu, nu^N) statement theorem}
        \P\left(W_p(\widehat\nu^N,\nu)>\varepsilon\right)
        \le
        \exp\left(-\frac{\gamma_p\lambda'}{2}N\varepsilon^2\right),
\end{equation}
    where
    \begin{equation*}
        \gamma_p = \begin{cases}
            1 & \text{if } 1\le p<2; \\
            3 - 2\sqrt{2} & \text{if } p=2.
        \end{cases}
    \end{equation*}
\end{theorem}

\begin{remark}[Comparison with concentration around the expectation]
\label{rem:centered-vs-BGV}
    Under the hypotheses of \cref{cor:squareexp-expectation},
    when $n>2p$, the bound
    \[
        \P\left(W_p(\widehat\nu^N,\nu)>\varepsilon\right)
        \le \exp\left(-\frac{\lambda'}{2}N\varepsilon^2\right)
    \]
    holds as soon as
    \[
        N\ge
        \left(
            \frac{C}
            {(1-\sqrt{\lambda'/\lambda})\varepsilon}
        \right)^n
    \]
    with the constant $C$ in \eqref{eq:optimal estimate MCP E alpha} for an admissible choice of $\alpha$. This improves the sample-size exponent in
    \cref{thm:W_p quantitative concentration} from $\varepsilon^{-(n'+2)}$ to $\varepsilon^{-n}$.
    Moreover, the loss in the constant $\gamma_p$
    when $p=2$ in \eqref{eq:W_p(nu, nu^N) statement theorem} is an artifact of the proof rather than
    an intrinsic limitation of the transport-cost inequality,
    as the concentration estimate around the expectation
    holds with the same Gaussian constant for every
    $1\le p\le2$.
    Nevertheless, \cref{thm:W_p quantitative concentration}
    remains relevant in the complementary regime $n\le2p$,
    where the expectation estimate does not yield the rate
    $N^{-1/n}$ used above.
\end{remark}
We also generalize \cite[Theorem 2.7]{bolley_quantitative_2007} which covers the case of measure only satisfying polynomial moments.
\begin{theorem}[Concentration under polynomial moments]\label{thm:W_p quantitative concentration polynomial}
    Let $(E, d, \m)$ be an unbounded metric measure space satisfying $\MCP(K,n)$ for some $K \le 0$ and $n > 1$. Let $q \ge 1$ and $\nu$ be a probability measure that satisfies the polynomial moment condition
    \begin{equation*}
        M_q \coloneqq \int_E d^q(x_0, x) \dd\nu(x) < +\infty. 
    \end{equation*}
    Set, for $p \in [1, q)$, the constant
    \[C_{p,q} \coloneqq
        \max\left\{
            (16M_q^{1/p})^{p/(q-p)},
            (2^{2p+1}M_q)^{1/(q-p)},
            1
        \right\}.\]
    Then:
    \begin{enumerate}[label=(\roman*)]
        \item For all $p \in [1, q/2), 0 < \delta < q/p - 2$, and  $n' > n$, there is a constant $N_0$ depending only on $\delta,p,q,n,n'$ and $M_q$, such that
        \begin{equation*}\label{eq:W_p(nu, nu^N) polynomial case (i) final}
            \P\left(W_p(\widehat\nu^N,\nu) > \varepsilon\right) \le \varepsilon^{-q}N^{-q/(2p)+\delta/2}
        \end{equation*}
        for all $0 < \varepsilon < 1$ and
        \begin{equation*}
            N \ge N_0 \varepsilon^{-\frac{q(n'+2p)}{q-p}}\exp\left(2C_{p,q}\sqrt{|K|(n-1)}\varepsilon^{-p/(q-p)}\right).
        \end{equation*}
        
        \item For all $p \in [q/2, q), 0 < \delta < q/p - 1$, and $n'>n$, there is a constant $N_0$ depending only on $\delta,p,q,n,n'$ and $M_q$ such that
        \begin{equation*}\label{eq:W_p(nu, nu^N) polynomial case (ii) final}
            \P\left(W_p(\widehat\nu^N,\nu) > \varepsilon\right) \le \varepsilon^{-q}N^{1-q/p + \delta}
        \end{equation*}
        for all $0 < \varepsilon < 1$ and
        \begin{equation*}
            N \ge N_0 \varepsilon^{-\frac{q(n'+2p)}{q-p}}\exp\left(2C_{p,q}\sqrt{|K|(n-1)}\varepsilon^{-p/(q-p)}\right).
        \end{equation*}
    \end{enumerate}
\end{theorem}

\begin{remark}[Exponential and polynomial concentration regimes]
    The regime in \cref{thm:W_p quantitative concentration} is the same as in 
    \cite[Theorem 2.1]{bolley_quantitative_2007}, namely
    \begin{equation*}
        N\ge N_0\varepsilon^{-(n'+2)},
    \end{equation*}
    up to replacing the Euclidean dimension by any effective dimension $n' > n$ on $(E, d, \m)$ satisfying $\MCP(K,n)$.

    In contrast, under only a polynomial moment assumption, \cref{thm:W_p quantitative concentration polynomial} holds in the regime
    \begin{equation}\label{eq:this regime}
        N \ge
        N_0\varepsilon^{-\frac{q(n'+2p)}{q-p}}
        \exp\left(
            2C_{p,q}\sqrt{|K|(n-1)}
            \varepsilon^{-p/(q-p)}
        \right).
    \end{equation}
    When $K=0$, the exponential factor disappears and one recovers the exact same polynomial regime in $1/\varepsilon$, as in the Euclidean case \cite[Theorem 2.7]{bolley_quantitative_2007}. 
    When $K<0$, however, the covering estimate introduces an additional exponential factor in \eqref{eq:this regime}, reflecting the exponential volume growth allowed under a negative curvature lower bound. This is consistent with the quantization phenomena discussed above: on negatively curved unbounded spaces, polynomial moments alone are in general not sufficient to compensate for exponential volume growth. In particular, as suggested by the counterexamples to the asymptotic quantization formula without suitable exponential moments, one should not expect estimates of the same strength as in the Euclidean case to hold in a purely polynomial regime in $1/\varepsilon$ under only polynomial moment assumptions.
\end{remark}

\subsection{Concentration for aggregation--diffusion}
    
    As an important application, we study mean-field
    interacting particle systems of McKean--Vlasov type
    on Riemannian manifolds. The particle system is
    \begin{equation}\label{eq:SDE interacting system particles 1}
        dX_t^i = \sqrt2\,dB_t^{i,M}-\nabla V(X_t^i)\,dt
        -\frac1N\sum_{j=1}^N\nabla_xW(X_t^i,X_t^j)\,dt,
        \qquad i=1,\ldots,N, 
    \end{equation}
    where $dB_t^M$ denotes the Riemannian Brownian motion on $M$, $V$ a confinement potential and $W$ an interaction potential. The associated empirical measure
    \begin{equation*}
        \widehat\mu_t^N \coloneqq \tfrac{1}{N}\sum_{i=1}^N \delta_{X_t^i}
    \end{equation*}
    has limiting law $\mu_t$ whose evolution is governed by a nonlinear Fokker--Planck equation, also called the aggregation--diffusion equation,
    \begin{equation*}
        \partial_t \mu_t = \Delta \mu_t + \nabla \cdot \Big(\mu_t\left[\nabla V + \nabla \left(W * \mu_t\right)\right]\Big),
    \end{equation*}
    where we define the convolution as
    \[  W * \mu_t(x) \coloneqq \int_{M} W(x,y) \dd\mu_t(y).\]
    As $N$ becomes large, the empirical interaction approaches its
    mean-field counterpart. We apply the independent-sample concentration
    estimate to the noninteracting proxy and compare the two systems by
    parallel coupling. The following result extends the strategy of
    \cite[Theorem 2.9]{bolley_quantitative_2007} to Riemannian manifolds with Ricci curvature bounded below.
    \begin{theorem}[Concentration for aggregation--diffusion on manifolds]
    \label{thm:main aggregation_diffusion}
        Let $(M,g)$ be a complete, connected Riemannian manifold of dimension
        $m\ge2$ satisfying $\Ric\ge Kg$ for some $K\le0$,
        and fix $1\le p\le2$.
        Assume that $V$ and $W$ are smooth and satisfy the Hessian,
        growth, and Lipschitz conditions
        \eqref{eq:assumptions V IPS}--
        \eqref{eq:assumption W second variable Lipschitz IPS}.
        Let $\mu_0$ have a finite square-exponential moment:
        \begin{equation*}
            E_{\alpha_0}
            \coloneqq\int_M e^{\alpha_0d^2(o,x)}\dd\mu_0(x)<\infty
        \end{equation*}
        for some $o\in M$ and $\alpha_0>0$. If $1<p\le2$,
        assume in addition that $\mu_0$ satisfies
        $T_2(\lambda_{2,0})$ for some $\lambda_{2,0}>0$.
        Let $(X_t^i)_{i=1}^N$ be the interacting particle system
        \eqref{eq:SDE interacting system particles 1} with
        independent initial positions of law $\mu_0$.
        Write $\mu_t$ for the McKean--Vlasov law and
        $\widehat\mu_t^N=N^{-1}\sum_{i=1}^N\delta_{X_t^i}$.

        Then, for every $T>0$ and $m'>m$, there exist constants
        $C_{T,p},c_{T,p},N_0>0$, independent of $N$ and
        $\varepsilon$, such that
        \begin{equation}\label{eq:main aggregation concentration}
            \P\left(\sup_{0\le t\le T}
                W_p(\widehat\mu_t^N,\mu_t)>\varepsilon\right)
            \le C_{T,p}(1+T\varepsilon^{-2})
                e^{-c_{T,p}N\varepsilon^2}
        \end{equation}
        for every $0<\varepsilon<1$ and
        $N\ge N_0\varepsilon^{-(m'+2)}$.
    \end{theorem}

\section{Proofs for independent random variables}\label{sec:Proofs for independent random variables}
\subsection{Expected Wasserstein distance}
Before detailing the proof of \cref{thm:sharp E nu n n K neg}, we give some outline. \cref{thm:BLG} yields bounds in expectation on spaces with finite diameter depending on their covering number. In order to obtain a sharp rate when the diameter is infinite, one cannot simply use a single truncation and compare the rate on the compact with the error term outside of that compact which depends on exponential decay of the tail. Instead, one uses an increasing sequence of truncations covering the whole space and rewrite the empirical measure as a series of empirical measure on each annuli up to some probability measure. On each of these truncations one applies \cref{thm:BLG}. The choice of the sequence of annuli depends on optimization considerations.

\begin{proof}[Proof of \cref{thm:sharp E nu n n K neg}]
    We decompose the space into annuli as in the proof of Boissard--Le Gouic \cite[Corollary 1.3]{boissard_mean_2014}. 
    
    Set $r_0 \coloneqq 0$, and
    \begin{equation*}
         A_i \coloneqq B(o,r_i)\setminus B(o,r_{i-1}), \qquad i \ge 1 \text{ and } r_i > r_{i-1}.
    \end{equation*}
    Let $X^1, \dots, X^N$ be i.i.d. with common law $\nu$ and denote by 
    \[\widehat\nu^N \coloneqq \frac{1}{N}\sum_{j=1}^N \delta_{X^j}\]
    the associated empirical measure. For each $i \ge 1$, we define $N_i \coloneqq N \widehat\nu^N(A_i)$, and whenever $N_i \neq 0$, we let
    \begin{equation*}
        \widehat\nu_i^N \coloneqq \frac{\mathbf{1}_{A_i}\widehat\nu^N}{\widehat\nu^N(A_i)}, \qquad \nu_i \coloneqq \frac{\mathbf{1}_{A_i}\nu}{\nu(A_i)}.
    \end{equation*}
    Notice that $N_i \in \N$ and $\widehat\nu_i^N$ is a conditional empirical measure. Then there exist probability measures $\mu_1$ and $\mu_2$ such that
    \begin{equation}\label{eq:decomposition}
        \widehat\nu^N = \sum_{i \ge 1}\lambda_i \widehat\nu_i^N + \lambda \mu_1, \qquad \nu = \sum_{i \ge 1}\lambda_i \nu_i + \lambda \mu_2,
    \end{equation}
    where \begin{equation*} 
        \lambda_i \coloneqq \widehat\nu^N(A_i)\wedge \nu(A_i), \qquad \lambda \coloneqq \sum_{i \ge 1}\left(\widehat\nu^N(A_i) - \nu(A_i)\right)_+.
    \end{equation*}
    Indeed, one can check, under the convention that whenever $\widehat\nu^N(A_i) = 0$ then $0 \times \widehat\nu_i^N \coloneqq 0$, that 
    \begin{align*}
        \mu_1&\coloneqq\frac1\lambda\sum_{i\ge1}
            \big(\widehat\nu^N(A_i)-\nu(A_i)\big)_+
            \widehat\nu_i^N,\\
        \mu_2&\coloneqq\frac1\lambda\sum_{i\ge1}
            \big(\nu(A_i)-\widehat\nu^N(A_i)\big)_+
            \nu_i.
    \end{align*}
    By \cite[Lemma 2.4]{boissard_mean_2014}, namely convexity of $W_p^p$ and the triangle inequality, and the decomposition \eqref{eq:decomposition}, we therefore have \begin{equation}\label{eq:annuli-decomposition-proof} 
        \begin{aligned}
            W_p(\widehat\nu^N,\nu)
            \le{}& \sum_{i\ge1}\lambda_i^{1/p}
                W_p(\widehat\nu_i^N,\nu_i)\\
            &+\lambda^{1/p}W_p(\nu,\mu_1)
             +\lambda^{1/p}W_p(\nu,\mu_2).
        \end{aligned} 
        \end{equation} 
    We estimate separately the shell contribution and the two remainder terms. 
        
    \noindent\textbf{Step 1 (Estimate on the remainder terms).} We only treat the term involving $\mu_1$, the one involving $\mu_2$ being identical. Again by \cite[Lemma 2.4]{boissard_mean_2014},
    \begin{equation*}
        \lambda^{1/p}W_p(\nu,\mu_1) \le \sum_{i \ge 1} \left(\widehat\nu^N(A_i)-\nu(A_i)\right)_+^{1/p} W_p(\nu,\widehat\nu_i^N). 
    \end{equation*} 
    By triangle inequality, $ W_p(\nu,\widehat\nu_i^N) \le W_p(\nu, \delta_o) + W_p(\delta_o, \widehat\nu_i^N)$, and we deduce that 
    \begin{equation*} 
        W_p(\nu,\widehat\nu_i^N) \le \left(\int_E d^p(o,x) \dd\nu(x)\right)^{1/p} + \left( \frac{1}{\widehat\nu^N(A_i)} \int_{A_i}d^p(o, x) \dd\widehat\nu^N(x)\right)^{1/p} \le M_p^{1/p} + r_i, 
    \end{equation*} 
    where we denote by $M_p \coloneqq \int_E d^p(o, x) \, d\nu(x)$ the $p$-th moment of $\nu$ which is finite, since $\nu$ admits some exponential moment $\mathcal E_\alpha$. 
    Taking expectations, we obtain 
    \begin{equation}\label{eq:to become}
        \begin{aligned}
            \E\left[\lambda^{1/p}W_p(\nu,\mu_1)\right]
            \le{}&\sum_{i\ge1}
            \E\left[\big(\widehat\nu^N(A_i)-\nu(A_i)\big)_+^{1/p}\right]\\
            &\quad\times\big(M_p^{1/p}+r_i\big).
        \end{aligned}
\end{equation}
    By Jensen's inequality, 
    \begin{equation*} 
        \E\left[\left(\widehat\nu^N(A_i)-\nu(A_i)\right)_+^{1/p}\right] \le \E\left[\left|\widehat\nu^N(A_i)-\nu(A_i)\right|^{1/p}\right] \le \E\left[\left|\widehat\nu^N(A_i)-\nu(A_i)\right|^2\right]^{1/(2p)}. 
    \end{equation*} 
    Since $N_i = N\widehat\nu^N(A_i) \sim \mathrm{Bin}(N,\nu(A_i))$, $\Var\left[\widehat\nu^N(A_i)\right] \le \tfrac{\nu(A_i)}{N}$, and thus \eqref{eq:to become} becomes \begin{equation}\label{eq:remainder-m1}
        \begin{aligned}
            \E\left[\lambda^{1/p}W_p(\nu,\mu_1)\right]
            \le{}&N^{-1/(2p)}\sum_{i\ge1}
            \nu(A_i)^{1/(2p)}\big(M_p^{1/p}+r_i\big).
        \end{aligned}
\end{equation} 
    The same estimate holds with $\mu_2$ in place of $\mu_1$. 
    We now estimate the shell masses using the exponential moment. 
    For every $i \ge 1$, 
    \begin{equation*} 
        \nu(A_i) \le \nu\left(E\setminus B(o,r_{i-1})\right) \le e^{-\alpha r_{i-1}}\int_E e^{\alpha d(o,x)}\dd\nu(x) = \mathcal E_\alpha e^{-\alpha(r_{i-1})}. \end{equation*} 
    Therefore, 
    \begin{equation}\label{eq:sufficient sum}
        \begin{aligned}
            \sum_{i\ge1}\nu(A_i)^{1/(2p)}\big(M_p^{1/p}+r_i\big)
            \le{}&\mathcal E_\alpha^{1/(2p)}
                \sum_{i\ge1}e^{-\alpha r_{i-1}/(2p)}
                \big(M_p^{1/p}+r_i\big).
        \end{aligned}
\end{equation} 
    We denote the right-hand side by $C_1$, which is finite for sequences $\{r_i\}_{i\ge0}$ increasing sufficiently fast.
    We deduce from \eqref{eq:remainder-m1} that for such sequences, 
    \begin{equation}\label{eq:remainder-final-proof}
        \E\left[\lambda^{1/p}W_p(\nu,\mu_1)
            +\lambda^{1/p}W_p(\nu,\mu_2)\right]
        \le 2C_1N^{-1/(2p)}.
\end{equation}

    \noindent\textbf{Step 2 (Conditional estimate on one shell).}
    Fix $i \ge 1$ and condition on the event $\{N_i = k\}$ for some $k \ge 1$. Conditionally on this event, $\widehat\nu_i^N$ is the empirical measure of $k$ independent random variables distributed according to $\nu_i$.
    
    Since $\nu_i$ is supported on $A_i \subset B(o,r_i)$, the uniform-in-$K$ covering estimate \eqref{eq:covering estimate uniform K} yields, for every $0 < \delta \le r_i/2$,
    \begin{equation}\label{eq:covering recall}
        N(A_i,\delta) \le N(B(o,r_i),\delta) \le\left(\frac{5r_i}{\delta}\right)^n\exp\left(2\sqrt{|K|(n-1)}r_i\right).
    \end{equation}
    Applying \cref{thm:BLG} with sample size $k$ together with \eqref{eq:covering recall}, we obtain that for every $0 < \delta < r_i/2$,
    \begin{align*}
        \E\left[W_p(\widehat\nu_i^N,\nu_i)\,\middle|\,N_i = k \right] 
        &\le \frac{64}{3}\left(\delta +
            k^{-1/(2p)}\int_\delta^{r_i/2} N(A_i,s)^{1/(2p)}\dd s
        \right) \\
        &\le \frac{64}{3}\left(
            \delta
            +
            \frac{5^{n/(2p)}}{\frac{n}{2p}-1}
            k^{-1/(2p)} r_i^{n/(2p)}e^{\sqrt{|K|(n-1)}r_i/p}\delta^{1-n/(2p)}
        \right). \numberthis\label{eq:previous estimate}
    \end{align*}
    We now choose $\delta$ to equate the two terms of the right-hand side;
    \begin{equation*}
        \delta \coloneqq \delta_{i,k}
        \coloneqq
        \frac{5}{\left(\tfrac{n}{2p}-1\right)^{2p/n}}
        r_i e^{\tfrac{2\sqrt{|K|(n-1)}\,r_i}{n}} k^{-1/n}.
    \end{equation*}
    If $\delta_{i,k} \le r_i/2$, then \eqref{eq:previous estimate} yields
    \begin{equation}\label{eq:one-shell-final-short}
        \E\left[W_p(\widehat\nu_i^N,\nu_i)\,\middle|\,N_i = k \right]
        \le C_2\, r_i e^{\tfrac{2\sqrt{|K|(n-1)}\,r_i}{n}} k^{-1/n}
    \end{equation}
    with $C_2 \coloneqq \tfrac{640}{3} \times (\tfrac{n}{2p}-1)^{-2p/n}$. Note that when $\delta_{i,k} > r_i/2$, we can estimate
    \begin{equation*}
        W_p(\widehat\nu_i^N,\nu_i) \le \operatorname{diam}(A_i) \le 2r_i < 4\delta_{i,k},
    \end{equation*}
    and thus the bound \eqref{eq:one-shell-final-short} is valid for every $k \ge 1$.

    \medskip

    \noindent\textbf{Step 3 (Averaging the shell terms).}
    Combining $\lambda_i \le \widehat\nu^N(A_i)=N_i/N$ with the estimate \eqref{eq:one-shell-final-short}, we obtain
    \begin{equation*}
        \E\left[\lambda_i^{1/p}W_p(\widehat\nu_i^N,\nu_i)\,\middle|\, N_i = k\right]
        \le
        C_2 N^{-1/p} r_i e^{\tfrac{2\sqrt{|K|(n-1)}\,r_i}{n}} k^{1/p-1/n}.
    \end{equation*}
    By the tower property,
    \begin{align*}
        \E\left[\lambda_i^{1/p}W_p(\widehat\nu_i^N,\nu_i)\right]
        &\le
        C_2 N^{-1/p} r_i e^{\tfrac{2\sqrt{|K|(n-1)}\,r_i}{n}}
        \E\left[N_i^{1/p-1/n}\right] \\
        &\le C_2 N^{-1/n} r_i e^{\tfrac{2\sqrt{|K|(n-1)}\,r_i}{n}} \nu(A_i)^{1/p-1/n};
    \end{align*}
    we used Jensen's inequality in the second line and $N_i = N\widehat\nu^N(A_i)$.
    Summing over $i \ge 1$ and recalling that $\nu(A_i) \le \mathcal E_\alpha e^{-\alpha r_{i-1}}$ with $\mathcal{E}_\alpha = \int_E e^{\alpha d(o, x)}\, d\nu(x)$,
    \begin{multline}\label{eq:shell-sum-before-choice}
        \sum_{i \ge 1} \E\left[\lambda_i^{1/p}W_p(\widehat\nu_i^N,\nu_i)\right]
        \\ \le
        C_2 \mathcal E_\alpha^{1/p-1/n} N^{-1/n} \sum_{i \ge 1} r_i \exp\left(\frac{2\sqrt{|K|(n-1)}}{n}r_i - \alpha\left(\frac{1}{p}-\frac{1}{n}\right)r_{i-1}\right).
    \end{multline}

    \noindent\textbf{Step 4 (Choice of the radii).}
    Denoting $h_i \coloneqq r_i-r_{i-1}$, the exponent in the last sum can be rewritten as
    \begin{equation*}
        -\beta r_{i-1} + \frac{2\sqrt{|K|(n-1)}}{n}h_i, \qquad \beta \coloneqq \alpha\left(\frac{1}{p}-\frac{1}{n}\right)-\frac{2\sqrt{|K|(n-1)}}{n}.   
    \end{equation*}
    In order for \eqref{eq:shell-sum-before-choice} to be summable, we therefore require $\beta > 0$, which is the condition \eqref{eq:assumption on alpha}.
    
    A convenient choice of radii is $r_i = i$ so that $h_i = 1$ for $i \ge 1$, and 
    \begin{equation*}
        \sum_{i \ge 1} r_i \exp\left(\frac{2\sqrt{|K|(n-1)}}{n}r_i - \alpha\left(\frac{1}{p}-\frac{1}{n}\right)r_{i-1}\right)
        =
        e^{\tfrac{2\sqrt{|K|(n-1)}}{n}}\frac{1}{(1-e^{-\beta})^2}.
    \end{equation*}
    Inserting this into \eqref{eq:shell-sum-before-choice}, we obtain
    \begin{equation}\label{eq:shell-final}
        \sum_{i \ge 1} \E\left[\lambda_i^{1/p}W_p(\widehat\nu_i^N,\nu_i)\right]
        \le
        C_3 N^{-1/n},
    \end{equation}
    where
    \begin{equation*}
        C_3
        \coloneqq
        C_2 \mathcal E_\alpha^{1/p-1/n}\,e^{\tfrac{2\sqrt{|K|(n-1)}}{n}}\frac{1}{(1-e^{-\beta})^2}.
    \end{equation*}

    Combining \eqref{eq:annuli-decomposition-proof}, \eqref{eq:remainder-final-proof}, and \eqref{eq:shell-final}, we obtain
    \begin{equation*}
        \E\left[W_p(\widehat\nu^N,\nu)\right]
        \le
        2C_1N^{-1/(2p)} + C_3 N^{-1/n},
    \end{equation*}
    where $C_1$ is defined in \eqref{eq:sufficient sum}; 
    \begin{equation*}
        \begin{aligned}
            C_1
            &=\mathcal E_\alpha^{1/(2p)}\sum_{i\ge1}
                e^{-\alpha r_{i-1}/(2p)}\big(M_p^{1/p}+r_i\big)\\
            &=\mathcal E_\alpha^{1/(2p)}
                \left(\frac{M_p^{1/p}}{1-e^{-\alpha/(2p)}}
                +\frac{1}{\big(1-e^{-\alpha/(2p)}\big)^2}\right).
        \end{aligned}
    \end{equation*}
    Since $n>2p$, one has $N^{-1/(2p)}\le N^{-1/n}$ for every $N\ge1$. Hence
    \begin{equation*}
        \E\left[W_p(\widehat\nu^N,\nu)\right]
        \le
        (2C_1+C_3)N^{-1/n}.
    \end{equation*}
    This concludes the proof.
\end{proof}

\subsection{Concentration inequalities}\label{sec:Concentration inequalities for independent observables}
We now adapt the quantitative Sanov proof of Bolley--Guillin--Villani, which does not require the dimensional restriction $n>2p$. We carefully track the constants in order to capture the precise effect of the lower curvature bound, or more generally of the measure contraction property.

\begin{proof}[Proof of \cref{thm:W_p quantitative concentration}]
    We adapt the strategy of Bolley, Guillin and Villani \cite[Theorem 2.1]{bolley_quantitative_2007}.

    \noindent \textbf{Step 1 (Truncation).} Denote
    \begin{equation*}
        \nu_R \coloneqq \frac{\mathbf{1}_{B_R(o)}\nu}{\nu\left[B_R(o)\right]}.
    \end{equation*}
    Then we sample independent random variables $(X^i)_{1 \le i \le N}$ and $(Y^i)_{1\le i \le N}$ according to $\nu$ and its truncation $\nu_R$ respectively in order to define
    \begin{equation*}
        X_R^i = \begin{cases}
            X^i & \text{if } d(o, X^i) \le R; \\
            Y^i & \text{if } d(o, X^i) > R.
        \end{cases}
    \end{equation*}
    The $p$-Wasserstein distance between the measure and its truncation can be estimated as
    \begin{align}
        W_p^p(\nu_R,\nu)
        &\le \E\left[d^p(X^1,Y^1)
            \mathbf{1}_{\{d(o,X^1)>R\}}\right] \notag\\
        &\le 2^p\int_{E\setminus B_R(o)}d^p(o,x)\dd\nu(x)
        \le 2^p E_\alpha R^p e^{-\alpha R^2}.
        \label{eq:W_p^p nuR nu}
    \end{align}
    for $R$ large enough, say $R \ge R_0 \coloneqq \sqrt{p/(2\alpha)}$, since $r \mapsto r^p e^{-\alpha r^2}$ is nonincreasing for $r \ge R$; recall that $\nu$ satisfies in particular $T_1(\lambda)$, and hence there is some $\alpha > 0$ such that
    \begin{equation*}
        E_\alpha \coloneqq \int_E e^{\alpha d^2(o, x)} \dd\nu(x) < +\infty.
    \end{equation*}
    
    We now compare the $p$-Wasserstein distance between the empirical measure and its truncation defined as
    \begin{equation*}
        \widehat\nu^N := \frac{1}{N}\sum_{i=1}^N \delta_{X^i}, \qquad \widehat\nu_R^N := \frac{1}{N}\sum_{i=1}^N \delta_{X_R^i}.
    \end{equation*}
    Set 
    \begin{equation*}
        Z^i \coloneqq 2^p d^p(o, X^i) \mathbf{1}_{\{d(o,X^i)>R\}}\mathbf{1}_{\{1\le p<2\}} + d^2(X^i, X_R^i)\mathbf{1}_{\{p=2\}}
    \end{equation*}
    and let $\varepsilon, \theta > 0$. The exponential Chebyshev inequality yields
    \begin{align*}
        \P\left(W_p(\widehat\nu^N, \widehat\nu_R^N) > \varepsilon\right) &\le \P\left(\frac{1}{N}\sum_{i=1}^N d^p(X^i, X_R^i) > \varepsilon^p\right) \\
        &\le \exp\Big(-N\left\{\theta\varepsilon^p - \log \E\left[\exp\left(\theta Z^1\right)\right]\right\}\Big) \numberthis\label{eq:expectation Z_1}
    \end{align*}
    by independence of the $Z^i$'s. 

    In the case $1 \le p < 2$, let $0 < \alpha_1 < \alpha$. Then for $R$ large enough, say $R \ge R_1 \coloneqq (2^p\theta/\alpha_1)^{1/(2-p)}$, using the inequalities $2^p\theta r^p \le \alpha_1 r^2$ valid for $r \ge R$ and $\log(1 + x) \le x$, we obtain
    \begin{equation}\label{eq:upper bound log E exp theta Z1}
        \begin{aligned}
            \log\E\left[e^{\theta Z^1}\right]
            &\le \log\left(1+\int_{E\setminus B_R(o)}
                e^{2^p\theta d^p(o,x)}\dd\nu(x)\right)\\
            &\le E_\alpha e^{-(\alpha-\alpha_1)R^2}.
        \end{aligned}
\end{equation}
    Inserting \eqref{eq:upper bound log E exp theta Z1} into \eqref{eq:expectation Z_1}, we deduce that
    \begin{equation}\label{eq:W_p(nu^N, nu_R^N)}
        \P\left(W_p(\widehat\nu^N,\widehat\nu_R^N)>\varepsilon\right)
        \le \exp\left(-N\left\{
            \theta\varepsilon^p-E_\alpha e^{-(\alpha-\alpha_1)R^2}
        \right\}\right).
\end{equation}

    In the case $p=2$, set $\theta \coloneqq \alpha_1/2$ for some $0<\alpha_1<\alpha$. Then for $R$ large enough, say $R \ge R_2 \coloneqq \sqrt{\log(2E_\alpha)/\alpha}$, so that $\nu[B_R(o)] \ge 1 - E_\alpha e^{-\alpha R^2} \ge \tfrac{1}{2}$, it holds
    \begin{align*}
        \E\left[\exp\left(\theta Z^1\right)\right] &\le \nu[B_R(o)]
        +
        \frac{1}{\nu[B_R(o)]}
        \int_{B_R(o)}\int_{E\setminus B_R(o)} e^{\frac{\alpha_1}{2}d^2(x,y)}\dd\nu(x)\dd\nu(y) \\
        &\le 1 + 2E_\alpha^2 e^{-(\alpha - \alpha_1)R^2}. \numberthis\label{eq:upper bound log E exp theta Z1 p2}
    \end{align*}
    Inserting \eqref{eq:upper bound log E exp theta Z1 p2} into \eqref{eq:expectation Z_1}, we deduce that
    \begin{equation}\label{eq:W_2(nu^N, nu_R^N)}
        \P\left(W_2(\widehat\nu^N,\widehat\nu_R^N)>\varepsilon\right)
        \le \exp\left(-N\left\{
            \frac{\alpha_1}2\varepsilon^2-2E_\alpha^2e^{-(\alpha-\alpha_1)R^2}
        \right\}\right).
\end{equation}

    We now gather the estimates obtained so far. The triangle inequality yields
    \begin{equation*}
        \begin{aligned}
            \P\left(W_p(\widehat\nu^N,\nu)>\varepsilon\right)
            \le \P\left(\begin{aligned}[t] &W_p(\widehat\nu^N,\widehat\nu_R^N)
                +W_p(\widehat\nu_R^N,\nu_R)\\[-1pt]
                &+W_p(\nu_R,\nu)>\varepsilon\end{aligned}\right).
        \end{aligned}
    \end{equation*}
    Let $\eta \in (0,1)$ and $R \ge \max\{R_0, R_1, R_2\}$. In the case $1 \le p < 2$, combining \eqref{eq:W_p^p nuR nu} with \eqref{eq:W_p(nu^N, nu_R^N)}, we obtain
    \begin{equation}\label{eq:Wp nu^N nu vareps step 1}
        \begin{aligned}
        \P\left(W_p(\widehat\nu^N,\nu) > \varepsilon\right)
        &\le
        \P\left(W_p(\widehat\nu_R^N,\nu_R) > \eta\varepsilon - 2E_\alpha^{1/p}R e^{-\alpha R^2/p}\right) \\
        &\quad+
        \exp\left(-N\left\{\theta(1-\eta)^p\varepsilon^p-E_\alpha e^{-(\alpha-\alpha_1)R^2}\right\}\right).
        \end{aligned}
    \end{equation}
    In the case $p=2$, combining \eqref{eq:W_p^p nuR nu} with \eqref{eq:W_2(nu^N, nu_R^N)}, we similarly obtain
    \begin{equation}\label{eq:W2 nu^N nu vareps step 1}
        \begin{aligned}
        \P\left(W_2(\widehat\nu^N,\nu) > \varepsilon\right)
        &\le\P\left(W_2(\widehat\nu_R^N,\nu_R)>\eta\varepsilon - 2E_\alpha^{1/2}R e^{-\alpha R^2/2}\right) \\
        &\quad+
        \exp\left(-N\left\{\frac{\alpha_1}{2}(1-\eta)^2\varepsilon^2-2E_\alpha^2 e^{-(\alpha-\alpha_1)R^2}\right\}\right).
        \end{aligned}
    \end{equation}
    
    \noindent \textbf{Step 2 (Covering number estimate).} Bolley--Guillin--Villani obtain the following quantitative Sanov estimate on compacts (see Step 2 in the proof of \cite[Theorem 2.1]{bolley_quantitative_2007}). For every measurable subset $\mathcal A$ of $\mathcal P(B_R(o))$ and $0 < \delta < R$,
    \begin{equation}\label{eq:quantitative Sanov A}
        \begin{aligned}
            \P\left(\widehat\nu_R^N\in\mathcal A\right)
            &\le \left(\frac{16eR}{\delta}\right)^{pN(B_R(o),\delta/2)}\\
            &\quad\times\exp\left(-N\inf_{\nu'\in\mathcal A_\delta}
                H(\nu'\mid\nu_R)\right).
        \end{aligned}
\end{equation}
    where $\mathcal{A}_\delta$ denotes the $\delta$-thickening in $W_p$ distance of $\mathcal{A}$ given by
    \begin{equation*}
        \mathcal{A}_\delta \coloneqq \Big\{\nu' \in \mathcal{P}(B_R(o)) \: : \: \exists \nu_\mathcal{A} \in \mathcal{A},  W_p(\nu', \nu_\mathcal{A}) \le \delta \Big\}.
    \end{equation*}
    In light of the last two estimates \eqref{eq:Wp nu^N nu vareps step 1} and \eqref{eq:W2 nu^N nu vareps step 1} of Step 1, we set
    \begin{equation*}
        \mathcal A := \left\{ \nu' \in \mathcal{P}(B_R(o))\: : \: W_p(\nu', \nu_R) \ge \eta\varepsilon - 2E_\alpha^{1/p} R e^{-\alpha R^2/p} \right\}.
    \end{equation*}
    
    Now, it remains to find a lower bound to the relative entropy in \eqref{eq:quantitative Sanov A} with respect to the truncated measure. For every $\nu' \in \mathcal P(B_R(o))$ and $\nu[B_R(o)] \ge 1 - E_\alpha e^{-\alpha R^2}$, it holds
    \begin{align*}
        H(\nu'\mid\nu_R) =
        H(\nu'\mid\nu) + \log\left(\nu[B_R(o)]\right) \\\ge H(\nu'\mid\nu) - 2E_\alpha e^{-\alpha R^2} \\ \ge \frac{\lambda}{2}\,W_p^2\left(\nu', \nu\right) - 2E_\alpha e^{-\alpha R^2},
    \end{align*}
    where in the last inequality we used the fact that $\nu$ satisfies a $T_p(\lambda)$ inequality. Moreover, by the reverse triangle inequality and \eqref{eq:W_p^p nuR nu},
    \begin{equation*}
        W_p(\nu', \nu) \ge
        W_p(\nu', \nu_R) - W_p(\nu_R,\nu)
        \ge W_p(\nu', \nu_R) - 2E_\alpha^{1/p}R e^{-\alpha R^2/p}.
    \end{equation*}
    Therefore,
    \begin{equation*}
        H(\nu'\mid\nu_R)
        \ge
        \frac{\lambda}{2}
        \left(
            W_p(\nu', \nu_R)
            -
            2E_\alpha^{1/p}R e^{-\alpha R^2/p}
        \right)_+^2
        -
        2E_\alpha e^{-\alpha R^2}.
    \end{equation*}
    The positive part is essential here: the reverse triangle inequality gives a potentially negative lower bound. Fix now $\lambda_1 \in (0,\lambda)$. For $a,b\ge0$ we use the elementary inequality
    \begin{equation}\label{eq:elementary ineq}
        (a-b)_+^2 \ge \frac{\lambda_1}{\lambda}a^2 - \frac{\lambda_1}{\lambda-\lambda_1}b^2,
    \end{equation}
    which follows for $a\ge b$ from
    \begin{equation*}
        \begin{aligned}
            &(a-b)^2-\frac{\lambda_1}{\lambda}a^2
                +\frac{\lambda_1}{\lambda-\lambda_1}b^2\\
            &\qquad=\left(1-\frac{\lambda_1}{\lambda}\right)
                \left(a-\frac{b}{1-\lambda_1/\lambda}\right)^2\ge0.
        \end{aligned}
    \end{equation*}
    When $a<b$, the right-hand side of \eqref{eq:elementary ineq} is nonpositive, so the inequality is immediate. We deduce that
    \begin{align*}
        H(\nu'\mid\nu_R) &\ge
        \frac{\lambda_1}{2}W_p^2\left(\nu', \nu_R\right) - \frac{2\lambda\lambda_1}{\lambda-\lambda_1}
        E_\alpha^{2/p}R^2 e^{-\frac{2\alpha}{p}R^2} - 2E_\alpha e^{-\alpha R^2} 
        \\
        &\ge \frac{\lambda_1}{2}W_p^2\left(\nu', \nu_R\right) - C_{1}\,R^2 e^{-\alpha R^2} \numberthis\label{eq:entropy lower bound truncated}
    \end{align*}
    with
    \begin{equation*}
        C_1 \coloneqq\frac{2\lambda\lambda_1}{\lambda-\lambda_1}E_\alpha^{2/p} + 2E_\alpha.
    \end{equation*}
    From \eqref{eq:entropy lower bound truncated}, we have, for any $\nu \in \mathcal A_\delta$,
    \begin{equation*}
        H(\nu'\mid\nu_R) \ge \frac{\lambda_1}{2}W_p^2\left(\nu', \nu_R\right) - C_1 R^2 e^{-\alpha R^2} \ge \frac{\lambda_1}{2}\xi^2 - C_1 R^2 e^{-\alpha R^2},
    \end{equation*}
    where $\xi \coloneqq (\eta\varepsilon - 2E_\alpha^{1/p} Re^{-\alpha R^2/p} - \delta)_+$, and recalling the elementary inequality \eqref{eq:elementary ineq}, we find
    \begin{equation}\label{eq:H(nu', nu_R) lambda_2}
        H(\nu'\mid\nu_R) \ge \frac{\lambda_1}{2}\xi^2 - C_1R^2 e^{-\alpha R^2} \ge \frac{\lambda_2}{2}\varepsilon^2 - C_2 R^2 e^{-\alpha R^2} 
    \end{equation}
    in the regime $R \ge \max\{1,R_0,R_1,R_2\}$, for every $\lambda_2 < \lambda_1 \le \lambda$ such that $\lambda_2 < \lambda_1(\eta_1 - \delta_1)^2$, with the constant
    \begin{equation*}
        C_2 \coloneqq C_1 + \frac{2\lambda_1\lambda_2 E_\alpha^{2/p}}{\lambda_1\left(\eta_1 - \delta_1\right)^2 - \lambda_2},
    \end{equation*}
    for the choice $\delta \coloneqq \delta_1\varepsilon$ and $\eta \coloneqq \eta_1$. While \eqref{eq:H(nu', nu_R) lambda_2} is valid for $p=2$, in that case, we instead use
    \begin{equation}\label{eq:H(nu', nu_R) lambda_2 p=2}
        H(\nu'\mid\nu_R) \ge \frac{\lambda_1}{2}\xi^2 - C_1R^2 e^{-\alpha R^2} \ge \frac{\lambda_2}{2}\eta^2\varepsilon^2 - C_2^{(2)} R^2 e^{-\alpha R^2}
    \end{equation}
    for $\lambda_2 < \lambda_1(\eta - \delta_1)^2/\eta^2$ with the constant
    \[
    C_2^{(2)}
    \coloneqq
    C_1
    +
    \frac{2\lambda_1\lambda_2\eta^2E_\alpha}
    {\lambda_1(\eta-\delta_1)^2-\lambda_2\eta^2}.
    \]
    for the choice $\delta \coloneqq \delta_1\varepsilon$. 

    We now estimate the metric entropy term in \eqref{eq:quantitative Sanov A}.
    Use the covering estimate \eqref{eq:covering estimate uniform K} with $\delta/2$ instead of $\delta$, we get
    \begin{equation*}\label{eq:metric entropy raw}
        pN\left(B_R(o),\frac{\delta}{2}\right)
        \log\left(\frac{16eR}{\delta}\right)
        \le
        p\left(\frac{10R}{\delta}\right)^n
        \exp\left(2\sqrt{|K|(n-1)}\,R\right)
        \log\left(\frac{16eR}{\delta}\right).
    \end{equation*}
    Let now $n_1 > n$. Since $\log r \le \tfrac{r^{n_1-n}}{e(n_1-n)}$ for $r \ge 1$, we deduce that
    \begin{equation}\label{eq:metric entropy bound final}
        \begin{aligned}
            pN\left(B_R(o),\frac\delta2\right)
            \log\left(\frac{16eR}{\delta}\right)
            \le{}&C_3\exp\left(2\sqrt{|K|(n-1)}\,R\right)
                \left(\frac{R}{\delta}\right)^{n_1}.
        \end{aligned}
\end{equation}
    with
    \begin{equation*}
        C_3 \coloneqq \frac{p\,10^n(16e)^{n_1-n}}{e(n_1-n)}.
    \end{equation*}

    We insert this estimate into \eqref{eq:quantitative Sanov A}. In the case $1 \le p < 2$, combining \eqref{eq:quantitative Sanov A}, \eqref{eq:H(nu', nu_R) lambda_2}, and \eqref{eq:metric entropy bound final}, we obtain
    \begin{equation}\label{eq:step2 compact bound p<2 delta1}
        \P\left(\widehat\nu_R^N \in \mathcal A\right)
        \le
        \exp\left(
        \frac{C_3}{\delta_1^{n_1}}\exp\left(2\sqrt{|K|(n-1)}\,R\right)\left(\frac{R}{\varepsilon}\right)^{n_1}
            -
            N\left(\frac{\lambda_2}{2}\varepsilon^2 - C_2R^2e^{-\alpha R^2}\right)
        \right).
    \end{equation}
    In the case $p=2$, combining \eqref{eq:quantitative Sanov A}, \eqref{eq:H(nu', nu_R) lambda_2 p=2}, and \eqref{eq:metric entropy bound final}, we similarly obtain
    \begin{multline}\label{eq:step2 compact bound p=2}
        \P\left(\widehat\nu_R^N \in \mathcal A\right)
        \\ \le
        \exp\left(
            \frac{C_3}{\delta_1^{n_1}}
            \exp\left(2\sqrt{|K|(n-1)}\,R\right)
            \left(\frac{R}{\varepsilon}\right)^{n_1}
            -
            N\left(\frac{\lambda_2}{2}\eta^2\varepsilon^2 - C_2^{(2)}R^2e^{-\alpha R^2}\right)
        \right).
    \end{multline}

    Finally, recalling \eqref{eq:Wp nu^N nu vareps step 1} and \eqref{eq:W2 nu^N nu vareps step 1}, we conclude that in the case $1 \le p < 2$,
    \begin{equation}\label{eq:before choice parameters p<2}
    \begin{aligned}
        &\P\left(W_p(\widehat\nu^N,\nu) > \varepsilon\right) \\
        &\qquad\le
        \exp\left(
            \frac{C_3}{\delta_1^{n_1}}
            \exp\left(2\sqrt{|K|(n-1)}\,R\right)
            \left(\frac{R}{\varepsilon}\right)^{n_1}
            -
            N\left(\frac{\lambda_2}{2}\varepsilon^2 - C_2R^2e^{-\alpha R^2}\right)
        \right) \\
        &\qquad\quad+
        \exp\left(-N\left\{\theta(1-\eta_1)^p\varepsilon^p-E_\alpha e^{-(\alpha-\alpha_1)R^2}\right\}\right), 
    \end{aligned}
    \end{equation}
    while in the case $p=2$,
    \begin{equation}\label{eq:before choice parameters p=2}
    \begin{aligned}
        &\P\left(W_2(\widehat\nu^N,\nu) > \varepsilon\right) \\
        &\qquad \le
        \exp\left(
            \frac{C_3}{\delta_1^{n_1}}
            \exp\left(2\sqrt{|K|(n-1)}\,R\right)
            \left(\frac{R}{\varepsilon}\right)^{n_1}
            -
            N\left(\frac{\lambda_2}{2}\eta^2\varepsilon^2 - C_2^{(2)}R^2e^{-\alpha R^2}\right)
        \right) \\
        &\qquad \quad+
        \exp\left(-N\left\{\frac{\alpha_1}{2}(1-\eta)^2\varepsilon^2-2E_\alpha^2 e^{-(\alpha-\alpha_1)R^2}\right\}\right). 
    \end{aligned}
    \end{equation}

    \noindent \textbf{Step 3 (Choice of parameters).} 

    The last step consists of choosing parameters and is quite technical. Let us first recall the
parameters at play in the case \(1\le p<2\). The number \(n_1>n\) is an
auxiliary exponent introduced in the metric-entropy estimate, after absorbing
the logarithmic factor in the covering bound. The constants
\(0<\lambda_2<\lambda_1<\lambda\) come from the entropy lower bound after
truncation, where \(\nu\) satisfies \(T_p(\lambda)\). The parameters
\(\eta_1,\delta_1\in(0,1)\) are chosen so that
\[
    \lambda_2<\lambda_1(\eta_1-\delta_1)^2.
\]
Here \(\eta_1\) determines how much of the deviation \(\varepsilon\) is assigned
to the compact empirical-measure term, while \(\delta_1\) is the relative
thickening parameter in the quantitative Sanov estimate, namely
\(\delta=\delta_1\varepsilon\). Finally, \(\theta>0\), \(\alpha_1\in(0,\alpha)\),
and the cutoff radius \(R \gtrsim 1\) are still free. They will be chosen as functions of
\(\varepsilon\) so that the truncation errors and the metric-entropy contribution
are all absorbed into the leading Gaussian term.
    
    Now, fix $\lambda'\in(0,\lambda)$, and choose
    \begin{equation}\label{eq:parameters defined}
        \lambda_1 \coloneqq \frac{\lambda+\lambda'}{2}, 
        \qquad \lambda_2 \coloneqq \frac{\lambda_1+\lambda'}{2},
        \qquad \Delta \coloneqq \lambda_2-\lambda' >0.
    \end{equation}
    We treat separately the cases $1\le p<2$ and $p=2$.

    \noindent \textbf{Case 1: $1 \le p < 2$.}
    Set
    \begin{equation*}
        \delta_1 \coloneqq \frac14\left(1-\sqrt{\frac{\lambda_2}{\lambda_1}}\right),
        \qquad \eta_1 \coloneqq 1-\delta_1.
    \end{equation*}
    Then $0<\delta_1<\eta_1<1$ and $\lambda_2<\lambda_1(\eta_1-\delta_1)^2$; the constant in \eqref{eq:H(nu', nu_R) lambda_2} is well-defined.
    Choose 
    \begin{equation*}
        \theta(\varepsilon) \coloneqq \frac{\lambda_2}{(1-\eta_1)^p}\varepsilon^{2-p}, \quad
        \alpha_1 \coloneqq \frac{\alpha}{2}, 
    \end{equation*}
    and, denoting by $\log_+$ the nonnegative part of the logarithm,
    \begin{equation*}
        R(\varepsilon)
        \coloneqq
        \max\left\{
        R_0,\,
        R_1,\,
        R_2,\,
        1,\,
        \sqrt{
            \frac{2}{\alpha}
            \log_+\left(
                \frac{8E_\alpha}{\Delta\varepsilon^2}
            \right)
        },\,
        \sqrt{
            \frac{2}{\alpha}
            \log_+\left(
                \frac{16C_2}{\alpha e\Delta\varepsilon^2}
            \right)
        }
        \right\},
    \end{equation*}
    so that $\theta(1-\eta_1)^p\varepsilon^p = \lambda_2\varepsilon^2$,
    \begin{equation*}
         C_2R^2e^{-\alpha R^2} \le \frac{2C_2}{\alpha e}e^{-\alpha R^2/2} \le \frac{\Delta}{8}\varepsilon^2,
    \end{equation*}
    and
    \begin{equation*}
        E_\alpha e^{-(\alpha-\alpha_1)R^2} = E_\alpha e^{-\alpha R^2/2} \le \frac{\Delta}{8}\varepsilon^2.
    \end{equation*}
    Therefore, the right-hand side of \eqref{eq:before choice parameters p<2} is bounded by
    \begin{equation}\label{eq:both terms}
        \exp\left(
            A_K(\varepsilon)
            -
            N\left(\frac{\lambda_2}{2}-\frac{\Delta}{8}\right)\varepsilon^2
        \right)
        +
        \exp\left(
            -N\left(\lambda_2 -\frac{\Delta}{8}\right)\varepsilon^2
        \right),
    \end{equation}
    where
    \begin{equation}\label{eq:A_K def}
        A_K(\varepsilon) \coloneqq \frac{C_3}{\delta_1^{n_1}} \exp\left(2\sqrt{|K|(n-1)}R(\varepsilon)\right)\left(\frac{R(\varepsilon)}{\varepsilon}\right)^{n_1}.
    \end{equation}
    If in addition we require
    \begin{equation}\label{eq:required on A_K}
        A_K(\varepsilon) \le \frac{\Delta}{8}N\varepsilon^2,
    \end{equation}
    then both terms of \eqref{eq:both terms} are bounded by $e^{-\tfrac{\lambda_2 + \lambda'}{4} N\varepsilon^2}$.
    Whenever
    \begin{equation}\label{eq:second on N}
        N\varepsilon^2\ge \frac{4\log 2}{\Delta},
    \end{equation}
    this yields
    \begin{equation*}
        \P\left(W_p(\widehat\nu^N,\nu)>\varepsilon\right)
        \le
        \exp\left(-\frac{\lambda'}{2}N\varepsilon^2\right)
    \end{equation*}
    recalling \eqref{eq:parameters defined}.

    It remains to make the condition on $N$ explicit so that both regimes \eqref{eq:required on A_K} and \eqref{eq:second on N} hold. 
    Fix $n'>n$ and choose $n_1\coloneqq\tfrac12(n+n')$ and
    $\sigma\coloneqq n'-n_1=\tfrac12(n'-n)>0$.
    We choose a constant $C_4$ such that
    \begin{equation}\label{eq:Regime R C_4}
        R(\varepsilon)\le C_4\sqrt{\log(e/\varepsilon)}
        \qquad\text{for all }0<\varepsilon<1.
    \end{equation}
    Such a constant exists, for example
    \begin{align*}
        C_4
        \coloneqq
        \max&\Bigg\{
            \sqrt{\frac{p}{2\alpha}},\,
            \left(\frac{2^p\lambda_2}{\alpha_1(1-\eta_1)^p}\right)^{1/(2-p)},\,
            \sqrt{\frac{\log_+(2E_\alpha)}{\alpha}}, \\
            &\quad 1,\,\sqrt{
                \frac{2}{\alpha}
                \left(2+\log_+\frac{8E_\alpha}{\Delta} \vee \log_+\frac{16C_2}{\alpha e\Delta}\right)
            }
        \Bigg\}.
    \end{align*}
    Using
    \begin{equation*}
        \exp\left(2\sqrt{|K|(n-1)}C_4\sqrt{\log(e/\varepsilon)}\right) \le \exp\left(\frac{2|K|(n-1)C_4^2}{\sigma} + \frac{\sigma}{2}\right)
        \varepsilon^{-\tfrac{\sigma}{2}}
    \end{equation*}
    and
    \begin{equation*}
        \log\left(e/{\varepsilon}\right)^{n_1/2}
        \le
        \left(\frac{n_1}{e\sigma}\right)^{n_1/{2}}
        \left(\frac{e}{\varepsilon}\right)^{\tfrac{\sigma}{2}}
    \end{equation*}
    from the elementary inequality $u^a e^{-bu} \le (\tfrac{a}{be})^a$ valid for $u \ge 0, a,b, > 0$, recalling \eqref{eq:Regime R C_4}, we obtain
    \begin{equation}\label{eq:A_K upper bound B_K vareps}
        A_K(\varepsilon) = \frac{C_3}{\delta_1^{n_1}}e^{2\sqrt{|K|(n-1)}R(\varepsilon)}\left(\frac{R(\varepsilon)}{\varepsilon}\right)^{n_1} \le B_K\varepsilon^{-n'}
    \end{equation}
    with
    \begin{equation*}
        B_K
        \coloneqq
        C_3\left(\frac{C_4}{\delta_1}\right)^{n_1}
        \exp\left(\frac{2|K|(n-1)C_4^2}{\sigma} + \sigma\right)\left(\frac{n_1}{e\sigma}\right)^{n_1/2}.
    \end{equation*}
    Since $n_1=(n+n')/2$ and $\sigma=(n'-n)/2$, the geometric factor $C_3$ in $B_K$ depends only on $n,n',p$. The remaining factors depend on the transport parameters, the chosen Gaussian moment $\alpha,E_\alpha$, and explicitly on $K$ through $\sqrt{|K|(n-1)}$. In particular, there is no dependence on $\varepsilon$ hidden in $B_K$.\par
    Therefore, using \eqref{eq:A_K upper bound B_K vareps}, the condition \(
        A_K(\varepsilon) \le \tfrac{\Delta}{8}N\varepsilon^2,
    \), namely \eqref{eq:required on A_K},
    is implied by
    \begin{equation*}
        N \ge \frac{8B_K}{\Delta} \varepsilon^{-(n'+2)}.
    \end{equation*}
    For the condition $N\varepsilon^2\ge 4\log 2/\Delta$, namely \eqref{eq:second on N}, to be verified as well, since $0 < \varepsilon < 1$, it is enough to take
    \begin{equation*}
        N\ge N_0\,\varepsilon^{-(n'+2)}, \qquad N_0
        \coloneqq
        \max\left\{
            \frac{8B_K}{\Delta},
            \frac{4\log 2}{\Delta}
        \right\}.
    \end{equation*}
    \noindent \textbf{Case 2: $p=2$.}
Choose
    \begin{equation*}
        \lambda'<\lambda_3<\lambda_2<\lambda_1<\lambda,
        \qquad
        \Delta_2\coloneqq \lambda_3-\lambda'>0.
    \end{equation*}
    Recall that since $\nu$ satisfies $T_2(\lambda)$, it has a square-exponential moment of every order $\alpha<\lambda/2$. Set $\alpha_1\coloneqq \lambda_3/2$ for $\alpha_1 < \alpha < \lambda/2$ and
    \begin{equation*}
        \eta\coloneqq \sqrt2-1,
        \qquad
        \delta_1
        \coloneqq
        \frac{\eta}{4}
        \left(1-\sqrt{\frac{\lambda_2}{\lambda_1}}\right).
    \end{equation*}
    Then $0<\delta_1<\eta<1$ and $\lambda_2\eta^2<\lambda_1(\eta-\delta_1)^2$; the constant in \eqref{eq:H(nu', nu_R) lambda_2 p=2} is well-defined. 
    Arguing as in the case $1 \le p < 2$, for
    \begin{equation*}
        R(\varepsilon)
        \coloneqq
        \max\left\{
        R_0,\,
        R_2,\,
        1,\,
        \sqrt{
            \frac{1}{\alpha-\alpha_1}
            \log_+\left(
                \frac{16E_\alpha^2}{\Delta_2\eta^2\varepsilon^2}
            \right)
        },\,
        \sqrt{
            \frac{2}{\alpha}
            \log_+\left(
                \frac{16C_2^{(2)}}{\alpha e\Delta_2\eta^2\varepsilon^2}
            \right)
        }
        \right\},
    \end{equation*}
    since
    \begin{equation*}
        \frac{\alpha_1}{2}(1-\eta)^2\varepsilon^2 = \alpha_1\eta^2\varepsilon^2 = \frac{\lambda_3}{2}\eta^2\varepsilon^2
    \end{equation*}
    with $\eta = \sqrt{2} - 1$, then the right-hand side of \eqref{eq:before choice parameters p=2} is bounded by
    \begin{equation}\label{eq:both terms p=2}
        \exp\left(
            A_K(\varepsilon)
            -
            N\left(
                \frac{\lambda_2}{2}\eta^2\varepsilon^2
                -
                \frac{\Delta_2}{8}\eta^2\varepsilon^2
            \right)
        \right)
        +
        \exp\left(
            -N\left(
                \frac{\lambda_3}{2}\eta^2\varepsilon^2
                -
                \frac{\Delta_2}{8}\eta^2\varepsilon^2
            \right)
        \right).
    \end{equation}
    If in addition
    \begin{equation}\label{eq:required on A_K p=2}
        A_K(\varepsilon)
        \le
        \frac{\Delta_2}{8}\eta^2N\varepsilon^2,
    \end{equation}
    then
    \begin{equation*}
        \P\left(W_2(\widehat\nu^N,\nu)>\varepsilon\right)
        \le
        2\exp\left(
            -\frac{\lambda_3+\lambda'}{4}\eta^2N\varepsilon^2
        \right).
    \end{equation*}
    Whenever
    \begin{equation}\label{eq:second on N p=2}
        N\varepsilon^2 \ge
        \frac{4\log 2}{\Delta_2\eta^2},
    \end{equation}
    this yields
    \begin{equation}\label{eq:before substituting eta W2}
        \P\left(W_2(\widehat\nu^N,\nu)>\varepsilon\right)
        \le
        \exp\left(
            -\frac{\lambda'}{2}\eta^2N\varepsilon^2
        \right).
    \end{equation}
    Fix $n'>n$ and choose $n_1\coloneqq\tfrac12(n+n')$ and
    $\sigma\coloneqq n'-n_1=\tfrac12(n'-n)>0$.
    We choose $C_4^{(2)}$ such that
    \begin{equation}\label{eq:Regime R C_4 p=2}
        R(\varepsilon)
        \le
        C_4^{(2)}\sqrt{\log(e/\varepsilon)}
        \qquad\text{for all }0<\varepsilon<1.
    \end{equation}
    For instance, one may take
    \begin{equation*}
        C_4^{(2)}
        \coloneqq
        \max\left\{
            R_0,\,
            R_2,\,
            1,\,
            \sqrt{
                \frac{1}{\alpha-\alpha_1}
                \left(
                    2+\log_+\frac{16E_\alpha^2}{\Delta_2\eta^2}
                \right)
            },\,
            \sqrt{
                \frac{2}{\alpha}
                \left(
                    2+\log_+\frac{16C_2^{(2)}}{\alpha e\Delta_2\eta^2}
                \right)
            }
        \right\}.
    \end{equation*}
    Arguing as in Case 1, from \eqref{eq:Regime R C_4 p=2} we obtain
    \begin{equation}\label{eq:A_K upper bound B_K p=2}
        A_K(\varepsilon)
        \le
        B_K^{(2)}\varepsilon^{-n'},
    \end{equation}
    with 
    \begin{equation*}
        B_K^{(2)}
        \coloneqq
        C_3
        \left(\frac{C_4^{(2)}}{\delta_1}\right)^{n_1}
        \exp\left(
            \frac{2|K|(n-1)(C_4^{(2)})^2}{\sigma}
            +
            \sigma
        \right)
        \left(\frac{n_1}{e\sigma}\right)^{n_1/2}.
    \end{equation*}
    Here as well $n_1=(n+n')/2$ and $\sigma=(n'-n)/2$: the dependence on the auxiliary exponent $n'$ is entirely through these displayed constants and the factor $C_3=C_3(n,n',p)$.\par
    Thus, by \eqref{eq:A_K upper bound B_K p=2}, the condition \eqref{eq:required on A_K p=2} is implied by
    \begin{equation*}
        N
        \ge
        \frac{8B_K^{(2)}}{\Delta_2\eta^2}
        \varepsilon^{-(n'+2)}.
    \end{equation*}
    For the condition $N\varepsilon^2\ge 4\log 2/(\Delta_2\eta^2)$, namely \eqref{eq:second on N p=2}, to be verified as well, since $0<\varepsilon<1$, it is enough to take
    \begin{equation*}
        N\ge N_0^{(2)}\,\varepsilon^{-(n'+2)},
        \qquad
        N_0^{(2)}
        \coloneqq
        \max\left\{
            \frac{8B_K^{(2)}}{\Delta_2\eta^2},
            \frac{4\log 2}{\Delta_2\eta^2}
        \right\}.
    \end{equation*}
    Since $\eta^2=3-2\sqrt2$, \eqref{eq:before substituting eta W2} gives
    \begin{equation*}
        \P\left(W_2(\widehat\nu^N,\nu)>\varepsilon\right)
        \le
        \exp\left(
            -(3-2\sqrt2)\frac{\lambda'}{2}N\varepsilon^2
        \right),
    \end{equation*}
    which concludes the proof.
\end{proof}

\begin{proof}[Proof of \cref{thm:W_p quantitative concentration polynomial}]
    Since the proof follows the same line as the one of \cref{thm:W_p quantitative concentration}, we only highlight the differences.
    
    \noindent \textbf{Step 1 (Truncation).} Let $s\in[1,q/p)$ and set $\overline s\coloneqq \max\{s/2,1\}$. Since now only the polynomial moment $M_q \coloneqq \int_E d^q(o, x) d\nu(x)$ is finite with $p < q$, and instead of \eqref{eq:W_p^p nuR nu} we have
    \begin{equation}\label{eq:poly truncation measure}
        W_p^p\left(\nu_R,\nu\right) \le 2^pM_qR^{p-q}.
    \end{equation}
    Recall that
    \begin{equation*}
        W_p^p(\widehat\nu^N,\widehat\nu_R^N) \le \frac1N\sum_{i=1}^N Z^i,
        \qquad
        Z^i\coloneqq 2^pd^p(o,X^i)\mathbf{1}_{\{d(o,X^i)>R\}}.
    \end{equation*}
    Using
    \begin{equation*}
        \E\left[\left|\frac{1}{N}\sum_{i=1}^N(Z^i-\E\left[Z^i\right])\right|^s\right] \le C_s N^{\overline s-s}\E\left[(Z^1)^s\right],
        \qquad
        C_s\coloneqq
        \begin{cases}
            2^{s+1}, & 1\le s\le 2;\\
            (2s)^s, & s>2;
        \end{cases}
    \end{equation*}
    together with
    \begin{equation*}
        \E Z_1\le 2^pM_qR^{p-q},
        \qquad
        \E Z_1^s\le 2^{ps}M_qR^{ps-q},
    \end{equation*}
    whenever $t^p > 2^p M_q R^{p-q}$, Markov's inequality then yields
    \begin{align*}
        \P\left(W_p^p\left(\widehat\nu^N,\widehat\nu_R^N\right)>t^p\right)
        &\le
        \P\left(\abs{\frac{1}{N}\sum_{i=1}^N\left(Z_i - \E\left[Z^i\right]\right)}^s > \abs{t^p - \E\left[Z^1\right]}^s\right) \\
        &\le
        C_s2^{ps}M_qN^{\overline s-s}
        \frac{R^{ps-q}}{(t^p-2^pM_qR^{p-q})^s}. \numberthis \label{eq:empirical truncation polynomial}
    \end{align*}
    For $R \ge 1$ large enough so that $W_p(\nu_R,\nu) \le \varepsilon/8$, using the triangle inequality and \eqref{eq:empirical truncation polynomial} we obtain
    \begin{align*}
        \P\left(W_p(\widehat\nu^N,\nu)>\varepsilon\right)
        &\le
        \P\left(W_p(\widehat\nu^N,\widehat\nu_R^N)>\frac{\varepsilon}{2}\right) +
        \P\left(W_p(\widehat\nu_R^N,\nu_R)>\frac{3\varepsilon}{8}\right) \\ &\quad +
        \P\left(W_p(\nu_R, \nu)>\frac{\varepsilon}{8}\right)
        \\
        &\le \P\left(W_p(\widehat\nu_R^N,\nu_R)>\frac{3\varepsilon}{8}\right) + C_s2^{ps}M_qN^{\overline s-s}
        \frac{R^{ps-q}}{\big(\!\left(\tfrac{\varepsilon}{2}\right)^p-2^pM_qR^{p-q}\big)^s}. \numberthis \label{eq:Wp nun nu vareps triangle step 1}
    \end{align*}
    Letting $0 < \varepsilon < 1$,
    \begin{equation}\label{eq:condition on R}
        \begin{gathered}
            R\ge R_0\varepsilon^{-p/(q-p)},\\
            R_0\coloneqq
            \max\left\{
                \big(16M_q^{1/p}\big)^{p/(q-p)},
                \big(2^{2p+1}M_q\big)^{1/(q-p)},
                1
            \right\}.
        \end{gathered}
\end{equation}
    not only makes $W_p(\nu_R, \nu) \le \varepsilon/8$ hold recalling \eqref{eq:poly truncation measure}, but also
    \begin{equation*}
        \left(\frac{\varepsilon}{2}\right)^p-2^pM_qR^{p-q}
        \ge 2^{-(p+1)}\varepsilon^p,
    \end{equation*}
    and \eqref{eq:Wp nun nu vareps triangle step 1} becomes
    \begin{equation}\label{eq:Wp nun nu vareps triangle step 1 end}
        \begin{aligned}
            \P\left(W_p(\widehat\nu^N,\nu)>\varepsilon\right)
            \le{}&\P\left(W_p(\widehat\nu_R^N,\nu_R)>
                \frac{3\varepsilon}{8}\right)\\
            &+C_1N^{\overline s-s}\varepsilon^{-ps}.
        \end{aligned}
\end{equation}
    where $C_1\coloneqq C_s2^{s(2p+1)}M_q$.

    \noindent \textbf{Step 2 (Covering number estimates).}
    We apply the quantitative Sanov estimate \eqref{eq:quantitative Sanov A},
    \begin{equation}\label{eq:quantitative Sanov A second instance}
        \begin{aligned}
            \P\left(\widehat\nu_R^N\in\mathcal A\right)
            &\le \left(\frac{16eR}{\delta'}\right)^{{pN(B_R(o),\delta'/2)}}\\
            &\quad\times\exp\left(-N
                \inf_{\nu'\in\mathcal A_{\delta'}}H(\nu'\mid\nu_R)\right).
        \end{aligned}
\end{equation}
    to
    \begin{equation*}
        \mathcal A\coloneqq
        \left\{
            \nu'\in\mathcal P(B_R(o))\::\:
            W_p(\nu', \nu_R) > \frac{3\varepsilon}{8}\right\}, \qquad \delta' \coloneqq \frac{\varepsilon}{8}.
    \end{equation*}
    If $\nu'\in\mathcal A_{\delta'}$, then $W_p(\nu', \nu_R)\ge \varepsilon/4$.
    Note that $\nu_R$ satisfies a penalized $T_p$ inequality thanks to the Csisz\'ar--Kullback--Pinsker inequality \cite[Particular case 2.5]{bolley_weighted_2005};
    \begin{equation*}
        W_p^{2p}\left(\nu',\nu_R\right) \le 2^{2p-1}R^{2p}H(\nu'\mid\nu_R), \qquad \forall \nu' \in \mathcal{P}(B_R(o)).
    \end{equation*}
    Thus for $\nu' \in \mathcal{A}_{\delta'}$,
    \begin{equation}\label{eq:Csizar-Kullback-Pinsker inequality}
        H(\nu'\mid\nu_R)
        \ge 2^{1-6p}\left(\frac\varepsilon R\right)^{2p}.
\end{equation}
    Let $n_1 > n$ and recall the estimate \eqref{eq:metric entropy bound final},
    \begin{equation}\label{eq:metric entropy bound final poly}
        \begin{aligned}
            pN\left(B_R(o),\frac{\delta'}2\right)
            \log\left(\frac{16eR}{\delta'}\right)
            \le{}&C_2\exp\left(2\sqrt{|K|(n-1)}\,R\right)
                \left(\frac{R}{\delta'}\right)^{n_1}.
        \end{aligned}
\end{equation}
    with
    \begin{equation*}
        C_2 \coloneqq \frac{p\,10^n(16e)^{n_1-n}}{e(n_1-n)}.
    \end{equation*}
    Combining \eqref{eq:quantitative Sanov A second instance}, \eqref{eq:Csizar-Kullback-Pinsker inequality}, and \eqref{eq:metric entropy bound final poly}, we obtain
    \begin{equation}\label{eq:compact polynomial}
        \P\left(W_p(\widehat\nu_R^N,\nu_R)>\frac{3\varepsilon}{8}\right)
        \le
        \exp\left(
            C_2 8^{n_1}e^{2\sqrt{|K|(n-1)}R}
            \left(\frac{R}{\varepsilon}\right)^{n_1}
            -
            2^{1-6p}N\left(\frac{\varepsilon}{R}\right)^{2p}
        \right).
    \end{equation}
    
    \noindent \textbf{Step 3 (Choice of parameters).}
    Fix $n'\coloneqq n_1>n$ and recall the condition \eqref{eq:condition on R} on $R$; set 
    \begin{equation*}
        R \coloneqq R_0 \varepsilon^{-p/(q-p)}.
    \end{equation*}
   We can estimate the exponential term in \eqref{eq:compact polynomial} as
\begin{align*}
    \exp\left(
        C_2 8^{n'}e^{2\sqrt{|K|(n-1)}R}
        \left(\frac{R}{\varepsilon}\right)^{n'}
        -
        2^{1-6p}N\left(\frac{\varepsilon}{R}\right)^{2p}
    \right) 
    &\le
    \exp\left(
        -2^{-6p}N\left(\frac{\varepsilon}{R}\right)^{2p}
    \right) \\
    &\le
    \exp\left(
        -C_3 N \varepsilon^{2p\left(1+p/(q-p)\right)}
    \right),
\end{align*}
where \(C_3\coloneqq (8R_0)^{-2p}\), provided
\begin{equation}\label{eq:N regime polynomial}
    N
    \ge
    N_0'
    \varepsilon^{-(n'+2p)\left(1+p/(q-p)\right)}
    \exp\left(2\sqrt{|K|(n-1)}R\right),
    \qquad
    N_0'
    \coloneqq
    2^{6p}C_2 8^{n'}R_0^{n'+2p}.
\end{equation}
    Combining with \eqref{eq:compact polynomial}, we obtain
    \begin{equation}\label{eq:mixed after compact polynomial}
        \P\left(W_p(\widehat\nu^N, \nu)>\varepsilon\right)
        \le
        \exp\left(-C_3 N\varepsilon^{2p\left(1+p/(q-p)\right)}\right) + C_1 N^{\overline{s}-s}\varepsilon^{-ps}.
    \end{equation}
    Now set
    \begin{equation*}
        s\coloneqq \frac qp-\frac{\delta}{2}.
    \end{equation*}
    For the case (i), $1\le p<q/2$, since $s>2$ and $\overline s=s/2$, then
    \begin{equation}\label{eq:truncation absorbed case i}
        C_1N^{\overline s-s}\varepsilon^{-ps}
        \le
        \frac12\,\varepsilon^{-q}N^{-q/(2p)+\delta/2}.
    \end{equation}
    if in addition $N_0\ge (2C_1)^{4/\delta}$.
    For the case (ii), $q/2\le p<q$, since $1<s<2$ and $\overline s=1$, then
    \begin{equation}\label{eq:truncation absorbed case ii}
        C_1N^{\overline s-s}\varepsilon^{-ps}
        \le
        \frac12\,\varepsilon^{-q}N^{1-q/p+\delta}.
    \end{equation}
    provided the same additional condition $N_0\ge (2C_1)^{4/\delta}$.
    
    It remains to dominate the exponential term of \eqref{eq:mixed after compact polynomial}. We claim that, under \eqref{eq:N regime polynomial}, after possibly enlarging $N_0'$,
    \begin{equation}\label{eq:compact exp dominated}
        \exp\left(-C_3N\varepsilon^{2p\left(1+p/(q-p)\right)}\right)
        \le
        \frac12\,\varepsilon^{-q}N^{-\beta}, \qquad \beta\coloneqq
        \begin{cases}
            \frac{q}{2p}-\frac{\delta}{2}, & \text{case (i)};\\
            \frac qp-1-\delta, & \text{case (ii)}.
        \end{cases}
    \end{equation}
    Set
    \begin{equation*}
        a\coloneqq (n'+2p)\left(1+\frac{p}{q-p}\right)=\frac{q(n'+2p)}{q-p},
        \qquad \gamma\coloneqq 2p\left(1+\frac{p}{q-p}\right)=\frac{2pq}{q-p}.
    \end{equation*}
    A possible choice is given by
    \begin{equation}\label{eq:N0 compact domination}
        N_0''
        \coloneqq
        \max\left\{
            \frac{4\beta}{C_3},
            \frac{4}{C_3}\left[ \log 2+\frac{c}{e} + \beta\log_+\left(\frac{4\beta}{C_3}\right) \right]
        \right\}, \qquad
        c\coloneqq \frac{(a\beta-q)_+}{a-\gamma}.
    \end{equation}
    Then $N_0''\ge N_0'$ implies \eqref{eq:compact exp dominated} for all $N \ge N_0''$. To see this, observe first that for $0 < \varepsilon < 1$,
    $N\ge N_0'\varepsilon^{-a}$, and $N_0'\ge \beta/C_3$, it holds
    \begin{equation*}
        N \ge \frac{\beta}{C_3}\varepsilon^{-\gamma}.
    \end{equation*}
    Therefore, the function $k\mapsto e^{-C_3k\varepsilon^\gamma}k^\beta$ is nonincreasing for $k \ge N$, and it is enough to prove \eqref{eq:compact exp dominated} at $N=N_0''\varepsilon^{-a}$. This reduces to
    \begin{equation*}
        e^{-C_3N_0'' \varepsilon^{-(a-\gamma)}}
        \le
        \frac12 (N_0'')^{-\beta}\varepsilon^{-(q-a\beta)}.
    \end{equation*}
    Equivalently, taking logarithms, it is enough that
    \begin{equation*}
        C_3N_0''\varepsilon^{-(a-\gamma)}
        \ge \log 2+\beta\log N_0''+ c\log \left(\varepsilon^{-(a-\gamma)}\right),
    \end{equation*}
    but $\log (\varepsilon^{-(a-\gamma)}) \le \varepsilon^{-(a-\gamma)}/e$ for $\varepsilon \le 1$, and since \eqref{eq:N0 compact domination} implies
    \begin{equation*}
        C_3N_0'' \ge \log 2+ \beta\log N_0'' + \frac{c}{e},
    \end{equation*}
    the claim follows.

    We thus require
    \begin{equation}\label{eq:lower bound N final poly}
        N \ge N_0 \varepsilon^{-\frac{q(n' + 2p)}{q-p}}\exp\left(2\sqrt{|K|(n-1)}R\right), \qquad N_0
        \coloneqq
        \max\left\{
            N_0', N_0'', (1\vee 2C_1)^{4/\delta}\right\}.
    \end{equation}
    Combining \eqref{eq:Wp nun nu vareps triangle step 1 end} and \eqref{eq:mixed after compact polynomial}--\eqref{eq:compact exp dominated} under \eqref{eq:lower bound N final poly}, we obtain in case (i)
    \begin{equation*}
        \P\left(W_p(\widehat\nu^N,\nu)>\varepsilon\right)
        \le
        \varepsilon^{-q}
        N^{-q/(2p)+\delta/2},
    \end{equation*}
    while in case (ii)
    \begin{equation*}
        \P\left(W_p(\widehat\nu^N,\nu)>\varepsilon\right)
        \le
        \varepsilon^{-q}
        N^{1-q/p+\delta}.
    \end{equation*}
    This concludes the proof.
\end{proof}

\section{An interacting particle system}\label{sec:An interacting particle system}

\subsection{Description of the interacting particle system}

Let \((M,g)\) be a complete, connected Riemannian manifold of dimension \(m\ge2\).
Let \(V:M\to\mathbb R\) be a confinement potential and
\(W:M\times M\to\mathbb R\) an interaction potential. For \(N\ge1\), we consider
the McKean--Vlasov interacting particle system
\begin{equation}\label{eq:SDE interacting system particles}
    dX_t^i
    =
    \sqrt2\,dB_t^{i,M}
    -
    \nabla V(X_t^i)\,dt
    -
    \frac1N\sum_{j=1}^N
    \nabla_xW(X_t^i,X_t^j)\,dt,
    \qquad i=1,\dots,N,
\end{equation}
where \(B^{1,M},\dots,B^{N,M}\) are independent Brownian motions on \(M\).
We denote by
\[
    \widehat\mu_t^N
    \coloneqq
    \frac1N\sum_{i=1}^N\delta_{X_t^i}
\]
the associated empirical measure.

Formally, as \(N\to\infty\), the empirical measure \(\widehat\mu_t^N\)
converges to a deterministic probability measure \(\mu_t\). The limiting
equation is obtained by replacing the empirical interaction in
\eqref{eq:SDE interacting system particles} by its average against \(\mu_t\).
Thus \(\mu_t\) solves the nonlinear Fokker--Planck equation, also called the
aggregation--diffusion equation,
\begin{equation}\label{eq:aggregation diffusion}
    \partial_t\mu_t
    =
    \Delta\mu_t
    +
    \nabla\cdot
    \left(
        \mu_t
        \left[
            \nabla V+\nabla(W*\mu_t)
        \right]
    \right),
\end{equation}
where
\[
    W*\mu_t(x)
    \coloneqq
    \int_M W(x,y)\dd\mu_t(y),
    \qquad
    \nabla(W*\mu_t)(x)
    =
    \int_M \nabla_xW(x,y)\dd\mu_t(y).
\]
This is the usual propagation-of-chaos picture for McKean--Vlasov diffusions.
Under the assumptions stated below, we construct the nonlinear law for
\eqref{eq:aggregation diffusion} without assuming a smooth initial
probability density. We then derive concentration estimates for its
approximation by interacting empirical measures.

As in Bolley--Guillin--Villani \cite{bolley_quantitative_2007}, we introduce a
noninteracting proxy system \((Y_t^i)_{1\le i\le N}\). Its particles are
independent and have common law \(\mu_t\):
\begin{equation}\label{eq:proxy noninteracting system}
    dY_t^i
    =
    \sqrt2\,d\overline B_t^{i,M}
    -
    \nabla V(Y_t^i)\,dt
    -
    \int_M \nabla_xW(Y_t^i,z)\dd\mu_t(z)\,dt.
\end{equation}
We write
\[
    \widehat\nu_t^N
    \coloneqq
    \frac1N\sum_{i=1}^N\delta_{Y_t^i}.
\]
Since \(\widehat\nu_t^N\) is the empirical measure of independent variables
with common law \(\mu_t\), the concentration estimate of
\cref{thm:W_p quantitative concentration} can be applied to it once suitable
moment estimates are propagated along \eqref{eq:aggregation diffusion}.

The proof has two parts. First, we prove square-exponential moment bounds and
time-increment estimates for the proxy process \(Y_t\). These estimates yield
concentration estimates for
\[
    \sup_{0\le t\le T}W_p(\widehat\nu_t^N,\mu_t).
\]
Second, we couple the interacting system \(X_t^i\) with the proxy system
\(Y_t^i\). The coupling is constructed by parallel displacement of the Brownian
noises. For this reason, it is convenient to write the stochastic equations on
the orthonormal frame bundle \(O(M)\). The coupling estimate then reduces the
control of \(W_p(\widehat\mu_t^N,\mu_t)\) to the concentration of the
independent proxy empirical measure \(\widehat\nu_t^N\).

\subsection{Assumptions on the potentials}

Fix \(o\in M\) and write \(r(x)\coloneqq d(o,x)\). We
assume that there exist constants
\[
    \beta,\widetilde\beta>0,
    \qquad
    A_V,A_W \ge 0, \qquad L_V, L_{W,1}, L_{W,2}>0,
\]
such that, for all \(x,y,z\in M\),
\begin{equation}\label{eq:assumptions V IPS}
    \Hess V\ge \beta g,
    \qquad
    |\nabla V(x)|\le A_V+L_Vr(x),
\end{equation}
and
\begin{equation}\label{eq:assumptions W growth IPS}
    \Hess_xW(\cdot,y)\ge \widetilde\beta g,
    \qquad
    |\nabla_xW(x,y)|
    \le
    A_W+L_{W,1}\big(d(o,x) + d(o,y)\big).
\end{equation}
We also assume that the dependence of \(\nabla_xW\) on the second variable is
Lipschitz, uniformly in the first variable; i.e., for every \(x,y,z\in M\) and every
unit vector \(u\in T_xM\),
\begin{equation}\label{eq:assumption W second variable Lipschitz IPS}
    \left|
        \big(
            \nabla_xW(x,y)-\nabla_xW(x,z)\big) \cdot
            u
    \right|
    \le
    L_{W,2}d(y,z).
\end{equation}

The lower Hessian bound in the first variable gives dissipation, the linear
growth assumption controls the moment estimates, and the Lipschitz assumption in
the second variable is used in the coupling estimate through Kantorovich
duality. Throughout the sequel, we write
\[
    \Lambda\coloneqq\beta+\widetilde\beta.
\]

As an example, on a Hadamard manifold, the quadratic potentials
\[
    V(x)=\frac{\ell}{2}d^2(o,x),\qquad
    W(x,y)=\frac{\ell}{2}d^2(x,y),\qquad \ell>0,
\]
satisfy \eqref{eq:assumptions V IPS}--
\eqref{eq:assumption W second variable Lipschitz IPS}.
Indeed, Hessian comparison \cite[Theorem 11.7]{lee_introduction_2018}
gives $\Hess_x(d^2(\cdot,y)/2)\ge g$, while the metric comparison
theorem \cite[Theorem 11.10]{lee_introduction_2018} implies that
$y\mapsto\log_x(y)$ is $1$-Lipschitz for every fixed $x\in M$.
Since $\nabla_xW(x,y)=-\ell\log_x(y)$, the growth and Lipschitz
assumptions hold with $L_{W,1}=L_{W,2}=\ell$.

This example also illustrates a difference with the recent
work of Chae--Cho--Jang--Yang
\cite{chae_cho_jang_yang_2026} whose interaction
assumptions impose a uniform control on the variation
of $\nabla_xW$ in the first variable, which is not
required here. Indeed, on the hyperbolic space of
constant sectional curvature $-1$, the tangential
eigenvalues of $\Hess_xW(\cdot,y)$ are
$\ell r\coth r$, where $r=d(x,y)$, and are
therefore unbounded as $r\to\infty$.
Consequently, the squared-distance interaction
does not satisfy their Hypothesis 3.2(ii),
although it satisfies our assumptions.

\subsection{Construction of the nonlinear law}

The assumptions above also give a direct construction of the McKean--Vlasov law. We record it here so that the moment computations below do not rely on assuming a smooth solution of \eqref{eq:aggregation diffusion}.

\begin{proposition}[Well-posedness of the nonlinear diffusion]
\label{prop:wellposed-mv}
    Let $(M,g)$ be a complete Riemannian manifold with
    $\Ric\ge Kg$ for some $K\le0$. Suppose that
    \eqref{eq:assumptions V IPS}--
    \eqref{eq:assumption W second variable Lipschitz IPS}
    hold. For every $\mu_0\in\mathcal P_2(M)$, there exists
    a unique continuous law curve
    \[
        \mu\in C([0,\infty);\mathcal P_2(M))
    \]
    of a nonexplosive McKean--Vlasov diffusion starting
    from $\mu_0$. It solves \eqref{eq:aggregation diffusion}
    against compactly supported smooth test functions.
\end{proposition}

The proof is given in \cref{app:wellposed-mv}, using the radial
Lyapunov and parallel-coupling estimates established below.

\subsection{Preliminary results}\label{sec:preliminary results}
We first establish propagation of square-exponential moments for \eqref{eq:aggregation diffusion}; these imply a $T_1$ inequality at each time, with a common constant on every finite interval. The radial cutoff lemma below justifies the integrations by parts in the presence of the cut locus; its proof is deferred to \cref{app:radial-cutoff}. 
While the standard PDE approach presented here requires
these justifications, one could instead work directly at
the SDE level using Kendall's formula, as we do later.
We nevertheless believe that the PDE approach is more
readable at first, once these technicalities have been addressed.

\begin{lemma}[Radial cutoff and integration by parts]\label{lem:radial-cutoff}
    Let $(\mu_t)_{0\le t\le T}$ be the law from
    \cref{prop:wellposed-mv}, and set
    \[
        r(x)\coloneqq d(o,x),\qquad
        L_t\coloneqq\Delta+b_t\cdot\nabla,\qquad
        b_t\coloneqq-\nabla V-\nabla(W*\mu_t).
    \]
    Let $\Psi_t=\Phi(t,r^2)$ be either $r^2$ or
    $e^{\alpha(t)r^2}$, where $\alpha>0$ is continuous and
    piecewise $C^1$ on $[0,T]$.
    Choose smooth, increasing, concave functions
    $\chi_j:[0,\infty)\to[0,\infty)$ satisfying
    \[
        \begin{gathered}
            \chi_j(s)=s\quad(s\le j),\qquad
            \chi_j'(s)=0\quad(s\ge2j),\\
            0\le\chi_j'\le1,\qquad\chi_j''\le0,
            \qquad\chi_j(s)\uparrow s.
        \end{gathered}
    \]
    Then, for $0\le s\le t\le T$,
    \begin{equation}\label{eq:localized-generator}
        \begin{aligned}
            &\int_M\chi_j(\Psi_t)\dd\mu_t
                -\int_M\chi_j(\Psi_s)\dd\mu_s\\
            &\qquad\le\int_s^t\!\int_M
                \chi_j'(\Psi_u)
                \bigl(\partial_u\Psi_u+(L_u\Psi_u)_{\mathrm{reg}}\bigr)
                \dd\mu_u\dd u.
        \end{aligned}
    \end{equation}
    Here $(L_u\Psi_u)_{\mathrm{reg}}$ denotes the classical
    generator evaluated outside $\{o\}\cup\operatorname{Cut}(o)$.
\end{lemma}

\begin{lemma}[Propagation of square-exponential moments]\label{thm:propagation of moments mu_t}
    Let $(M, g)$ be a complete noncompact Riemannian manifold of dimension $m \ge 2$ satisfying $\Ric\ge Kg$ for some $K \le 0$. Let $\mu_t$ solve the aggregation--diffusion equation \eqref{eq:aggregation diffusion} with potentials $V, W$ satisfying respectively \eqref{eq:assumptions V IPS} and \eqref{eq:assumptions W growth IPS}.
    Assume that $\mu_0$ has finite square-exponential moment; i.e., there is a constant $\alpha(0) > 0$ such that
    \begin{equation*}
        E_{\alpha(0)} \coloneqq \int_M e^{\alpha(0) d^2(o, x)}\dd\mu_0(x) < +\infty.
    \end{equation*}
    Then for every $t > 0$, there exists a continuous function $\alpha(t) > 0$ such that
    \begin{equation*}
        E_{\alpha(t)} \coloneqq \int_M e^{\alpha(t) d^2(o, x)}\dd\mu_t(x) < +\infty.
    \end{equation*}
    If in addition $\beta + \widetilde\beta > L_{W,1}$, then $\sup_{t \ge 0} E_{\alpha(t)} < +\infty$.
\end{lemma}
\begin{proof}
    \noindent \textbf{(Step 1) Propagation of second moment.}
    Set $r(x)\coloneqq d(o,x)$ and
    \begin{equation*}
        M_2(t)\coloneqq \int_M r^2(x) \dd\mu_t(x),
    \end{equation*}
    A formal integration by parts with $r^2$ yields
    \begin{align*}
        \frac{\dd}{\dd t}M_2(t)
        ={}&\int_M \Delta r^2(x)\dd\mu_t(x)
        -\int_M \nabla r^2(x)\cdot\nabla V(x)\dd\mu_t(x)\\
        &-\int_{M\times M}
        \nabla r^2(x)\cdot\nabla_xW(x,y)
        \dd\mu_t(y)\dd\mu_t(x).
    \end{align*}
    However, this identity is not directly justified since
    $r^2$ is unbounded and not smooth at the cut locus. Instead, we apply \cref{lem:radial-cutoff} to $\Psi=r^2$, and only remove the cutoffs after obtaining the bounds below. In particular, the singular part of the distributional Laplacian contributes with a nonpositive sign.

    By Laplacian comparison theory \cite[Theorem 11.15]{lee_introduction_2018}, outside of $\{o\} \cup \operatorname{Cut}(o)$,
    \begin{equation*}
        \Delta r(x) \le (m-1)\sqrt{\frac{|K|}{m-1}}\coth \left( \sqrt{\frac{|K|}{m-1}} r(x)\right) \le \frac{m-1}{r(x)} + \sqrt{|K|(m-1)},
    \end{equation*}
    where in the last inequality we used $\coth u \le 1 + 1/u$. Then, since $|\nabla r(x)| = 1$,
    \begin{equation}\label{eq:comparison Laplacian}
        \begin{aligned}
            \Delta r^2(x)
            &=2r(x)\Delta r(x)+2|\nabla r(x)|^2\\
            &\le2m+2\sqrt{|K|(m-1)}\,r(x).
        \end{aligned}
    \end{equation}
    Recall we assumed that
    \begin{equation*}
        \Hess V\ge\beta g,\qquad
        \Hess_xW(\cdot,y)\ge\widetilde\beta g, \qquad \abs{\nabla_x W(o,y)}\le A_W+L_{W,1}r(y).
    \end{equation*}
    With this notation, if $\gamma : [0,r(x)]\to M$ is a unit-speed minimizing geodesic from $o$ to $x$, then outside of $\{o\} \cup \operatorname{Cut}(o)$,
\begin{align*}
    \nabla r(x)\cdot\nabla V(x)
    &=
    \frac{\dd}{\dd s}\bigg|_{s=r(x)} V(\gamma(s)) \\
    &=
    \frac{\dd}{\dd s}\bigg|_{s=0} V(\gamma(s))
    +
    \int_0^{r(x)}
    \frac{\dd^2}{\dd s^2}V(\gamma(s))\dd s \\
    &=
    \langle \nabla V(o),\dot\gamma(0)\rangle
    +
    \int_0^{r(x)}
    \Hess V_{\gamma(s)}(\dot\gamma(s),\dot\gamma(s))\dd s \\
    &\ge
    -|\nabla V(o)|+\beta r(x).
\end{align*}
    Therefore, for every $x \in M\setminus(\{o\} \cup \operatorname{Cut}(o))$, it holds
    \begin{equation}\label{eq:r2 V upper}
        \nabla r^2(x)\cdot\nabla V(x) \ge 2\beta r^2(x)-2|\nabla V(o)|r(x).
    \end{equation}
    Similarly, for every $x \in M\setminus(\{o\} \cup \operatorname{Cut}(o))$ and every $y \in M$,
    \begin{equation}\label{eq:r2 W upper}
        \nabla r^2(x)\cdot\nabla_x W(x,y) \ge 2\widetilde\beta r^2(x)-2|\nabla_xW(o,y)|r(x).
    \end{equation}
    Combining \eqref{eq:comparison Laplacian}--\eqref{eq:r2 V upper}, and integrating \eqref{eq:r2 W upper} in $d\mu_t(y)$, the localized inequality of \cref{lem:radial-cutoff} gives, after removing the cutoffs,
    \begin{equation}\label{eq:d dt M_2 before Gronwall}
        \frac{\dd}{\dd t}M_2(t) \le 2m + 2C_1 \sqrt{M_2(t)} - 2\left(\beta + \widetilde\beta - L_{W,1}\right)M_2(t)
    \end{equation}
    with
    \[C_1 \coloneqq \sqrt{|K|(m-1)} + |\nabla V(o)| + A_W.\]
    Here the passage from $\chi_j(r^2)$ to $r^2$ uses that $\mu\in C([0,T];\mathcal P_2(M))$, so the drift upper bound has at most quadratic growth and is integrable on finite time intervals. The concavity term and the cut-locus term have already been discarded in the localized estimate.\par
    By Gronwall's inequality, given $T > 0$, from \eqref{eq:d dt M_2 before Gronwall} we conclude
    \begin{equation*}
        \sup_{t \in [0, T]} M_2(t) \coloneqq M_{2,T} < +\infty,
    \end{equation*}
    and if in addition $\Lambda>L_{W,1}$, then $M_{2,T} < C_{\mathrm{M}_2}$ for some constant $C_{\mathrm{M}_2} > 0$, independently of $T$.

    \noindent \textbf{Step 2 (Propagation of square-exponential moments).}
    We again first use $\chi_j(\Psi_t)$ with $\Psi_t=e^{\alpha(t)r^2}$ and apply \cref{lem:radial-cutoff}. The following differential relations are first established for the cutoff integrals, and the estimates obtained from them are independent of $j$. The limit $j\to\infty$ is then justified by monotone convergence.\par
    Let
    \begin{equation*}
        E_{\alpha(t)}\coloneqq \int_M e^{\alpha(t)r^2(x)}\dd\mu_t(x),
    \end{equation*}
    where $\alpha(t)>0$ will be chosen below. Since outside of $\{o\} \cup \operatorname{Cut}(o)$,
    \begin{equation*}
        \nabla e^{\alpha r^2}=\alpha e^{\alpha r^2}\nabla r^2,
        \qquad
        \Delta e^{\alpha r^2}
        =
        e^{\alpha r^2}\left(\alpha\Delta r^2+\alpha^2|\nabla r^2|^2\right),
        \qquad
        |\nabla r^2|^2=4r^2,
    \end{equation*}
    we obtain, using \eqref{eq:comparison Laplacian}--\eqref{eq:r2 W upper},
    \begin{align*}
        \frac{d}{d t}E_{\alpha(t)} &\le
        \int_M e^{\alpha(t)r^2} \Big[ 2m\alpha(t) + \left(\frac{d}{d t}\alpha(t)+4\alpha(t)^2-2\Lambda\alpha(t)\right)r^2 +2\alpha(t)C_2r\Big]\dd\mu_t \numberthis\label{eq:d dt E alpha diff} \\
        &\le (2m + C_2^2)\alpha(t) E_{\alpha(t)} + \left(\frac{\dd}{\dd t}\alpha(t)+4\alpha(t)^2+\left(1-2\Lambda\right)\alpha(t)\right)\int_M r^2 e^{\alpha r^2}\dd\mu_t, \numberthis\label{eq:d dt E alpha diff 2}
    \end{align*}
    where we used $2\alpha C_2 r \le \alpha C_2^2 + \alpha r^2$, with
    \begin{equation*}
        C_2\coloneqq \sqrt{|K|(m-1)}+|\nabla V(o)|+A_W+L_{W,1}\sqrt{M_{2,T}}.
    \end{equation*}
    We choose $\alpha$ as the solution of
    \begin{equation*}
        \frac{d}{d t}\alpha(t)+4\alpha(t)^2+\left(1-2\Lambda\right)\alpha(t)=0,
        \qquad
        \alpha(0)=\alpha_0,
    \end{equation*}
    so that the term in brackets in \eqref{eq:d dt E alpha diff 2} vanishes. Since $\alpha(t)>0$ on $[0, T]$, \eqref{eq:d dt E alpha diff 2} becomes
    \begin{equation*}
        \frac{\dd}{\dd t}E_{\alpha(t)}(t)
        \le
        \alpha(t)(2m+C_2^2)E_{\alpha(t)}(t).
    \end{equation*}
    Hence, for every $T>0$, Gronwall's inequality yields
    \begin{equation*}
        \sup_{0\le t\le T}E_{\alpha(t)}(t)
        \le
        E_{\alpha_0}(0)
        \exp\left((2m+C_2^2)\int_0^T\alpha(s)\dd s\right)<+\infty.
    \end{equation*}

    Assume now furthermore that $\Lambda > L_{W,1}$ and set
    \begin{equation*}
        C_2^\infty
        \coloneqq
        \sqrt{|K|(m-1)}+|\nabla V(o)|+A_W+L_{W,1} C_{\mathrm{M}_2}^{1/2}.
    \end{equation*}
    In this case, we choose $\alpha(t)$ as the solution of
    \begin{equation*}
        \frac{d}{dt}\alpha(t)=\Lambda\alpha(t)-4\alpha(t)^2,
        \qquad
        \alpha(0)=\alpha_0
    \end{equation*}
     so that
     \begin{equation*}
         \alpha_{-} \coloneqq \min\left\{\alpha_0, \frac{\Lambda}{4}\right\} \le \alpha(t) \le \max\left\{\alpha_0, \frac{\Lambda}{4}\right\} \eqqcolon \alpha_+.
     \end{equation*}
    Using $2C_2^{\infty} r\le \tfrac{\Lambda}{2}r^2+\tfrac{2(C_2^\infty)^2}{\Lambda}$ in the differential inequality \eqref{eq:d dt E alpha diff} gives
    \begin{equation*}
        \frac{\dd}{\dd t}E_{\alpha(t)}
        \le
        \alpha(t)\int_M
        \left(C_3-\frac{\Lambda}{2}r^2\right)e^{\alpha(t)r^2}\dd\mu_t,
        \qquad
        C_3\coloneqq 2m+\frac{2(C_2^\infty)^2}{\Lambda}.
    \end{equation*}
    Splitting the integral over $\{r\le R_*\}$ and $\{r>R_*\}$ for $R_*^2\coloneqq 4C_3/\Lambda$ yields
    \begin{equation*}
        \frac{\dd}{\dd t}E_{\alpha(t)}
        \le
        2C_3\alpha_+e^{\alpha_+R_*^2}
        -
        C_3\alpha_- E_{\alpha(t)},
    \end{equation*}
    and thus $\sup_{t\ge0}E_{\alpha(t)}<+\infty$. In this last step, we first use the finite-time exponential estimate, pass $j\to\infty$ in the localized inequality with its favorable negative quadratic term, and then apply Gronwall.
\end{proof}

We now use the propagated square-exponential moments to estimate the increments of the noninteracting proxy on finite time intervals. In this auxiliary system, the empirical mean-field interaction is replaced by the limiting law \(\mu_t\). Equivalently, the dynamics are written in Stratonovich form on the orthonormal frame bundle \(\pi:O(M)\to M\) as
\begin{equation}\label{eq:linear SDE}
dU_t^i
=
\sqrt{2}\sum_{k=1}^m H_k(U_t^i)\circ dB_t^{i,k}
+
\widetilde b_t(U_t^i)dt,
\qquad
Y_t^i=\pi(U_t^i), \qquad i = 1, \dots, N,
\end{equation}
where \(B_t=(B_t^1,\dots,B_t^m)\) is a standard Brownian motion on \(\mathbb R^m\), \(\widetilde b_t\) denotes the horizontal lift of the time-dependent vector field
\begin{equation*}
b_t(x)
\coloneqq
-\nabla V(x)-\nabla(W*\mu_t)(x),
\end{equation*}
and \(H_k\) is the fundamental horizontal vector field associated with the canonical basis vector \(e_k\in\mathbb R^m\); see \cite[Chapters 3--4]{hsu_stochastic_2002} for the notion of horizontal lift to the orthonormal frame bundle and the monograph of Hsu \cite{hsu_stochastic_2002}.
We obtain that the empirical measure of $(Y^i_t)_{1\le i\le N}$ given by
\begin{equation*}
    \widehat\nu_t^N \coloneqq \frac{1}{N} \sum_{i=1}^N \delta_{Y^i_t}
\end{equation*}
satisfies the same quantitative estimate as \cite[Proposition 5.1]{bolley_quantitative_2007}, under some additional assumptions on the potentials $V$ and $W$.
\begin{proposition}[Time increments of the independent empirical measure]\label{prop:time-increment-proxy}
    In the setting of \cref{thm:propagation of moments mu_t}, for every $T > 0$, there are constants $C_T, a_T > 0$ such that, for every $1\le p\le2$, $t_0\in[0,T-\delta]$, $0<\delta<T$, and $\varepsilon>0$,
    \begin{equation*}
        \P\left(\sup_{t_0\le s\le t\le t_0+\delta}
            W_p(\widehat\nu_s^N,\widehat\nu_t^N)>\varepsilon\right)
        \le\exp\!\left[-N(a_T\varepsilon^2-C_T\delta)\right].
    \end{equation*}
\end{proposition}
\begin{proof}
    Pairing $Y_s^i$ with $Y_t^i$ and using the power-mean
    inequality for $1\le p\le2$, we obtain
    \begin{equation*}
        \sup_{t_0\le s\le t\le t_0+\delta}
        W_p^2(\widehat\nu_s^N,\widehat\nu_t^N)
        \le\frac1N\sum_{i=1}^N
        \sup_{t_0\le s\le t\le t_0+\delta}
        d^2(Y_s^i,Y_t^i).
    \end{equation*}
    The processes $Y^1,\ldots,Y^N$ are independent. Thus the
    exponential increment estimate \eqref{eq:square exponential moment quantitative ind}
    below and exponential Chebyshev yield
    \begin{align*}
        &\P\left(\sup_{t_0\le s\le t\le t_0+\delta}
            W_p(\widehat\nu_s^N,\widehat\nu_t^N)>\varepsilon\right)\\
        &\qquad\le e^{-Na_T\varepsilon^2}
        (1+C_T\delta)^N
        \le e^{-N(a_T\varepsilon^2-C_T\delta)},
    \end{align*}
    as claimed.
\end{proof}

The required exponential increment bound is contained in the following lemma.

\begin{lemma}[Estimates for the noninteracting proxy]\label{thm:Estimates for the noninteracting proxy}
    In the setting of \cref{thm:propagation of moments mu_t}, let $(Y_t)_{t\ge0}$ solve \eqref{eq:linear SDE}.
    Then for every $T > 0$, there are constants $C_T, a_T > 0$ such that, for every
    $0 \le t_0 \le s \le t \le t_0 + \delta \le T$, $T > 0$,
    \begin{equation}\label{eq:second moment quantitative ind}
        \E\left[d^2(Y_t,Y_s)\right] \le C_T(t-s)
    \end{equation}
    and
    \begin{equation}\label{eq:square exponential moment quantitative ind}
        \E\left[
            \sup_{s,t \, : \, t_0 \le s\le t\le t_0+\delta}
            \exp\left(a_T d^2(Y_t,Y_s)\right)
        \right] \le 1 + C_T\delta.
    \end{equation}
\end{lemma} 
\begin{proof}[Proof of \cref{thm:Estimates for the noninteracting proxy}]
    Fix $0\le s \le t\le T$ and denote $r \coloneqq d(\cdot, Y_s)$.
    We use Kendall's radial It\^o formula
    \cite[Theorem 1.1]{kendall_radial_1987}
    (see also \cite[Theorem 3.5.1]{hsu_stochastic_2002}
    and \cite{arnaudon_harnack_2006}). For $t>s$, it gives
    \begin{equation}\label{eq:Kendall r}
        \begin{aligned}
            \dd r(Y_t)
            ={}&\sqrt2\,\dd\beta_t
                +\bigl((\Delta r)_{\mathrm{reg}}(Y_t)
                +b_t(Y_t)\cdot\nabla r(Y_t)\bigr)\dd t\\
                &-\dd L_t^{\mathrm{cut}},
        \end{aligned}
    \end{equation}
    where $\beta_t$ is a one-dimensional Brownian motion, $(\Delta r)_{\mathrm{reg}}$ the density of the absolutely continuous part of the distributional Laplacian (see \cite{mantegazza_mascellani_uraltsev_2014}), and
    $L^{\mathrm{cut}}$ is a nondecreasing process supported on
    visits to the cut locus. The derivatives of $r$ in the drift
    agree with the classical derivatives away from
    $\{Y_s\}\cup\operatorname{Cut}(Y_s)$; the subscript
    $\mathrm{reg}$ denotes the absolutely continuous part of
    the distributional Laplacian. Its singular contribution
    has the favorable nonpositive sign, as in
    \cite{mantegazza_mascellani_uraltsev_2014}.
    We may therefore discard $-\dd L_t^{\mathrm{cut}}$ in
    upper estimates.
    Hence, applying Itô's formula to $r^2$, noting that $d\langle r(Y)\rangle_t = 2dt$ and recalling \eqref{eq:comparison Laplacian},
    \begin{align}
        dr^2(Y_t)
        &\le
        2\sqrt2 r(Y_t) d\beta_t
        +
        \left[ 2m+2\sqrt{|K|(m-1)}\,r(Y_t)
                +2r(Y_t)b_t(Y_t)\cdot\nabla r(Y_t)
        \right]dt. \label{eq:increment R2 Kendall}
    \end{align}

    We now estimate the drift term. 
    From the Hessian lower bounds, $\Hess V\ge\beta g$ and $\Hess_xW(\cdot,y)\ge\widetilde\beta g$, we obtain
    \begin{equation*}
        b_{t}(Y_{t})\cdot\nabla r(Y_t) \le |\nabla V(Y_s)| + |\nabla(W*\mu_{t})(Y_s)| - \left(\beta+\widetilde\beta\right)r(Y_t).
    \end{equation*}
    Under the growth assumptions on $V$ and $W$, and recalling that $M_{2,T} \coloneqq \sup_{t \in [0,T]} M_2(t) < +\infty$ by \cref{thm:propagation of moments mu_t}, we have
    \begin{equation*}
        |\nabla V(Y_s)| + |\nabla(W*\mu_t)(Y_s)| \le A_V + L_V d(o, Y_s) + A_W + L_{W,1}\left(d(o, Y_s) + \sqrt{M_{2,T}}\right).
    \end{equation*}
    Inserting this into \eqref{eq:increment R2 Kendall} and using Cauchy--Schwarz inequality yields
    \begin{equation}\label{eq:R2 differential increment}
        d r^2(Y_t) \le 2\sqrt{2} r(Y_t) d\beta_t + C_T'\left(1 + d^2(o, Y_s) + r^2(Y_t)\right) dt 
    \end{equation}
    with
    \begin{equation*}
        C_T' \coloneqq \max\left\{3, (L_V + L_{W,1})^2, 2m + |K|(m-1) + \left(A_V + A_W + L_{W,1}\sqrt{M_{2,T}}\right)^2\right\}.
    \end{equation*}
    After localizing and taking conditional expectation with respect to $Y_s$, we obtain
    \begin{equation*}
        \frac{\dd}{\dd t}\E\left[r^2(Y_t) \,\Big|\, Y_s\right]
        \le
        C_T'\left(1+d^2(o, Y_s)\right) + C_T'\E\left[r^2(Y_t) \,\Big|\, Y_s\right].
    \end{equation*}
    Since $r(Y_s) = d(Y_s, Y_s) = 0$, by Gronwall's inequality,
    \begin{equation}\label{eq:conditional second increment}
        \E\left[r^2(Y_t) \,\Big|\, Y_s\right]
        \le
        C_T'\left(1+d^2(o, Y_s)\right)
        \int_s^t e^{C_T'(t-u)}\,du
        \le
        C_T'e^{C_T'T}
        \left(1+d^2(o, Y_s)\right)(t-s).
    \end{equation}
    Taking the expectation and recalling that $\sup_{0\le s\le T} \E\left[d^2(o,Y_s)\right] = M_{2,T}<+\infty$ by \cref{thm:propagation of moments mu_t} since $\operatorname{Law}(Y_s) = \mu_s$, we obtain \eqref{eq:second moment quantitative ind};
    \begin{equation*}
        \E\left[d^2(Y_t,Y_s)\right] = \E\left[r^2(Y_t)\right] \le C_T'e^{C_T'T} (1 + M_{2,T})(t-s).
    \end{equation*}
    
    We now estimate $e^{a(t)r^2(Y_t)}$, where $a(t)>0$ will be chosen below.
    From \eqref{eq:R2 differential increment}, we have
    \begin{equation*}
        d r^2(Y_t)
        \le
        2\sqrt2\,r(Y_t)\,d\beta_t
        +
        C_T'\left(1+d^2(o,Y_s)+r^2(Y_t)\right)\,dt, \qquad d\langle r^2(Y)\rangle_t = 8 r^2(Y_t) dt.
    \end{equation*}
    By Itô's formula,
    \begin{multline}\label{eq:Ito squared exponential moment estimate}
        d e^{a(t)r^2(Y_t)} \\ \le dM_t + e^{a(t)r^2(Y_t)}
        \Big[
            a(t)C_T'\left(1+d^2(o,Y_s)\right)
            +
            \left(\frac{\dd}{\dd t}a(t)+C_T'a(t)+4a(t)^2\right)r^2(Y_t)
        \Big]dt
    \end{multline}
    for some local martingale $M_t$. We choose \(a\) as the solution of
\begin{equation}\label{eq:increment-a-ode}
    a'(u)+C_T'a(u)+4a(u)^2=0,
    \qquad
    a(s)\coloneqq a_0>0 .
\end{equation}
Then \(a\) is nonincreasing on \([s,T]\). Hence
\eqref{eq:Ito squared exponential moment estimate} becomes
\begin{equation}\label{eq:Ito squared exponential moment estimate became}
    d e^{a(u)r^2(Y_u)}
    \le
    dM_u
    +
    C_T'a(u)\left(1+d^2(o,Y_s)\right)
    e^{a(u)r^2(Y_u)}\,du .
\end{equation}
After localization and conditioning on \(Y_s\), combining Gronwall's inequality with $a(u) \le a_o$ gives
\begin{equation}\label{eq:conditional square-exp}
    \E\left[
        e^{a(t)d^2(Y_t,Y_s)}
        \,\middle|\,Y_s
    \right]
    \le
    \exp\Big(
        C_T'a_0(t-s)\left(1+d^2(o,Y_s)\right)
    \Big).
\end{equation}
Taking expectations yields
\[
    \E\left[
        e^{a(t)d^2(Y_t,Y_s)}
    \right]
    \le
    e^{C_T'a_0(t-s)}
    \E\left[
        e^{C_T'a_0(t-s)d^2(o,Y_s)}
    \right].
\]

    In particular, by choosing $a_o$ small enough, say
    \begin{equation}\label{eq:regime on a zero}
        0 < a_0 < \frac{\inf_{0 \le \tau \le T} \alpha(\tau)}{T C_T'},
    \end{equation}
    where $\alpha$ is the continuous function of \cref{thm:propagation of moments mu_t}, we obtain
    \begin{equation}\label{eq:to be uniform in time}
        \E\left[e^{a(t) d^2(Y_t, Y_s)}\right] \le 1 + C_T(t-s)
    \end{equation}
    for some constant $C_T > 0$; one may take
\[
    C_T
    \coloneqq
    e^{C'_Ta_0T}
    C'_Ta_0
    \left(
        1+\frac{E_T-1}{\underline{\alpha}_T}
    \right)
\]
with
\[\underline{\alpha}_T \coloneqq \inf_{0 \le \tau \le T} \alpha(\tau), \qquad E_T \coloneqq \sup_{0 \le \tau \le T} \int_M e^{{\underline{\alpha}_T}d^2(o,x)}\dd\mu_\tau(x).\]

We now aim to obtain \eqref{eq:to be uniform in time} with a supremum over
\(0\le t_0\le s\le t\le t_0+\delta\le T\) inside the expectation. Fix
\(t_0\in[0,T-\delta]\), set
\[
    r\coloneqq d(\cdot,Y_{t_0}),
    \qquad
    S_u\coloneqq
    \sup_{t_0\le v\le u}
    e^{a(v)r^2(Y_v)}.
\]
Applying \eqref{eq:Ito squared exponential moment estimate became} with
\(s=t_0\), integrating, and taking the supremum, we get, after localization,
\[
    S_u
    \le
    1
    +
    \sup_{t_0\le \tau\le u}\abs{M_\tau-M_{t_0}}
    +
    C_T'a_0\left(1+d^2(o,Y_{t_0})\right)
    \int_{t_0}^u S_\tau\dd\tau .
\]
Here the martingale part is given by
\[
    dM_u
    =
    2\sqrt2\,a(u)r(Y_u)e^{a(u)r^2(Y_u)}\,d\beta_u,
\]
and therefore
\[
    d\langle M\rangle_u
    \le
    8a_0^2r^2(Y_u)e^{2a(u)r^2(Y_u)}\,du.
\]
Hence, by the Burkholder--Davis--Gundy inequality,
\[
\begin{aligned}
    \E\left[
        \sup_{t_0\le \tau\le u}\abs{M_\tau-M_{t_0}}
        \,\middle|\,Y_{t_0}
    \right]
    &\le
    C_{\mathrm{BDG}}\sqrt{8}a_o
    \left(
        \int_{t_0}^u
        \E\left[
            r^2(Y_\tau)e^{2a(\tau)r^2(Y_\tau)}
            \,\middle|\,Y_{t_0}
        \right]\dd\tau
    \right)^{1/2}.
\end{aligned}
\]
Choosing $a_0>0$ smaller if necessary, the following bound is obtained by an additional It\^o calculation, which also explains the factor $\tau-t_0$. Set $A\coloneqq C'_T$, $H\coloneqq1+d^2(o,Y_{t_0})$, and $R_u \coloneqq r^2(Y_u)$. For $b_0>2a_0$ sufficiently small, let $b'+Ab+4b^2=0$, $b(t_0)=b_0$, with $b(u)\ge2a(u)$ on $[t_0,T]$. The semimartingale estimate \eqref{eq:R2 differential increment} and It\^o's formula for $R_u e^{b(u)R_u}$ give, after localization,
\[
    \frac{d}{du}\E\left[R_u e^{b(u)R_u}\mid Y_{t_0}\right]
    \le AH\E\left[e^{b(u)R_u}\mid Y_{t_0}\right]
    +(A+8b_0+Ab_0H)\E\left[R_u e^{b(u)R_u}\mid Y_{t_0}\right].
\]
The coefficient of $R_u^2e^{b(u)R_u}$ vanishes by the differential equation for $b$, and the conditional exponential estimate above bounds the first term by $AH e^{Ab_0(u-t_0)H}$. Gronwall's inequality, recalling that $R_{t_0}=0$, gives
\[
    \E\left[
        r^2(Y_\tau)e^{2a(\tau)r^2(Y_\tau)}
        \,\middle|\,Y_{t_0}
    \right]
    \le
    C_T(\tau-t_0)\left(1+d^2(o,Y_{t_0})\right)
    e^{C_T\delta(1+d^2(o,Y_{t_0}))}.
\]
Consequently, after increasing $C_T$ that we do not relabel,
\[
    \E\left[
        \sup_{t_0\le \tau\le u}\abs{M_\tau-M_{t_0}}
        \,\middle|\,Y_{t_0}
    \right]
    \le
    C_T\delta\left(1+d^2(o,Y_{t_0})\right)
    e^{C_T\delta(1+d^2(o,Y_{t_0}))}.
\]
Hence,
\begin{equation*}
    \E\left[S_u \,\middle|\, Y_{t_0} \right]
    \le
    1
    +
    C_T\delta\left(1+d^2(o,Y_{t_0})\right)
    e^{C_T\delta(1+d^2(o,Y_{t_0}))}
    +
    C_T\left(1+d^2(o,Y_{t_0})\right)
    \int_{t_0}^u \E\left[S_\tau \,\middle|\, Y_{t_0} \right]\dd\tau .
\end{equation*}
By Gronwall's inequality, after increasing \(C_T\),
\[
    \E\left[
        \sup_{t_0\le u\le t_0+\delta}
        e^{a(u)d^2(Y_u,Y_{t_0})}
        \,\middle|\,Y_{t_0}
    \right]
    \le
    1
    +
    C_T\delta\left(1+d^2(o,Y_{t_0})\right)
    e^{C_T\delta(1+d^2(o,Y_{t_0}))}.
\]
For any $s,t\in[t_0,t_0+\delta]$, the triangle inequality gives
\[
    d^2(Y_t,Y_s)
    \le
    2d^2(Y_t,Y_{t_0})
    +
    2d^2(Y_s,Y_{t_0})
    \le
    4\sup_{t_0\le u\le t_0+\delta}
    d^2(Y_u,Y_{t_0}).
\]
Moreover, the solution of \eqref{eq:increment-a-ode} admits a
strictly positive lower bound $a_*=a_*(T)>0$, uniform in
$t_0\in[0,T]$. Hence, setting \(
    a_T\coloneqq \frac{a_*}{4},
\)
we obtain
\[
    \sup_{t_0\le s\le t\le t_0+\delta}
    e^{a_Td^2(Y_t,Y_s)}
    \le
    \sup_{t_0\le u\le t_0+\delta}
    e^{a(u)d^2(Y_u,Y_{t_0})}.
\]
Choosing $a_0>0$ small enough, the propagated square-exponential moment, together with the fact that $a(u)$ is uniformly bounded below on intervals of length at most $T$ yield, after possibly decreasing $a_T > 0$,
\[
    \E\left[
        \sup_{t_0\le s\le t\le t_0+\delta}
        e^{a_T d^2(Y_t,Y_s)}
    \right]
    \le 1+C_T\delta .
\]
for large enough $C_T > 0$.

\end{proof}

\subsection{Coupling the interacting system and non interacting proxy}

We first derive a distance estimate for two diffusions with possibly different drifts, coupled by parallel transport of their Brownian noises under a lower Ricci curvature bound. We follow the parallel-coupling construction in \cite[Section 2.1]{wang_analysis_2013}.

Let \(X_t\) and \(Y_t\) solve
\begin{equation*}
    dX_t=\sqrt2\,dB_t^M+b_t^X(X_t)\,dt,
    \qquad
    dY_t=\sqrt2\,d\overline B_t^M+b_t^Y(Y_t)\,dt,
\end{equation*}
where $B_t^M$ and $\overline B_t^M$ are Brownian motions on $M$. Equivalently, on the orthonormal frame bundle $\pi : O(M) \to M$ in Stratonovich form, denoting by tildes horizontal lifts,
\begin{equation*}
    d\Xi_t
    =
    \sqrt2\sum_{k=1}^m H_k(\Xi_t)\circ dB_t^k
    +
    \widetilde{b_t^X}(\Xi_t)\,dt,
    \qquad
    X_t=\pi(\Xi_t),
\end{equation*}
and
\begin{equation*}
    dU_t
    =
    \sqrt2\sum_{k=1}^m H_k(U_t)\circ d\overline B_t^k
    +
    \widetilde{b_t^Y}(U_t)dt,
    \qquad
    Y_t=\pi(U_t).
\end{equation*}
Away from the joint cut locus, let \(\gamma_t:[0,d(X_t,Y_t)]\to M\) be the unit-speed minimizing geodesic from \(X_t\) to \(Y_t\), and denote by
\begin{equation*}
    P_{X_t\to Y_t}:T_{X_t}M\to T_{Y_t}M
\end{equation*}
parallel transport along \(\gamma_t\). We couple the Brownian noises by parallel displacement. More precisely, for each pair $(X_t^i,Y_t^i)$ we set
\[
    d\overline B_t^i = O_t^i\,dB_t^i,
    \qquad
    O_t^i \coloneqq (U_t^i)^{-1}
    P_{X_t^i\to Y_t^i}\Xi_t^i .
\]
This parallel coupling is extended through the cut locus via the standard approximation procedure; see
\cite[Section 2.1]{wang_analysis_2013}.
Since the block-diagonal matrix $\operatorname{diag}(O_t^1,\dots,O_t^N)$ is orthogonal, $(\overline B^1,\dots,\overline B^N)$ is again a standard $\mathbb R^{Nm}$-valued Brownian motion. In particular, $\overline B^1,\dots,\overline B^N$ are independent, and consequently so are $Y^1,\dots,Y^N$, with common law $\mu_t$. The Brownian part of the coupled process \((X_t,Y_t)\) is therefore driven along the vector fields
\[
    \sqrt2\,(E_k(t),\overline E_k(t)),
    \qquad k=1,\dots,m,
\]
where
\[
    E_k(t)\coloneqq \Xi_t e_k,
    \qquad
    \overline E_k(t)\coloneqq P_{X_t\to Y_t}E_k(t).
\]

Applying the distance It\^o formula
\cite[Theorem 2.3.2]{wang_analysis_2013} to
$\rho\coloneqq d(\cdot,\cdot)$ on $M\times M$
(see also \cite{arnaudon_harnack_2006,hsu_stochastic_2002}),
one obtains
\begin{equation}\label{eq:Kendall parallel coupling}
    d\rho(X_t,Y_t)
    \le
    \mathcal I_tdt
    +
    b_t^X(X_t)\cdot\nabla_x\rho(X_t,Y_t)dt
    +
    b_t^Y(Y_t)\cdot\nabla_y\rho(X_t,Y_t)dt,
\end{equation}
where 
\begin{equation}\label{eq:index term}
    \mathcal I_t
    \coloneqq
    \sum_{k=1}^m
    \Hess\rho_{(X_t,Y_t)} \Big(\left(E_k, \overline{E}_k\right)(t), \left(E_k, \overline{E}_k\right)(t)\Big).
\end{equation}
Away from the joint cut locus and diagonal, \eqref{eq:Kendall parallel coupling} is an equality. At points where the distance function is not smooth, we use the generalized Itô--Kendall formula, whose cut-locus contribution is a nonpositive finite-variation term that can be discarded in upper bounds; compare \eqref{eq:Kendall r}.

Recall that there is no martingale term due to the parallel coupling. Indeed,
\begin{equation*}
    dZ_t = \sqrt{2}\sum_{k=1}^m \left(E_k(t), \overline E_k(t)\right) \circ dB_t^k + \left(b_t^X(X_t), b_t^Y(Y_t)\right) dt,
\end{equation*}
and applying Itô's formula to $\rho$ yields then
\begin{equation*}
    d\rho(Z_t) =
    \sqrt{2}\sum_{k=1}^m d\rho_{Z_t} \left(E_k(t), \overline E_k(t)\right) dB_t^k +
    \mathcal I_tdt
    +
    b_t^X(X_t)\cdot\nabla_x\rho(Z_t)dt
    +
    b_t^Y(Y_t)\cdot\nabla_y\rho(Z_t)dt
\end{equation*}
But
\begin{equation*}
    \nabla_x\rho(Z_t)=-\dot\gamma_t(0),
    \qquad
    \nabla_y\rho(Z_t)=\dot\gamma_t(\rho(Z_t)),
\end{equation*}
and thus, since \(\overline E_k(t)=P_{X_t\to Y_t}E_k(t)\),
\begin{align*}
    d_x\rho(Z_t) E_k(t) + d_y\rho(Z_t) \overline E_k(t)
    &=
    - E_k(t) \cdot \dot\gamma_t(0) + \overline E_k(t) \cdot \dot\gamma_t(\rho(Z_t)) 
    \\ &=
    - E_k(t) \cdot \dot\gamma_t(0) +
    E_k(t) \cdot \dot\gamma_t(0) \\ &=0.
\end{align*}

We now estimate the index term \eqref{eq:index term} under the Ricci lower bound \(\Ric\ge Kg\) with \(K \le 0\).
We choose $E_1(t) \coloneqq \dot\gamma_t(0)$ so the radial direction does not contribute to the Hessian of \(\rho\).
Hence
\begin{equation*}
    \mathcal I_t
    =
    \sum_{k=2}^m
    \Hess_{(X_t,Y_t)}\rho
    \Big(
        (E_k(t),\overline E_k(t)),
        (E_k(t),\overline E_k(t))
    \Big).
\end{equation*}
For \(k=2,\dots,m\), let \(J_k\) be the parallel vector field along \(\gamma_t\) such that
\begin{equation*}
    J_k(0)=E_k(t),
    \qquad
    J_k(\rho(X_t,Y_t))=\overline E_k(t).
\end{equation*}
By the index lemma \cite[Lemma 6.7.1]{hsu_stochastic_2002},
\begin{equation*}
    \Hess_{(X_t,Y_t)}\rho
    \Big(
        (E_k(t),\overline E_k(t)),
        (E_k(t),\overline E_k(t))
    \Big)
    \le
    I_{\gamma_t}(J_k,J_k),
\end{equation*}
where
\begin{equation*}
    I_{\gamma_t}(J,J)
    \coloneqq
    \int_0^{\rho(X_t,Y_t)}
    \left(
        |\nabla_{\dot\gamma_t}J|^2
        -
        \langle R(J,\dot\gamma_t)\dot\gamma_t,J\rangle
    \right)\dd s.
\end{equation*}
But \(J_k\) is parallel, and therefore,
\begin{equation*}
    I_{\gamma_t}(J_k,J_k)
    = -
    \int_0^{\rho(X_t,Y_t)}
    \langle R(J_k,\dot\gamma_t)\dot\gamma_t,J_k\rangle\dd s.
\end{equation*}
Summing over \(k=2,\dots,m\), and using that \(\{J_k\}_{2\le k\le m}\) form a parallel orthonormal frame orthogonal to \(\dot\gamma_t\), we obtain
\begin{equation*}
    \mathcal I_t \le
    -
    \int_0^{\rho(X_t,Y_t)}
    \sum_{k=2}^m
    \langle R(J_k,\dot\gamma_t)\dot\gamma_t,J_k\rangle\dd s \\
    =
    -
    \int_0^{\rho(X_t,Y_t)}
    \Ric(\dot\gamma_t,\dot\gamma_t)\dd s \le -K\rho(X_t, Y_t).
\end{equation*}
Consequently,
\begin{equation}\label{eq:Kendall parallel coupling Ricci bound}
    d\rho(X_t,Y_t)
    \le
    \left[
        |K|\rho(X_t,Y_t)
        +
        b_t^X(X_t)\cdot\nabla_x\rho(X_t,Y_t)
        +
        b_t^Y(Y_t)\cdot\nabla_y\rho(X_t,Y_t)
    \right]dt.
\end{equation}

\begin{proposition}[Parallel-coupling estimate for empirical measures]\label{thm:coupling estimate empirical measures}
    Let
    \begin{equation*}
        \widehat\mu_t^N\coloneqq \frac1N\sum_{i=1}^N\delta_{X_t^{i}},
        \qquad
        \widehat\nu_t^N\coloneqq \frac1N\sum_{i=1}^N\delta_{Y_t^i},
    \end{equation*}
    where
    \begin{align*}
        dX_t^{i}
        &=
        \sqrt2\,dB_t^{i,M}
        -\nabla V(X_t^i)\,dt
        -\frac1N\sum_{j=1}^N\nabla_xW(X_t^i,X_t^j)\,dt,\\
        dY_t^i
        &=
        \sqrt2\,d\overline B_t^{i,M}
        -\nabla V(Y_t^i)\,dt
        -\int_M\nabla_xW(Y_t^i,z)\dd\mu_t(z)\,dt.
    \end{align*}
    Assume that $X_0^i=Y_0^i$ and couple each pair
    $(X_t^i,Y_t^i)$ by parallel coupling. Under the hypotheses
    of \cref{thm:propagation of moments mu_t} and
    \eqref{eq:assumption W second variable Lipschitz IPS},
    for every $p\ge1$ and $t\ge0$, it holds
    \begin{equation}\label{eq:coupling estimate empirical measures}
        \begin{aligned}
        W_p(\widehat\mu_t^N,\mu_t)
        \le{}& W_p(\widehat\nu_t^N,\mu_t)\\
        &+L_{W,2}\int_0^t e^{-\alpha_M(t-s)}
        W_p(\widehat\nu_s^N,\mu_s)\dd s,
        \end{aligned}
    \end{equation}
    where
    \begin{equation*}
        \alpha_M \coloneqq \beta + \widetilde\beta -|K| - L_{W,2}.
    \end{equation*}
\end{proposition}

\begin{proof}
    We apply the parallel-coupling Kendall formula
    \eqref{eq:Kendall parallel coupling Ricci bound}
    to the pair $(X_t^i,Y_t^i)$. Writing
    $\rho_i(t)\coloneqq\rho(X_t^i,Y_t^i)$, we obtain
    \begin{equation}\label{eq:rho_i Kendall bound}
        d\rho_i(t) \le
        \Big[
            |K|\rho_i(t)
            +
            b_t^{i}(X_t^i)\cdot\nabla_x\rho(X_t^i,Y_t^i)
            +
            \overline b_t(Y_t^i)\cdot\nabla_y\rho(X_t^i,Y_t^i)
        \Big]dt,
    \end{equation}
    where
    \begin{equation*}
        b_t^{i}(x)
        \coloneqq -\nabla V(x)
        - \frac1N\sum_{j=1}^N\nabla_xW(x,X_t^{j}),
        \qquad
        \overline b_t(x)
        \coloneqq -\nabla V(x)-\nabla_x(W*\mu_t)(x).
    \end{equation*}

    We first estimate the contribution of $V$.
    Since $\Hess V\ge \beta g$, integrating along the
    geodesic $\gamma_i$ from $X_t^i$ to $Y_t^i$ yields
    \begin{equation*}
        \nabla V(Y_t^i)\cdot\dot\gamma_i(\rho_i(t))
        -
        P_{X_t^i\to Y_t^i}\nabla V(X_t^i)
        \cdot\dot\gamma_i(\rho_i(t))
        \ge \beta\rho_i(t).
    \end{equation*}
    Equivalently,
    \begin{equation}\label{eq:V coupling convexity}
        -\nabla V(X_t^i)\cdot\nabla_x\rho(X_t^i,Y_t^i)
        -
        \nabla V(Y_t^i)\cdot\nabla_y\rho(X_t^i,Y_t^i)
        \le -\beta\rho_i(t).
    \end{equation}

    We now treat the interaction.
    Its contribution in \eqref{eq:rho_i Kendall bound}
    is given by
    \begin{align*}
        \mathcal W_i(t)
        \coloneqq
        &-\frac1N\sum_{j=1}^N
        \nabla_xW(X_t^i,X_t^j)
        \cdot\nabla_x\rho(X_t^i,Y_t^i)\\
        &-
        \int_M\nabla_xW(Y_t^i,z)\dd\mu_t(z)
        \cdot\nabla_y\rho(X_t^i,Y_t^i).
    \end{align*}
    After parallel transport, adding and subtracting
    vector fields in $T_{Y_t^i}M$, we deduce that
    \begin{align*}
        \mathcal W_i(t)
        &=
        \left(
            \frac1N\sum_{j=1}^N
            \left[
                P_{X_t^i\to Y_t^i}
                \nabla_xW(X_t^i,X_t^j)
                -
                \nabla_xW(Y_t^i,X_t^j)
            \right]
        \right)\cdot\nabla_y\rho(X_t^i,Y_t^i)\\
        &\quad+
        \left(
            \frac1N\sum_{j=1}^N
            \left[
                \nabla_xW(Y_t^i,X_t^j)
                -
                \nabla_xW(Y_t^i,Y_t^j)
            \right]
        \right)\cdot\nabla_y\rho(X_t^i,Y_t^i)\\
        &\quad+
        \left(
            \frac1N\sum_{j=1}^N
            \nabla_xW(Y_t^i,Y_t^j)
            -
            \int_M\nabla_xW(Y_t^i,z)\dd\mu_t(z)
        \right)\cdot\nabla_y\rho(X_t^i,Y_t^i).
    \end{align*}
    The lower Hessian bound
    $\Hess_xW(\cdot,z)\ge\widetilde\beta g$,
    for every $z\in M$, yields
    \begin{equation*}
        \left(
            P_{X_t^i\to Y_t^i}\nabla_xW(X_t^i,z)
            -
            \nabla_xW(Y_t^i,z)
        \right)
        \cdot\nabla_y\rho(X_t^i,Y_t^i)
        \le -\widetilde\beta\rho_i(t),
    \end{equation*}
    the assumption
    \eqref{eq:assumption W second variable Lipschitz IPS}
    on $\nabla_xW$ gives
    \begin{equation*}
        \abs{
            \left(
                \frac1N\sum_{j=1}^N
                \left[
                    \nabla_xW(Y_t^i,X_t^j)
                    -
                    \nabla_xW(Y_t^i,Y_t^j)
                \right]
            \right)\cdot\nabla_y\rho(X_t^i,Y_t^i)
        }
        \le
        \frac{L_{W,2}}N\sum_{j=1}^N\rho_j(t),
    \end{equation*}
    and furthermore using Kantorovich duality,
    \begin{equation*}
        \begin{aligned}
        &\abs{
            \left(
                \frac1N\sum_{j=1}^N
                \nabla_xW(Y_t^i,Y_t^j)
                -
                \int_M\nabla_xW(Y_t^i,z)\dd\mu_t(z)
            \right)\cdot\nabla_y\rho(X_t^i,Y_t^i)
        }\\
        &\qquad\le
        L_{W,2}W_1(\widehat\nu_t^N,\mu_t)
        \le L_{W,2}W_p(\widehat\nu_t^N,\mu_t).
        \end{aligned}
    \end{equation*}
    Hence, by H\"older's inequality,
    \begin{align}
        \mathcal W_i(t)
        &\le
        -\widetilde\beta\rho_i(t)
        +
        \frac{L_{W,2}}N\sum_{j=1}^N\rho_j(t)
        +
        L_{W,2}W_p(\widehat\nu_t^N,\mu_t)
        \notag\\
        &\le
        -\widetilde\beta\rho_i(t)
        +
        L_{W,2}R_N(t)
        +
        L_{W,2}W_p(\widehat\nu_t^N,\mu_t),
        \label{eq:interaction coupling estimate}
    \end{align}
    where
    \begin{equation*}
        R_N(t)\coloneqq
        \left(
            \frac1N\sum_{j=1}^N\rho_j^p(t)
        \right)^{1/p}.
    \end{equation*}

    Combining
    \eqref{eq:rho_i Kendall bound}--
    \eqref{eq:interaction coupling estimate},
    and since there is no martingale term due to
    the parallel coupling, It\^o's formula gives
    \begin{align*}
        d\rho_i^p(t)
        &=p\rho_i^{p-1}(t)\,d\rho_i(t)\\
        &\le
        \Big[
            -p(\beta+\widetilde\beta-|K|)
            \rho_i^p(t)\\
        &\qquad\quad+
            pL_{W,2}\rho_i^{p-1}(t)
            \left(
                R_N(t)+W_p(\widehat\nu_t^N,\mu_t)
            \right)
        \Big]dt.
    \end{align*}
    Averaging over $i=1,\dots,N$ and using
    H\"older's inequality,
    \begin{equation*}
        \frac1N\sum_{i=1}^N\rho_i^{p-1}(t)
        \le R_N^{p-1}(t),
    \end{equation*}
    we obtain
    \begin{equation*}
        dR_N^p(t)
        \le
        \Big[
            -p\alpha_M R_N^p(t)
            +
            pL_{W,2}R_N^{p-1}(t)
            W_p(\widehat\nu_t^N,\mu_t)
        \Big]dt,
    \end{equation*}
    where
    \begin{equation*}
        \alpha_M\coloneqq
        \beta+\widetilde\beta-|K|-L_{W,2}.
    \end{equation*}
    Applying the chain rule to $R_N$, with a
    regularization at $R_N=0$, yields
    \begin{equation*}
        dR_N(t)
        \le
        \Big[
            -\alpha_M R_N(t)
            +
            L_{W,2}W_p(\widehat\nu_t^N,\mu_t)
        \Big]dt.
    \end{equation*}
    Since $X_0^i=Y_0^i$, we have $R_N(0)=0$.
    Gronwall's inequality yields
    \begin{equation*}
        R_N(t)
        \le
        L_{W,2}\int_0^t e^{-\alpha_M(t-s)}
        W_p(\widehat\nu_s^N,\mu_s)\dd s.
    \end{equation*}
    Finally, pairing each $X_t^i$ with $Y_t^i$
    gives
    \begin{equation*}
        W_p(\widehat\mu_t^N,\widehat\nu_t^N)
        \le
        \left(
            \frac1N\sum_{i=1}^N
            d^p(X_t^i,Y_t^i)
        \right)^{1/p}
        =
        R_N(t).
    \end{equation*}
    Therefore, by the triangle inequality,
    \begin{equation*}
        \begin{aligned}
            W_p(\widehat\mu_t^N,\mu_t)
            &\le
            W_p(\widehat\mu_t^N,\widehat\nu_t^N)
            +
            W_p(\widehat\nu_t^N,\mu_t)\\
            &\le
            L_{W,2}\int_0^t e^{-\alpha_M(t-s)}
            W_p(\widehat\nu_s^N,\mu_s)\dd s
            +
            W_p(\widehat\nu_t^N,\mu_t),
        \end{aligned}
    \end{equation*}
    which proves \eqref{eq:coupling estimate empirical measures}.
\end{proof}

\begin{remark}[Propagation of $T_2$ along the nonlinear diffusion]
\label{rem:aggregation-higher-p}
    The propagation of square-exponential moments
    \cref{thm:propagation of moments mu_t} only gives a common
    $T_1$ inequality on each interval $[0,T]$. We recall that $T_p$ for $p > 1$ is not simply characterized by moments. Nevertheless we can propagate $T_2$ by identifying a common lower bound of the Bakry--\'Emery tensor. Set
    \[
        \mathcal U_t\coloneqq V+W*\mu_t,
        \qquad
        \kappa\coloneqq K+\beta+\widetilde\beta.
    \]
    The assumptions on the potentials give
    \[
        \Ric+\Hess\mathcal U_t\ge\kappa g
        \qquad\text{for every }t\ge0.
    \]
    If $\mu_0$ satisfies $T_2(\lambda_{2,0})$, then
    $\mu_t$ satisfies $T_2(\lambda_2(t))$, with the
    admissible constant
    \begin{equation}\label{eq:transport-propagation-T2}
        \lambda_2(t)^{-1}
        =e^{-2\kappa t}\lambda_{2,0}^{-1}
        +2\int_0^t e^{-2\kappa s}\dd s.
    \end{equation}

    For completeness, we recall the coupling--Girsanov argument,
    following the approach of \cite[Section~2.3]{cattiaux_guillin_2014}.
    Fix $t>0$ and $\nu\ll\mu_t$ with finite entropy. Tilt the
    path law $\mathbb P$ of the diffusion by
    \[
        \frac{\dd\mathbb Q}{\dd\mathbb P}
        =\frac{\dd\nu}{\dd\mu_t}(Y_t).
    \]
    Under $\mathbb Q$, Girsanov's theorem adds a drift
    $\sqrt2\,u_s$ to $-\nabla\mathcal U_s$, and the entropy
    decomposition reads
    \[
        H(\nu\mid\mu_t)
        =H(\nu_0\mid\mu_0)
        +\frac12\E_{\mathbb Q}
            \int_0^t|u_s|^2\dd s,
    \]
    where $\nu_0$ is the initial law under $\mathbb Q$.
    Couple the controlled diffusion $Y_s$ with an
    unperturbed diffusion $\widetilde Y_s$ started from $\mu_0$,
    choosing an optimal $W_2$ coupling of the initial laws
    and using parallel noise. Set
    $\rho_s\coloneqq d(Y_s,\widetilde Y_s)$.
    Applying \eqref{eq:Kendall parallel coupling Ricci bound}
    and using $\Hess\mathcal U_s\ge\Lambda g$, we obtain
    \[
        \dd\rho_s
        \le-\kappa\rho_s\dd s+\sqrt2\,|u_s|\dd s.
    \]
    Minkowski's and Cauchy--Schwarz inequalities give
    \begin{align*}
        W_2(\nu,\mu_t)
        &\le e^{-\kappa t}
            \sqrt{\frac{2H(\nu_0\mid\mu_0)}{\lambda_{2,0}}}\\
        &\quad+2\sqrt{\int_0^t e^{-2\kappa s}\dd s}
            \left(\frac12\E_{\mathbb Q}
                \int_0^t|u_s|^2\dd s\right)^{1/2}.
    \end{align*}
    Another Cauchy--Schwarz inequality gives
    \[
        W_2^2(\nu,\mu_t)
        \le 2\left(
            e^{-2\kappa t}\lambda_{2,0}^{-1}
            +2\int_0^t e^{-2\kappa s}\dd s\right)H(\nu\mid\mu_t),
    \]
    proving \eqref{eq:transport-propagation-T2}. As usual,
    the argument is first localized away from the cut locus
    and then passed to the limiting parallel coupling.

    In particular, all $\mu_t$, $0\le t\le T$, satisfy
    $T_2(\lambda_{2,T})$, with
    \begin{equation*}
        \lambda_{2,T}^{-1}
        =\max\left\{
            \lambda_{2,0}^{-1},
            e^{-2\kappa T}\lambda_{2,0}^{-1}
                +2\int_0^T e^{-2\kappa s}\dd s
        \right\}.
    \end{equation*}
    If $\kappa>0$, one may simply take
    $\lambda_{2,T}=\min\{\lambda_{2,0},\kappa\}$.
    Of course, this yields
    a common $T_p(\lambda_{2,T})$ inequality as well.
    Together with the separately propagated square-exponential
    moments, it supplies the uniform constants required in \cref{thm:W_p quantitative concentration} for $1< p\le2$.
\end{remark}

\begin{remark}[Propagation of chaos and comparison with recent work]
\label{rem:chaos}
    Let $1\le p\le2$ and suppose that $m>2p$, so that
    \cref{thm:sharp E nu n n K neg} applies. The moment estimate
    in \cref{thm:propagation of moments mu_t} gives a uniform
    exponential moment of $\mu_t$ for $t\in[0,T]$. Consequently,
    \cref{thm:sharp E nu n n K neg} and
    \eqref{eq:coupling estimate empirical measures} imply
    \begin{equation*}
        \sup_{0\le t\le T}\E W_p(\widehat\mu_t^N,\mu_t)
        \le C_{T,p}\left(
            1+L_{W,2}\int_0^T e^{-\alpha_M s}\dd s
        \right)N^{-1/m}.
    \end{equation*}
    In particular, when $p=1$ and $m>2$, exchangeability and
    the same particlewise coupling give, for any fixed $k\ge1$,
    \begin{equation*}
        \begin{aligned}
        &W_{1,d_k}\!\left(
            \operatorname{Law}(X_t^1,\ldots,X_t^k),
            \mu_t^{\otimes k}
        \right)\\
        &\qquad\le kC_{T,1}L_{W,2}
            \left(\int_0^T e^{-\alpha_M s}\dd s\right)N^{-1/m},
        \end{aligned}
    \end{equation*}
    uniformly over $t\in[0,T]$, where
    $d_k((x_i)_{i=1}^k,(y_i)_{i=1}^k)\coloneqq\sum_{i=1}^k d(x_i,y_i)$.

    This yields propagation of chaos from the full-empirical-measure
    approximation. The recent preprint of Chae--Cho--Jang--Yang
    \cite{chae_cho_jang_yang_2026} establishes a stronger
    $N^{-1/2}$ finite-particle chaos rate, uniformly in time
    under its hypotheses, in a modified Wasserstein distance
    using reflection coupling. Our dimension-dependent
    $N^{-1/m}$ rate concerns the entire empirical measure in
    the original Wasserstein distance, and is not claimed to
    be optimal for finitely many particles.
\end{remark}

\begin{proof}[Proof of \cref{thm:main aggregation_diffusion}]
    We follow the final discretization step of
    \cite[Theorem 2.9]{bolley_quantitative_2007}.
    Fix $1\le p\le2$. For $p=1$, the square-exponential
    moment bound gives a common $T_1$ constant on $[0,T]$.
    For $1<p\le2$, the additional initial $T_2$ assumption
    and \cref{rem:aggregation-higher-p} give a common
    $T_p(\lambda_{2,T})$ constant. The coupling estimate gives
    \[
        \sup_{t\le T}W_p(\widehat\mu_t^N,\mu_t)
        \le A_T\sup_{t\le T}W_p(\widehat\nu_t^N,\mu_t),\qquad
        A_T=1+L_{W,2}\int_0^T e^{-\alpha_M s}\dd s.
    \]
    Take a time grid of mesh $\delta=c_T\varepsilon^2$ in $[0,T]$. At each grid point, the fixed-time quantitative concentration theorem applies to the independent proxy particles with this common transport constant and the square-exponential moment bound supplied by \cref{thm:propagation of moments mu_t}. It gives a tail $\exp(-c_TN\varepsilon^2)$ once $N\ge N_0\varepsilon^{-(m'+2)}$. To control increments on each grid interval, the preceding
    proxy estimate and independence give, by exponential Chebyshev,
    \[
        \P\left(\sup_{s,t\text{ in one interval}}
        W_p(\widehat\nu_s^N,\widehat\nu_t^N)>u\right)
        \le\exp\!\left[-N\big(c_Tu^2-C_T\delta\big)\right].
    \]
    Moreover, for $1\le p\le2$, the second-moment increment bound gives
    \[
        W_p(\mu_s,\mu_t)
        \le \big(\E\left[d^p(Y_s,Y_t)\right]\big)^{1/p}
        \le \big(\E\left[d^2(Y_s,Y_t)\right]\big)^{1/2}
        \le C_T\sqrt{|t-s|}.
    \]
    Choosing $c_T>0$ small enough that the deterministic and random increments can each be absorbed into a fixed fraction of $\varepsilon/A_T$, and taking a union bound over at most $1+T/\delta$ intervals, yields
    \[
        \P\left(\sup_{t\le T}W_p(\widehat\mu_t^N,\mu_t)>\varepsilon\right)
        \le C_T(1+T\varepsilon^{-2})e^{-c_TN\varepsilon^2},
    \]
    as claimed.
\end{proof}

\clearpage
\appendix
\section{Auxiliary proofs for the nonlinear diffusion}
\subsection{Well-posedness of the nonlinear diffusion}
\label{app:wellposed-mv}

\begin{proof}[Proof of \cref{prop:wellposed-mv}]
    Fix $T>0$ and let
    $\nu\in C([0,T];\mathcal P_2(M))$ with $\nu_0=\mu_0$.
    We consider the frozen drift
    \begin{equation*}
        b_t^\nu(x)
        \coloneqq -\nabla V(x)
        -\int_M\nabla_xW(x,y)\dd\nu_t(y).
    \end{equation*}
    By \eqref{eq:assumptions V IPS} and
    \eqref{eq:assumptions W growth IPS}, the drift $b_t^\nu$
    has at most linear growth in $r(x)=d(o,x)$,
    uniformly over $t\in[0,T]$. Moreover, the Hessian
    lower bounds give one-sided dissipativity.
    More precisely, combining the radial estimates
    \eqref{eq:comparison Laplacian}--\eqref{eq:r2 W upper},
    we obtain, outside $\{o\}\cup\operatorname{Cut}(o)$,
    \begin{equation}\label{eq:frozen-radial-Lyapunov}
        (\Delta+b_t^\nu\cdot\nabla)r^2(x)
        \le 2m+2C_{\nu,T}r(x)-2\Lambda r^2(x),
    \end{equation}
    where $\Lambda\coloneqq\beta+\widetilde\beta$ and
    \begin{equation*}
        \begin{aligned}
            C_{\nu,T}
            \coloneqq{}& \sqrt{|K|(m-1)}+|\nabla V(o)|+A_W\\
            &+L_{W,1}\sup_{0\le t\le T}
                \int_M r(y)\dd\nu_t(y)<\infty.
        \end{aligned}
    \end{equation*}
    Let $X^\nu$ be the maximal frozen diffusion obtained
    by the standard localized construction on complete
    manifolds; see \cite{hsu_stochastic_2002,zhang_2005}.
    For $R>0$, let $\tau_R$ be its first exit time from
    $B(o,R)$. We apply the radial It\^o--Kendall formula
    \eqref{eq:Kendall r}, now with the fixed base point $o$,
    to $r^2(X^\nu_{t\wedge\tau_R})$.
    The nonpositive cut-locus contribution can be discarded;
    combining this with \eqref{eq:frozen-radial-Lyapunov} and
    $2C_{\nu,T}r-2\Lambda r^2\le C_{\nu,T}^2/(2\Lambda)$ gives
    \begin{equation*}
        \mathbb E\bigl[r^2(X^\nu_{t\wedge\tau_R})\bigr]
        \le \mathbb E\bigl[r^2(X^\nu_0)\bigr]
        +\left(2m+\frac{C_{\nu,T}^2}{2\Lambda}\right)t,
        \qquad 0\le t\le T.
    \end{equation*}
    In particular, the right-hand side is independent of
    $R$, and letting $R\to\infty$ rules out explosion.
    Denote the law curve of the resulting frozen diffusion
    by $\mathcal F(\nu)$. Path continuity and the standard
    localized increment estimates yield continuity
    in $\mathcal P_2(M)$, so that $\mathcal F$ maps
    $C([0,T];\mathcal P_2(M))$ into itself.

    We now compare two frozen law curves $\nu$ and $\eta$
    with the same initial law $\mu_0$. Let $X^\nu$ and
    $X^\eta$ denote the corresponding frozen diffusions,
    started from the same random variable with law $\mu_0$
    and coupled by parallel transport. Set
    $\rho_s\coloneqq d(X_s^\nu,X_s^\eta)$.
    Applying the parallel-coupling Kendall estimate
    \eqref{eq:Kendall parallel coupling Ricci bound},
    together with the Hessian lower bounds and
    \eqref{eq:assumption W second variable Lipschitz IPS},
    we obtain
    \begin{equation*}
        d\rho_s
        \le\bigl[-(\Lambda-|K|)\rho_s
            +L_{W,2}W_2(\nu_s,\eta_s)\bigr]\dd s.
    \end{equation*}
    Since $\rho_0=0$, Gronwall's inequality yields
    \begin{equation*}
        \rho_t\le L_{W,2}\int_0^t
        e^{-(\Lambda-|K|)(t-s)}W_2(\nu_s,\eta_s)\dd s.
    \end{equation*}
    Taking the $L^2$ norm of the coupled distance gives
    \begin{equation}\label{eq:frozen-law-stability}
        W_2\bigl(\mathcal F(\nu)_t,\mathcal F(\eta)_t\bigr)
        \le L_{W,2}\int_0^t
        e^{-(\Lambda-|K|)(t-s)}W_2(\nu_s,\eta_s)\dd s.
    \end{equation}

    Equip the space of continuous $\mathcal P_2(M)$-valued
    curves starting from $\mu_0$ with the weighted distance
    \begin{equation*}
        d_q(\nu,\eta)
        \coloneqq\sup_{0\le t\le T}
        e^{-qt}W_2(\nu_t,\eta_t).
    \end{equation*}
    This space is complete, and for
    \(
        q> (L_{W,2}-\Lambda+|K|)_+,
    \)
    it follows from \eqref{eq:frozen-law-stability} that
    $\mathcal F$ is a strict contraction.
    Hence, Banach's fixed-point theorem gives a unique law curve
    on $[0,T]$. Since $T$ is arbitrary and the fixed points
    agree on overlapping intervals, the solution is global
    and nonexplosive. Finally, testing the frozen equation
    at its fixed point against compactly supported smooth
    functions gives \eqref{eq:aggregation diffusion}.
\end{proof}

\subsection{Radial cutoff and integration by parts}
\label{app:radial-cutoff}

\begin{proof}[Proof of \cref{lem:radial-cutoff}]
    Let $(X_t)_{0\le t\le T}$ be the nonexplosive diffusion
    constructed in \cref{prop:wellposed-mv}, with
    $\operatorname{Law}(X_t)=\mu_t$. Set $r\coloneqq d(o,\cdot)$.
    Recall Kendall's formula
    \eqref{eq:Kendall r}, now with base point $o$.
    The distributional Laplacian of the distance satisfies
    \begin{equation*}
        \Delta r=(\Delta r)_{\mathrm{reg}}\,\dd\Vol-\sigma_r,
        \qquad \sigma_r\ge0,
    \end{equation*}
    where $\sigma_r$ is supported on $\operatorname{Cut}(o)$;
    see \cite{mantegazza_mascellani_uraltsev_2014}.
    Here $(\Delta r)_{\mathrm{reg}}$ agrees with the classical
    Laplacian outside $\{o\}\cup\operatorname{Cut}(o)$.
    For $\Psi_t=\Phi(t,r^2)$, the regular generator is
    \begin{equation*}
        \begin{aligned}
            (L_t\Psi_t)_{\mathrm{reg}}
            ={}&\partial_2\Phi(t,r^2)
                \bigl(2+2r(\Delta r)_{\mathrm{reg}}
                    +2r\,b_t\cdot\nabla r\bigr)\\
                &+4r^2\partial_{22}\Phi(t,r^2).
        \end{aligned}
    \end{equation*}
    Note that the singular part 
    $-2r\,\partial_2\Phi(t,r^2)\,\sigma_r$ is nonpositive since $\Psi_t$ is nondecreasing in $r$.

    Applying Itô's formula first to $\Psi_t(X_t)$ and then
    to $\chi_j(\Psi_t(X_t))$, we obtain, after localization,
    \begin{equation*}
        \begin{aligned}
            \dd\chi_j(\Psi_t(X_t))
            \le{}&\dd M_t^{(j)}
                +\chi_j'(\Psi_t(X_t))\\
                &\quad\times\bigl(\partial_t\Psi_t
                    +(L_t\Psi_t)_{\mathrm{reg}}\bigr)(X_t)\dd t,
        \end{aligned}
    \end{equation*}
    where $M^{(j)}$ is a local martingale. In addition
    to the nonpositive cut-locus term, we have discarded
    $\chi_j''(\Psi_t)|\nabla\Psi_t|^2\le0$.

    For fixed $j$, the support of $\chi_j'(\Psi_t)$ lies
    in a compact ball $B(o,R_{j,T})$, uniformly on $[0,T]$.
    For $\Psi=r^2$ this follows from $\chi_j'(u)=0$
    when $u\ge2j$; for $\Psi_t=e^{\alpha(t)r^2}$ it follows
    also from $\min_{[0,T]}\alpha>0$. Consequently
    $\chi_j(\Psi_t)$ is uniformly bounded and constant
    outside this ball, and
    \begin{equation*}
        \sup_{0\le t\le T}\sup_{x\in M}
        |\nabla\chi_j(\Psi_t(x))|\le C_{j,T}.
    \end{equation*}
    In particular, the quadratic variation of the
    martingale satisfies
    \begin{equation*}
        \langle M^{(j)}\rangle_T
        =2\int_0^T
            |\nabla\chi_j(\Psi_t(X_t))|^2\dd t
        \le2C_{j,T}^2T,
    \end{equation*}
    so $M^{(j)}$ is square-integrable and has mean zero.

    The drift $b_t$ is uniformly bounded on $B(o,R_{j,T})$
    for $t\in[0,T]$; indeed, the linear growth hypotheses
    give a bound in terms of $R_{j,T}$ and
    $\sup_{0\le t\le T}\int_M r(y)\dd\mu_t(y)$,
    which is finite by continuity in $\mathcal P_2(M)$.
    By the Laplacian comparison estimate
    \eqref{eq:comparison Laplacian} and the boundedness of
    $\Phi$ and its relevant derivatives on this ball,
    the regular drift has uniformly bounded positive part;
    \begin{equation*}
        \left[\chi_j'(\Psi_t)
            \bigl(\partial_t\Psi_t
                +(L_t\Psi_t)_{\mathrm{reg}}\bigr)\right]_+
        \le C'_{j,T}
    \end{equation*}
    for almost every $t\in[0,T]$ (on each interval where
    $\alpha$ is $C^1$). The cut locus has zero volume \cite[Theorem 10.34]{lee_introduction_2018} and the nondegenerate diffusion has absolutely
    continuous time marginals for positive times; thus the
    choice of the measurable extension of the regular part
    does not affect the time-integrated expectation.

    Finally, let $\tau_R$ be the first exit time of $X$
    from $B(o,R)$. Integrate the preceding inequality
    up to $t\wedge\tau_R$ and take expectations.
    The martingale has zero mean, while the singular
    cut-locus and concavity contributions are nonpositive.
    The positive regular drift is uniformly bounded and
    the negative finite-variation contributions are
    controlled by the stopped identity. Letting
    $R\to\infty$, nonexplosion gives bounded convergence
    for the cutoff terms, dominated convergence for the
    positive drift and Fatou's lemma for the negative part.
    Since $\operatorname{Law}(X_u)=\mu_u$, we obtain
    \eqref{eq:localized-generator}. No smooth initial
    density or compact support is required.
\end{proof}

\section{Covering estimate}\label{sec:appendix covering}
\begin{proof}[Proof of \cref{lem:packing}]
    Since the covering number $N(B(x,R), \delta)$ is bounded by the packing
    number $P(B(x,R), \delta)$, it is sufficient to estimate the latter.
    Let $\{x_i\}_{i=1}^{P(B(x,R), \delta)}$ be a maximal $\delta$-packing
    of $B(x,R)$; in particular, the balls $B(x_i,\delta/2)$ are pairwise disjoint
    and contained in $B(x,R+\delta/2)$ and therefore
    \begin{equation}\label{eq:packing-first}
        P(B(x,R),\delta)
        \le
        \max_{x_i} \frac{\m\left[B(x,R+\delta/2)\right]}{\m\left[B(x_i,\delta/2)\right]}.
    \end{equation}
    On the other hand, since $x_i\in B(x,R)$, we have
    \(
        B(x,R)\subseteq B(x_i,2R).
    \)
    Hence, by \cref{thm:Bishop Gromov},
    \begin{equation}\label{eq:packing-second}
        \frac{\m[B(x_i,\delta/2)]}{V_{K,n}(\delta/2)}
        \ge
        \frac{\m[B(x_i,2R)]}{V_{K,n}(2R)}
        \ge
        \frac{\m[B(x,R)]}{V_{K,n}(2R)},
    \end{equation}
    Inserting \eqref{eq:packing-second} into \eqref{eq:packing-first}, we obtain
    \begin{align*}
        P(B(x,R),\delta)
        &\le
        \frac{\m[B(x,R+\delta/2)]}{\m[B(x,R)]}
        \frac{V_{K,n}(2R)}{V_{K,n}(\delta/2)} \\
        &\le
        \frac{V_{K,n}(R+\delta/2)}{V_{K,n}(R)}
        \frac{V_{K,n}(2R)}{V_{K,n}(\delta/2)}.
        \numberthis \label{eq:packing-product}
    \end{align*}
    after another application of Bishop--Gromov.
    
    If $K=0$, then $V_{0,n}(r)$ is proportional to $r^n$, and therefore
    \begin{equation*}
        P(B(x,R),\delta)
        \le
        \left(1+\frac{\delta}{2R}\right)^n
        \left(\frac{4R}{\delta}\right)^n
        \le
        \left(\frac{5R}{\delta}\right)^n,
    \end{equation*}
    since $\delta\le R/2$. Otherwise, for $K<0$, using the inequality
    \[
        \frac{V_{K,n}(r_2)}{V_{K,n}(r_1)} = \frac{\int_0^{r_2} \sinh^{n-1}\left(\sqrt{\frac{|K|}{n-1}} s\right) ds}{\int_0^{r_1} \sinh^{n-1}\left(\sqrt{\frac{|K|}{n-1}} s\right) ds}
        \le 
        \left(\frac {r_2}{r_1}\right)^n
        \exp\!\left(\sqrt{|K|(n-1)}(r_2-r_1)\right)
    \]
    valid for all $0 < r_1 \le r_2$ in \eqref{eq:packing-product}, we obtain
    \begin{align*}
        P(B(x,R),\delta) \le
        \left(\frac{5R}{\delta}\right)^n
        \exp\!\left(2\sqrt{|K|(n-1)}\,R\right).
    \end{align*}
\end{proof}

\section*{Acknowledgments}
The author would like to thank Frank Redig for inspiring discussions and for proofreading an earlier version of this manuscript. The author is also grateful to Nathael Gozlan for suggesting the inclusion of concentration estimates around the expectation, and to François Bolley for his interest in this work. Finally, the author would like to thank Fabrice Baudoin for his explanations of Kendall's formula and Aarhus University for its support and hospitality during a period in which part of this work was carried out.

This work was supported by grant OCENW.M20.251
(\href{https://doi.org/10.61686/CJFWO82746}
{DOI: 10.61686/CJFWO82746})
from the Open Competition ENW research program, which is (partly) financed
by the Dutch Research Council (NWO) \href{https://www.nwo.nl/en/projects/cjfwo82746}{\includegraphics[height=\fontcharht\font`\B]{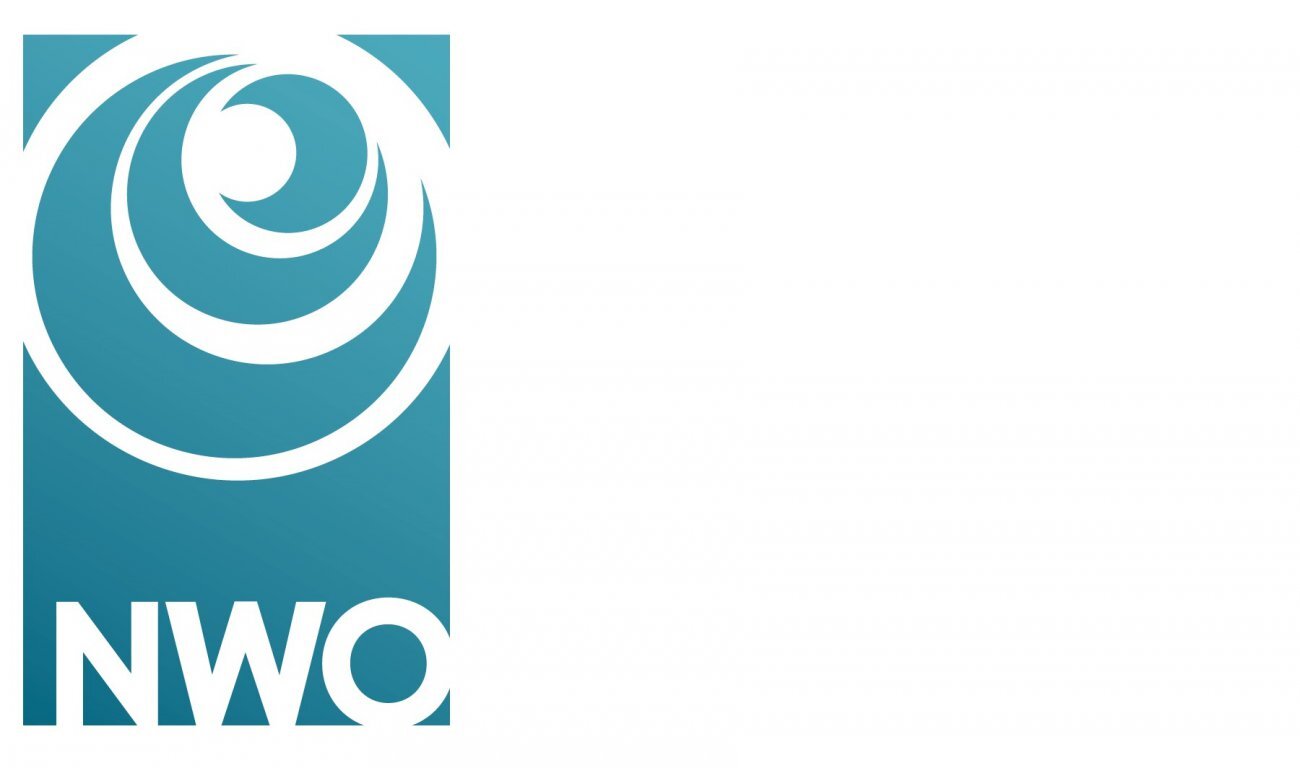}} \hspace{-10pt}.

\section*{Declaration of generative AI and AI-assisted technologies in the writing process}
The author is responsible for the architecture of the paper and proofs
of the main results. During the preparation of this manuscript, the
author used Codex and ChatGPT (GPT-6 Astra) to assist with editorial restructuring, proof
presentation, bibliographical research, symbolic
computations, and \LaTeX\ consistency. 
In particular, AI assistance was used to rigorously
justify the radial integration-by-parts computations
through the cutoff argument of \cref{lem:radial-cutoff}
(which can be avoided by working directly at the SDE
level, as explained in \cref{sec:preliminary results}),
and to establish propagation of higher-order
transport-cost inequalities discussed in
\cref{rem:aggregation-higher-p}.
The author reviewed and
edited all suggestions and takes full responsibility for the content
of the article.

\printbibliography

@article{boissard_mean_2014,
  author = {Boissard, Emmanuel and Le Gouic, Thibaut},
  title = {On the mean speed of convergence of empirical and occupation measures in {Wasserstein} distance},
  journal = {Annales de l'Institut Henri Poincar\'e, Probabilit\'es et Statistiques},
  year = {2014}, volume = {50}, number = {2}, pages = {539--563},
  doi = {10.1214/12-AIHP517}
}

@article{bolley_quantitative_2007,
  author = {Bolley, Fran\c{c}ois and Guillin, Arnaud and Villani, C\'edric},
  title = {Quantitative concentration inequalities for empirical measures on non-compact spaces},
  journal = {Probability Theory and Related Fields},
  year = {2007}, volume = {137}, number = {3--4}, pages = {541--593},
  doi = {10.1007/s00440-006-0004-7}
}

@article{fournier_rate_2015,
  author = {Fournier, Nicolas and Guillin, Arnaud},
  title = {On the rate of convergence in {Wasserstein} distance of the empirical measure},
  journal = {Probability Theory and Related Fields},
  year = {2015}, volume = {162}, number = {3--4}, pages = {707--738},
  doi = {10.1007/s00440-014-0583-7}
}

@article{weed_sharp_2019,
  author = {Weed, Jonathan and Bach, Francis},
  title = {Sharp asymptotic and finite-sample rates of convergence of empirical measures in {Wasserstein} distance},
  journal = {Bernoulli},
  year = {2019}, volume = {25}, number = {4A}, pages = {2620--2648},
  doi = {10.3150/18-BEJ1065}
}

@misc{dedecker_concentration_2026,
  author = {Dedecker, J\'er\^ome and Fischer, Aur\'elie and Michel, Bertrand},
  title = {Concentration of the empirical measure in {Wasserstein} distance: bounds involving the covering dimension},
  year = {2026}, eprint = {2601.11133}, archivePrefix = {arXiv},
  primaryClass = {math.PR}, url = {https://arxiv.org/abs/2601.11133},
  note = {Preprint}
}

@article{iacobelli_asymptotic_2016,
  author = {Iacobelli, Mikaela},
  title = {Asymptotic quantization for probability measures on {Riemannian} manifolds},
  journal = {ESAIM: Control, Optimisation and Calculus of Variations},
  year = {2016}, volume = {22}, number = {3}, pages = {770--785},
  doi = {10.1051/cocv/2015025}
}

@article{aydin_asymptotic_2025,
  author = {Ayd{\i}n, Ata Deniz and Iacobelli, Mikaela},
  title = {Asymptotic quantization of measures on {Riemannian} manifolds via covering growth estimates},
  journal = {Advances in Mathematics}, year = {2025}, volume = {474},
  pages = {110311}, doi = {10.1016/j.aim.2025.110311}
}

@article{gozlan_leonard_2007,
  author = {Gozlan, Nathael and L\'eonard, Christian},
  title = {A large deviation approach to some transportation cost inequalities},
  journal = {Probability Theory and Related Fields},
  year = {2007}, volume = {139}, number = {1--2}, pages = {235--283},
  doi = {10.1007/s00440-006-0045-y}
}

@article{kendall_radial_1987,
  author = {Kendall, Wilfrid S.},
  title = {The Radial Part of {Brownian} Motion on a Manifold: A Semimartingale Property},
  journal = {The Annals of Probability},
  volume = {15},
  number = {4},
  pages = {1491--1500},
  year = {1987},
  doi = {10.1214/aop/1176991988}
}

@book{hsu_stochastic_2002,
  author = {Hsu, Elton P.},
  title = {Stochastic Analysis on Manifolds},
  series = {Graduate Studies in Mathematics},
  volume = {38},
  publisher = {American Mathematical Society},
  address = {Providence, RI},
  year = {2002},
  doi = {10.1090/gsm/038}
}

@book{wang_analysis_2013,
  author = {Wang, Feng-Yu},
  title = {Analysis for Diffusion Processes on {Riemannian} Manifolds},
  series = {Advanced Series on Statistical Science \& Applied Probability},
  volume = {18},
  publisher = {World Scientific},
  address = {Singapore},
  year = {2013},
  doi = {10.1142/8737},
}

@article{zhang_2005,
  author = {Zhang, Xicheng},
  title = {Strong Solutions of {SDEs} with Singular Drift and {Sobolev} Diffusion Coefficients},
  journal = {Stochastic Processes and their Applications},
  volume = {115},
  number = {11},
  pages = {1805--1818},
  year = {2005},
  doi = {10.1016/j.spa.2005.06.003}
}

@article{mantegazza_mascellani_uraltsev_2014,
  author = {Mantegazza, Carlo and Mascellani, Giovanni and Uraltsev, Gennady},
  title = {On the Distributional {Hessian} of the Distance Function},
  journal = {Pacific Journal of Mathematics},
  volume = {270},
  number = {1},
  pages = {151--166},
  year = {2014},
  doi = {10.2140/pjm.2014.270.151},
}

@article{arnaudon_delmoral_2019,
  author = {Arnaudon, Marc and Del Moral, Pierre},
  title = {A Variational Approach to Nonlinear and Interacting Diffusions},
  journal = {Stochastic Analysis and Applications},
  volume = {37},
  number = {5},
  pages = {717--748},
  year = {2019},
  doi = {10.1080/07362994.2019.1609985}
}

@article{bolley_weighted_2005,
  author = {Bolley, Fran{\c c}ois and Villani, C{\'e}dric},
  title = {Weighted {Csisz{\'a}r--Kullback--Pinsker} Inequalities and Applications to Transportation Inequalities},
  journal = {Annales de la Facult{\'e} des Sciences de Toulouse. Math{\'e}matiques},
  series = {6},
  volume = {14},
  number = {3},
  pages = {331--352},
  year = {2005},
  doi = {10.5802/afst.1095},
}

@misc{chae_cho_jang_yang_2026,
  author = {Chae, Myeongju and Cho, Gunhee and Jang, Hyun Chul and Yang, Guang},
  title = {Propagation of Chaos on {Riemannian} Manifolds},
  year = {2026},
  eprint = {2609.12948},
  archivePrefix = {arXiv},
  primaryClass = {math.PR},
  note = {Preprint, submitted 11 September 2026},
  doi = {10.48550/arXiv.2609.12948},
  url = {https://arxiv.org/abs/2609.12948},
}

@book{lee_introduction_2018,
	address = {Cham},
	series = {Graduate {Texts} in {Mathematics}},
	title = {Introduction to {Riemannian} {Manifolds}},
	volume = {176},
	doi = {10.1007/978-3-319-91755-9},
	publisher = {Springer International Publishing},
	author = {Lee, John M.},
	year = {2018},
}

@article{ohta_measure_2007,
	title = {On the measure contraction property of metric measure spaces},
	volume = {82},
	doi = {10.4171/cmh/110},
	number = {4},
	journal = {Comment. Math. Helv.},
	author = {Ohta, Shin-ichi},
	year = {2007},
	pages = {805--828},
}

@article{wang_wasserstein_2022,
	title = {Wasserstein convergence rate for empirical measures on noncompact manifolds},
	volume = {144},
	doi = {10.1016/j.spa.2021.11.006},
	journal = {Stochastic Processes and their Applications},
	author = {Wang, Feng-Yu},
	year = {2022},
	pages = {271--287},
}

@article{boissard_simple_2011,
  author = {Boissard, Emmanuel},
  title = {Simple bounds for the convergence of empirical and occupation measures in 1-{Wasserstein} distance},
  journal = {Electronic Journal of Probability},
  year = {2011}, volume = {16}, number = {83}, pages = {2296--2333},
  doi = {10.1214/EJP.v16-958}
}

@article{fournier_convergence_2023,
  author = {Fournier, Nicolas},
  title = {Convergence of the empirical measure in expected {Wasserstein} distance: non-asymptotic explicit bounds in {$\mathbb{R}^d$}},
  journal = {ESAIM: Probability and Statistics},
  year = {2023}, volume = {27}, pages = {749--775},
  doi = {10.1051/ps/2023011}
}

@article{wang_precise_2021,
  author = {Wang, Feng-Yu},
  title = {Precise limit in {Wasserstein} distance for conditional empirical measures of {Dirichlet} diffusion processes},
  journal = {Journal of Functional Analysis},
  year = {2021}, volume = {280}, number = {11}, pages = {108998},
  doi = {10.1016/j.jfa.2021.108998}
}

@article{arnaudon_harnack_2006,
	title = {Harnack inequality and heat kernel estimates on manifolds with curvature unbounded below},
	volume = {130},
	copyright = {https://www.elsevier.com/tdm/userlicense/1.0/},
	issn = {00074497},
	url = {https://linkinghub.elsevier.com/retrieve/pii/S0007449705001004},
	doi = {10.1016/j.bulsci.2005.10.001},
	language = {en},
	number = {3},
	urldate = {2025-03-27},
	journal = {Bulletin des Sciences Mathématiques},
	author = {Arnaudon, Marc and Thalmaier, Anton and Wang, Feng-Yu},
	month = apr,
	year = {2006},
	pages = {223--233},
}

@book{villani_optimal_2009,
	address = {Berlin, Heidelberg},
	series = {Grundlehren der mathematischen {Wissenschaften}},
	title = {Optimal {Transport}},
	volume = {338},
	copyright = {http://www.springer.com/tdm},

	url = {http://link.springer.com/10.1007/978-3-540-71050-9},
	doi = {10.1007/978-3-540-71050-9},
	urldate = {2026-02-12},
	publisher = {Springer Berlin Heidelberg},
	author = {Villani, Cédric},
	editor = {Berger, M. and Eckmann, B. and De La Harpe, P. and Hirzebruch, F. and Hitchin, N. and Hörmander, L. and Kupiainen, A. and Lebeau, G. and Ratner, M. and Serre, D. and Sinai, Ya. G. and Sloane, N. J. A. and Vershik, A. M. and Waldschmidt, M.},
	year = {2009},
}

@article{wang_wasserstein_2021,
  author={Wang, Feng-Yu and Wu, Bingyao},
  title={Wasserstein convergence for empirical measures of subordinated diffusions on {Riemannian} manifolds},
  journal={Potential Analysis}, volume={59}, number={3}, pages={933--954},
  year={2023}, doi={10.1007/s11118-022-09989-6}
}

@article{shalova_schlichting_2026,
  author = {Shalova, Anna and Schlichting, Andr{\'e}},
  title = {Solutions of stationary {McKean--Vlasov} equation on a high-dimensional sphere and other {Riemannian} manifolds},
  journal = {Advances in Nonlinear Analysis},
  volume = {15}, number = {1}, pages = {20250141}, year = {2026},
  doi = {10.1515/anona-2025-0141}
}

@misc{peletier_shalova_2025,
  author = {Peletier, Mark A. and Shalova, Anna},
  title = {Nonlinear diffusion limit of non-local interactions on a sphere},
  year = {2025}, eprint = {2512.03185}, archivePrefix = {arXiv},
  primaryClass = {math.AP}, doi = {10.48550/arXiv.2512.03185},
  note = {Preprint}
}

@article{fetecau_park_vaidya_2026,
  author = {Fetecau, Razvan C. and Park, Hansol and Vaidya, Vishnu},
  title = {Global energy minimizers for a diffusion-aggregation model on sphere},
  journal = {Nonlinearity}, year = {2026}, volume = {39}, number = {5},
  pages = {055011}, doi = {10.1088/1361-6544/ae6961}
}

@misc{dimarco_quantization_2026,
  author = {Di Marco, Marco and Iacobelli, Mikaela},
  title = {Quantization of measures in {Carnot} groups},
  year = {2026},
  eprint = {2606.09373},
  archivePrefix = {arXiv},
  primaryClass = {math.MG},
  url = {https://arxiv.org/abs/2606.09373},
  note = {Preprint}
}

@incollection{cattiaux_guillin_2014,
  author    = {Cattiaux, Patrick and Guillin, Arnaud},
  title     = {Semi Log-Concave {Markov} Diffusions},
  booktitle = {Séminaire de Probabilités {XLVI}},
  series    = {Lecture Notes in Mathematics},
  volume    = {2123},
  pages     = {231--292},
  publisher = {Springer},
  year      = {2014},
  doi       = {10.1007/978-3-319-11970-0_9}
}
\end{document}